\documentclass[aop,seceqn,citesort,MSNbibl,noautosecdot,dvips]{arximspdf}
\usepackage{mathrsfs}
\usepackage{graphicx}

\doi{10.1214/10-AOP582}
\volume{39}
\issue{5}
\pubyear{2011}
\firstpage{1702}
\lastpage{1767}

\makeatletter

\newcommand{\iint}{\int\!\!\int}

\newtheorem{theorem}{Theorem}[section]
\newtheorem{corollary}[theorem]{Corollary}
\newtheorem{lemma}[theorem]{Lemma}
\newtheorem{proposition}[theorem]{Proposition}

\newproclaim{definition}[theorem]{Definition}
\newproclaim{remark}[theorem]{Remark}
\newproclaim{Example}{Example}

\newcommand{\Prob}{{\mathbb P}}
\newcommand{\R}{\mathbb{R}}
\newcommand{\Q}{\mathbb{Q}}
\newcommand{\C}{\mathbb{C}}
\newcommand{\Z}{\mathbb{Z}}
\newcommand{\N}{\mathbb{N}}
\newcommand{\T}{\mathbb T}

\newcommand{\Piv}{\mathscr{P}}
\newcommand{\dist}{\operatorname{dist}}
\newcommand{\SLE}{\operatorname{SLE}}
\newcommand{\Var}{\operatorname{Var}}
\newcommand{\Cov}{\operatorname{Cov}}
\newcommand{\eps}{\varepsilon}
\newcommand{\PP}{\Prob}

\newcommand{\p}{{\partial}}
\newcommand{\E}{{\mathbb E}}

\newcommand{\sign}{\operatorname{sign}}

\newcommand{\Spec}{\mathscr{S}}

\newcommand{\Inf}{\mathbf{I}}
\newcommand{\InfV}{\operatorname{\mathbf{Inf}}}
\newcommand{\II}{\mathbf{H}}
\newcommand{\Circ}{\operatorname{Circ}}
\newcommand{\Exc}{\mathcal{E}}
\newcommand{\MAJ}{\Phi}
\newcommand{\reveal}{\bolds{\delta}}

\makeatother

\begin{document}
\begin{frontmatter}

\title{Oded Schramm's contributions to noise sensitivity}
\runtitle{Oded Schramm's contributions to noise sensitivity}

\begin{aug}
\author[A]{\fnms{Christophe} \snm{Garban}\corref{}\thanksref{t2}\ead[label=e1]{christophe.garban@umpa.ens-lyon.fr}}
\runauthor{C. Garban}
\affiliation{ENS Lyon, CNRS}
\dedicated{This paper is dedicated to the memory of Oded Schramm.\break
I feel very fortunate to have known him.}
\address[A]{UMPA\\
Ecole Normale Sup\'erieure de Lyon\\
46, all\'ee d'Italie\\
69364 Lyon Cedex 07 \\
France\\
\printead{e1}} 
\end{aug}

\thankstext{t2}{Supported in part by ANR-06-BLAN-00058 and the
Ostrowski foundation.}

\received{\smonth{7} \syear{2010}}

%
\begin{abstract}
We survey in this paper the main contributions of Oded Schramm related
to \textit{noise sensitivity}. We will describe in particular his various
works which focused on the ``spectral analysis'' of critical
percolation (and more generally of Boolean functions), his work on the
shape-fluctuations of first passage percolation and finally his
contributions to the model of dynamical percolation.
\end{abstract}

%

\end{frontmatter}

A sentence which summarizes well Oded's work on noise sensitivity is
the following quote from
Jean Bourgain.

\textit{There is a general philosophy which claims that if $f$ defines a
property of `high complexity,'
then $\mathrm{Sup} \hat f$, the support of the Fourier Transform, has
to be `spread out.'}

Through his work on models coming from statistical physics (in
particular percolation), Oded Schramm was
often confronted with such functions of ``\textit{high complexity}.'' For
example, in percolation, any large-scale connectivity property
can be encoded by a Boolean function of the ``inputs'' (edges or sites).
At criticality, these large-scale connectivity functions turn out to
be of ``high frequency'' which gives deep information on the underlying model.
As we will see along this survey, Oded Schramm developed over the last
decade 
highly original and deep ideas to understand the ``complexity'' of
Boolean functions.

We will essentially follow the chronology of his contributions in the field;
it is quite striking that three distinct periods emerge from Oded's
work and they will constitute Sections \ref{s.hyper}, \ref{s.algo} and
\ref{s.geo}
of this survey, corresponding, respectively, to the papers
\cite{MR2001m60016,SchrammSteif,GPS}.

We have chosen to present numerous sketches of proof, since we believe
that the elegance of Oded's mathematics is best illustrated
through his proofs, which usually combined imagination, simplicity and
powerfulness.
This choice is the reason for the length of the present survey.

\vspace*{15pt}
\tableofcontents[level=2]

\section{\texorpdfstring{Introduction.}{Introduction}}\label{sec1}
In this Introduction, we will start by describing the scientific
context in the late nineties
which lead Oded Schramm to work on the sensitivity of Boolean
functions. We will then motivate and define
the notion of \textit{noise sensitivity} and finally we will review the
main contributions of
Oded that will concern us throughout this survey.

\subsection{\texorpdfstring{Historical context.}{Historical context}}

When Oded started working on Boolean functions (with the idea to use
them in statistical mechanics), there was already
important literature in computer science dedicated to the properties of
Boolean functions.

Here is an example of a related problem which was solved before Oded
came into the field:
in \cite{BenorLinial}, Ben-Or and Linial conjectured that if $f$ is
any Boolean function on $n$ variables (i.e., $f \dvtx \{0,1\}^n \to\{0,1\}$),
taking the value~1 for half of the configurations of the hypercube $\{
0,1\}^n$; then there exists some deterministic set $S=S_f\subset\{
1,\ldots, n\}$
of less than $c \frac{n} {\log n}$ variables, such that $f$ remains
undetermined as long as the variables in $S$ are not assigned
(the constant $c$ being a universal constant). This means that for any
such Boolean function $f$, there should always exist a set of small
size which is ``pivotal''
for the outcome.

This conjecture was proved in \cite{KKL}. Besides the proof of the
conjecture,
what is most striking about this paper (and which will concern us
throughout the rest of this
survey) is that, for the first time, techniques brought from harmonic
analysis were used in \cite{KKL} to study properties of Boolean functions.
At the time, the authors wrote, ``These new connections with harmonic
analysis are very promising.''

Indeed, as they anticipated, the main technique they borrowed from
harmonic analysis, namely \textit{hypercontractivity}, was later used in
many subsequent
works. In particular, as we will see later (Section \ref{s.hyper}),
hypercontractivity was one of the main ingredients in the landmark
paper on
\textit{noise sensitivity} written by Benjamini, Kalai and Schramm
\cite{MR2001m60016}.

Before going into \cite{MR2001m60016} (which introduced the concept of
\textit{noise sensitivity}), let us mention some of the related works in this
field which
appeared in the period from \cite{KKL} to \cite{MR2001m60016} and
which made
good use of hypercontractivity. We distinguish several directions of research.
\begin{itemize}
\item First of all, the results of \cite{KKL} have been extended to
more general cases: nonbalanced Boolean functions, other measures than
the uniform measure on the hypercube $\{0,1\}^n$ and finally,
generalizations to ``Boolean'' functions of continuous dependence $f
\dvtx
[0,1]^n \to\{0,1\}$. See \cite{MR1194785} as well as
\cite{MR1303654}. Note that both of these papers rely on
hypercontractive estimates.

\item Based on these generalizations, Friedgut and Kalai studied in
\cite{MR1371123} the phenomenon of ``sharp thresholds.''
Roughly speaking, they proved that any monotone ``graph property'' for
the so-called random graphs $G(n,p), 0\le p \le1$ satisfies a
\textit{sharp threshold} as the number of vertices $n$ goes to infinity (see
\cite{MR1371123} for a more precise statement).
In other words, they show that any monotone graph event $\mathcal{A}$
appears ``suddenly'' while increasing the edge intensity $p$ [whence a
``cut-off'' like shape for the function $p\mapsto\PP_{n,p}(\mathcal{A})$].
In some sense their work is thus intimately related to the subject of
this survey, since many examples of such ``graph properties'' concern
connectedness,
size of clusters and so on.

The study of \textit{sharp thresholds} began in percolation theory with
the seminal work of Russo \cite{MR0488383,MR671248}, where he
introduced the idea of sharp thresholds
in order to give an alternate proof of Kesten's result $p_c(\Z
^2)=\frac
1 2$. See also the paper by Margulis \cite{MR0472604}.

\item Finally, Talagrand made several remarkable contributions over
this period which highly influenced
Oded and his coauthors (as we will see in particular through Section
\ref{s.hyper}). To name a few of these: an important result
of~\cite{MR2001m60016} (Theorem \ref{th.ns} in this survey) was inspired by
\cite{MR1401897}; the study of fluctuations in first passage
percolation in \cite{MR2016607} (see Section \ref{s.hyper}) was
influenced by a result from \cite{MR1303654} (this paper by
Talagrand was already mentioned above since it somewhat overlapped with
\cite{MR1194785,MR1371123}). More generally the questions addressed
along this survey are related to the \textit{concentration of measure
phenomenon} which has been deeply understood by Talagrand (see
\cite{MR1361756}).
\end{itemize}

\subsection{\texorpdfstring{Concept of noise sensitivity.}{Concept of noise sensitivity}}

It is now time to motivate and then define what is \textit{noise
sensitivity}. This concept was introduced in the seminal work
\cite{MR2001m60016} whose content
will be described in Section \ref{s.hyper}. As one can see from its
title, \textit{noise sensitivity of Boolean functions and applications to
Percolation}, the authors
introduced this notion having in mind applications in percolation.
Before explaining these motivations in the next subsection, let 
us consider the following simple situation.

Assume we have some initial data $\omega=(x_1,\ldots,x_n)\in\{0,1\}
^n$, and we are interested in some functional of this data that we may represent
by a (real-valued) function $f \dvtx \{0,1\}^n \to\R$. Often, the
functional $f$ will be the indicator function of an event $A\subset\{
0,1\}^n$; in other words $f$
will be \textit{Boolean}. Now imagine that there are some errors in the
transmission of this data $\omega$ and that one receives
only the sightly perturbed data $\omega^{\eps}$; one would like that
the quantity we are interested in, that is, $f(\omega)$, is not ``too
far'' from what we actually
receive, that is, $f(\omega^\eps)$. (A similar situation would be: we
are interested in some \textit{real} data $\omega$ but there are some
inevitable errors in collecting
this data and one ends up with $\omega^\eps$.) To take a typical
example, if $f \dvtx \{0,1\}^n\to\{0,1\}$ represents a \textit{voting
scheme}, for example, majority,
then it is natural to wonder how robust $f$ is with respect to errors.

At that point, one would like to know how to properly model the
perturbed data $\omega^\eps$. The correct answer depends on the
probability law which governs
our ``random'' data $\omega$. In the rest of this survey, our data will
always be sampled according to the uniform measure; hence it will be
sufficient for us
to assume that $\omega\sim\PP=\PP^n=(\frac1 2 \delta_0 + \frac1 2
\delta_1)^{\otimes n}$, the uniform measure on $\{0,1\}^n$. Other
natural measures may be considered instead, but we will stick to this
simpler case. Therefore a natural way to model the perturbed
data~$\omega^\eps$ is to assume that each variable in $\omega=(x_1,\ldots
,x_n)$ is \textit{resampled} independently with small probability $\eps$.
Equivalently, if $\omega=(x_1,\ldots,x_n)$,
then~$\omega^\eps$ will correspond to the random configuration
$(y_1,\ldots,y_n)$, where independently for each $i\in\{1,\ldots,n\}$,
with probability $1-\eps$, $y_i:=x_i$ and with probability $\eps$,
$y_i$ is sampled according to a Bernoulli ($1/2$).
It is clear that such a ``noising'' procedure preserves the uniform
measure on $\{0,1\}^n$.

In computer science, one is naturally interested in the \textit{noise
stability} of~$f$ which, if $f$ is Boolean (i.e., with values in $\{
0,1\}$),
can be measured by the quantity ${\mathbb P}[{f(\omega)\neq
f(\omega^\eps)}]$,
where as above $\PP=\PP^n$ denotes the uniform measure on $\{0,1\}^n$
[in fact there is a slight abuse of notation here, since~$\PP$ samples
the coupling $(\omega, \omega^\eps)$; hence there is also some extra
randomness needed for the randomization procedure]. In
\cite{MajorityStablest}, it is shown using Fourier analysis that in some
sense, the most stable Boolean
function on $\{0,1\}^n$ is the majority function (under reasonable
assumptions which exclude dictatorship and so~on).

If a functional $f$ happens to be ``noise stable,''
this means that very little information is lost on the outcome
$f(\omega
)$ knowing the ``biased'' information~$\omega^\eps$.
Throughout this survey, we will be mainly interested in the opposite
property, namely \textit{noise sensitivity}. We will give precise
definitions later,
but roughly speaking, $f$ will be called ``noise sensitive'' if almost
ALL information is lost on the outcome $f(\omega)$ knowing the biased
information $\omega^\eps$.
This complete loss of information can be measured as follows: $f$ will
be ``noise sensitive'' if \mbox{$\Var[ {\mathbb E}[{f(\omega)
\mid \omega^\eps}]] =o(1)$}. It turns out
that if $\Var(f)=O(1)$, it is equivalent to consider the correlation
between $f(\omega)$ and $f(\omega^\eps)$ and equivalently, $f$ will be
called ``noise sensitive''
if $\Cov(f(\omega), f(\omega^\eps)) = {\mathbb E}[{f(\omega
)f(\omega^\eps)}] - {\mathbb E}[{f}]^2 = o(1)$.
\begin{remark}
Let us point out that \textit{noise stability} and \textit{noise sensitivity}
are two extreme situations.
Imagine our functional $f$ can be decomposed as $f= f_{\mathrm{stable}}
+ f_{\mathrm{sens}}$,
then after transmission of the data, we will keep the information on
$f_{\mathrm{stable}}$ but lose the information on the other component.
As such
$f$ will be neither stable nor sensitive. 
\end{remark}

We will end by an example of a noise sensitive functional in the
context of percolation. Let us consider a percolation configuration
$\omega_n$
on the rescaled lattice $\frac1 n \Z^2$ in the window $[0,1]^2$ (see
Section \ref{ss.perco} for background and references on percolation)
at criticality $p_c=\frac1 2$
(hence $\omega_n$ is sampled according to the uniform measure on $\{
0,1\}^{E_n}$, where $E_n$ is the set of edges of $\frac1 n \Z^2 \cap
[0,1]^2$).
In Figure \ref{f.perconoise}, we represented a percolation
configuration $\omega_{90}$ and its noised configuration $\omega
_{90}^\eps$ with $\eps=\frac1 {10}$.
In each picture, the three largest clusters are colored in red, blue
and green. As is suggested from this small-scale picture, the functional
giving the size of the largest cluster turns out to be asymptotically
``noise sensitive''; that is, for fixed $\eps>0$ when $n\to\infty$
the entire
information about the largest cluster is lost from $\omega_n$ to
$\omega
_n^\eps$ (even though most of the ``microscopic'' information is preserved).

For nice applications of the concept of noise sensitivity in the
context of computer science, see \cite{MR2208732} and
\cite{ODonnellThesis}.

\subsection{\texorpdfstring{Motivations from statistical physics.}{Motivations from statistical physics}}
Beside the obvious motivations in computer science, there were several
motivations
coming from statistical physics (mostly from percolation)
which lead Benjamini, Kalai and Schramm to introduce and study \textit{noise sensitivity} of Boolean functions.
We wish to list some of these in this subsection.

\subsubsection*{\texorpdfstring{Dynamical percolation.}{Dynamical
percolation}}
In ``real life,'' models issued from statistical phy\-sics
undergo \textit{thermal agitation}; consequently, their state
evolves in time. For example, in the case of the Ising model, the
natural dynamics associated to thermal agitation is
the so-called \textit{Glauber dynamics}.

In the context of percolation, a natural dynamics modeling this
\textit{thermal agitation} has been introduced in \cite{MR1465800}
under the name of \textit{dynamical percolation} (it was also invented
independently by Itai Benjamini). This model is defined in a very
simple way and we will describe it in detail in Section \ref{s.dp} (see
\cite{SurveySteif} for a very nice survey).

For percolation in dimension two, it is known that at
\textit{criticality}, there is almost surely no infinite cluster.
Nevertheless, the following open question was asked back in
\cite{MR1465800}: if one lets the critical percolation configuration
evolve according to this dynamical percolation process, is it the case
that there will exist \textit{exceptional times} where an infinite
cluster will appear? As we will later see, such exceptional times
indeed exist. This reflects the fact that the \textit{dynamical}
behavior of percolation is very different from its \textit{static}
behavior.


Dynamical percolation is intimately related to our concept of noise
sensitivity since if $(\omega_t)$ denotes the trajectory of a dynamical
percolation process, then the configuration at time $t$ is exactly a
``noised'' configuration of the initial configuration $\omega_0$
[$\omega_t \equiv \omega_0^\eps$ with an explicit correspondence
between $t$ and $\eps$ (see Section~\ref{s.dp})].

%
%
\begin{figure}[t]

\includegraphics{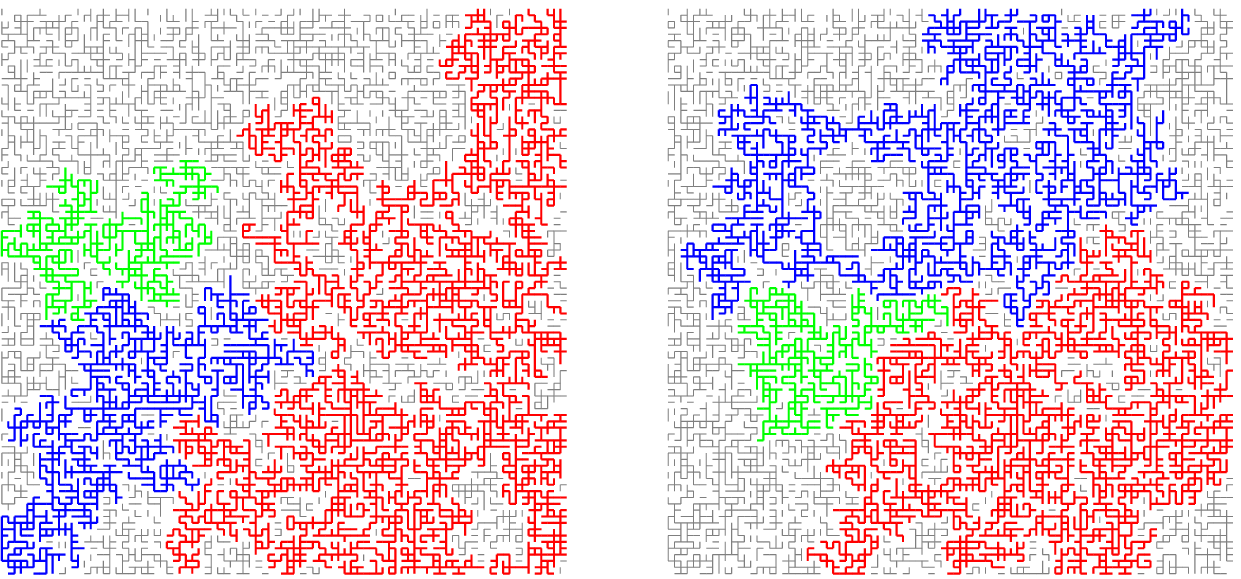}

\caption{}\label{f.perconoise}
\end{figure}

As is suggested by Figure \ref{f.perconoise}, large clusters of
$(\omega _t)_{t\geq0}$ ``move'' (or change) very quickly as the time
$t$ goes on. This \textit{rapid mixing} of the large scale properties
of $\omega_t$ is the reason for the appearance of ``exceptional''
infinite clusters along the dynamics. Hence, as we will see in Section
\ref{s.dp}, the above question addressed in \cite{MR1465800} needs an
analysis of the \textit{noise sensitivity} of large clusters. In cite
\cite{MR2001m60016}, the first results in this direction are shown
(they prove that in some sense any large-scale connectivity property is
``noise sensitive''). However, their control on the sensitivity of
large scale properties (later we will rather say the ``frequency'' of
large-scale properties) was not sharp enough to imply that exceptional
times do exist. The question was finally answered in the case of the
triangular lattice in \cite{SchrammSteif} thanks to a better
understanding of the sensitivity of percolation. The case of $\Z^2$ was
carried out in \cite{GPS}. This conjecture from \cite{MR1465800} on the
existence of \textit{exceptional times} will be one of our driving
motivations throughout this survey.

\subsubsection*{\texorpdfstring{Conformal invariance of
percolation.}{Conformal invariance of percolation}}
In the nineties, an important conjecture popularized by
Langlands, Pouliot and Saint-Aubin in \cite{MR1230963}
stated that critical percolation should satisfy some (asymptotic)
\textit{conformal invariance} (in the same way as Brownian motion does).
Using conformal field theory, Cardy made several important predictions
based on this conformal invariance principle (see \cite{MR92m82048}).

\textit{Conformal invariance} of percolation was probably the main
conjecture on planar percolation at that time and people working
in this field, including Oded, were working actively on this question
until Smirnov solved it on the triangular lattice in
\cite{MR1851632} (note that a major motivation
of the introduction of the $\mathrm{SLE}_\kappa$ processes by Oded in
\cite{MR1776084} was this conjecture).

At the time of \cite{MR2001m60016}, conformal invariance of
percolation was not
yet proved (it still remains open on $\Z^2$),
and the $\mathrm{SLEs}$ from \cite{MR1776084} were ``in construction,''
so the route toward conformal invariance was still vague.
One nice idea from \cite{MR2001m60016} in this direction was to randomly
perturb the lattice itself and then claim that the crossing
probabilities are almost not affected
by the random pertubation of the lattice. This setup is easily seen to
be equivalent to proving noise sensitivity. What was difficult and
remained to be done
was to show that one could use such random perturbations to go from one
``quad'' to another conformally equivalent ``quad'' (see \cite{MR2001m60016}
for more
details). Note that at the same period, a different direction to attack
conformal invariance was developed by Benjamini and Schramm in
\cite{MR1646475} where they proved a certain kind of
conformal invariance for Vorono\"i percolation.

\subsubsection*{\texorpdfstring{Tsirelson\textup{'}s \textup{``}black noise.\textup{''}}{Tsirelson's ``black noise''}}
In \cite{MR1606855}, Tsirelson and Vershik constructed
sig\-ma-fields which are not ``produced'' by white noise (they called
these sigma-fields
``nonlinearizable'' at that time). Tsirelson then realized that a~good
way to characterize
a ``Brownian filtration'' was to test its stability versus
perturbations. Systems intrinsically unstable led him to the notion of
\textit{black noise}
(see \cite{MR2079671} for more details). This instability
corresponds exactly to our notion of being ``noise sensitive,'' and
after \cite{MR2001m60016} appeared,
Tsirelson used the concept of \textit{noise sensitivity} to describe his
theory of black noises. Finally, black noises are related to
percolation, since according to
Tsirelson himself, percolation would provide the most important example
of a (two-dimensional) black noise. In some sense \cite{MR2001m60016} shows
that if
percolation can be seen asymptotically as a \textit{noise} (i.e., a
continuous sigma-field which ``factorizes'' (see
\cite{MR2079671})), then this noise has to be
\textit{black}, that is, all its observables or functionals are noise
sensitive. The remaining step (proving that percolation asymptotically
factorizes as a \textit{noise})
was proved recently by Schramm and Smirnov~\cite{SSblacknoise}.

\subsubsection*{\texorpdfstring{Anomalous fluctuations.}{Anomalous fluctuations}}
In a different direction, we will see that noise sensitivity
is related to the study of random shapes whose fluctuations are
much smaller than ``Gaussian'' fluctuations (this is what we mean by
``anomalous fluctuations'').
Very roughly speaking, if a random metric space is highly sensitive to
noise (in the sense that its geodesics are ``noise sensitive''),
then it induces a lot of independence within the system itself and the
metric properties of the system decorrelate fast (in space or
under perturbations). This instability implies in general very small
fluctuations for the macroscopic metric properties (like the shape of
large balls and so on).
In Section \ref{s.hyper}, we will give an example of a random metric
space, first passage percolation, whose fluctuations can be
analyzed with the same techniques as the ones used in \cite{MR2001m60016}.


\subsection{\texorpdfstring{Precise definition of noise
sensitivity.}{Precise definition of noise sensitivity}}
Let us now fix some notations and definitions that will be used
throughout the rest of the article,
especially the definition of \textit{noise sensitivity} which was only
informal so far.

First of all, henceforth, it will be more convenient to work with the
hypercube $\{-1,1\}^n$ rather than with $\{0,1\}^n$.
Let us then call $\Omega_n:= \{-1,1\}^n$. The advantage of this choice
is that the \textit{characters} of $\{-1,1\}^n$ have a more simple form.

In the remainder of the text a \textit{Boolean function} will denote a
function from~$\Omega_n$ into $\{0, 1\}$
(except in Section \ref{s.geo}, where when made precise it could also
be into $\{-1,1\}$) and as argued above, $\Omega_n$ will be endowed
with the uniform probability measure $\PP=\PP^n$ on $\Omega_n$.
Through this survey, we will sometimes extend the study to the larger
class of real-valued functions
from~$\Omega_n$ into $\R$. Some of the results will hold for this
larger class, but the Boolean hypothesis will be crucial
at several different steps.

In the above informal discussion, ``noise sensitivity'' of a function
$f \dvtx \Omega_n \to\{0,1\}$ corresponded
to $\Var[{\mathbb E}[{f(\omega)\mid \omega^\eps
}]]$ or $\Cov
(f(\omega
),f(\omega^\eps))$ being ``small.''
To be more specific, \textit{noise sensitivity} is defined in \cite{MR2001m60016}
as an asymptotic property.
\begin{definition}[\cite{MR2001m60016}]\label{d.NS}
Let $(m_n)_{n\geq0}$ be an increasing sequence in $\N$. A~sequence of
Boolean functions $f_{n} \dvtx \{-1,1\}^{m_n} \to\{0,1\}$
is said to be (asymptotically) \textit{noise sensitive} if for any level
of noise $\eps>0$,
%
%
\begin{equation}\label{e.NS}
\Cov(f_n(\omega), f_n(\omega^\eps)) = {\mathbb E}
[{f_n(\omega) f_n(\omega^\eps)}] -{\mathbb E}
[{f_n}]^2 \mathop{\longrightarrow}_{n\to\infty} 0 .
\end{equation}
\end{definition}

In \cite{MR2001m60016}, the asymptotic condition was rather that
\[
\Var[{\mathbb E}[{f_n(\omega) \mid \omega^\eps}
] ] \- \mathop
{\longrightarrow}_{n\to\infty} 0 ,
\]
but as we mentioned above, the definitions are easily seen to be
equivalent (using the Fourier expansions of $f_n$).
\begin{remark}
One can extend in the same fashion this definition to the class of
(real-valued) functions
$f_n \dvtx \Omega_{m_n} \to\R$ of bounded variance (bounded variance is
needed to guarantee the equivalence of the two above criteria).
\end{remark}
\begin{remark}
Note that if $\Var(f_n)$ goes to zero as $n\to\infty$, then $(f_n)$
will automatically satisfy the condition (\ref{e.NS}),
hence our definition of noise sensitivity is meaningful only for
nondegenerate asymptotic events.
\end{remark}

The opposite notion of \textit{noise stability} is defined in \cite{MR2001m60016}
as follows:
\begin{definition}[\cite{MR2001m60016}]
Let $(m_n)_{n\geq0}$ be an increasing sequence in $\N$. A~sequence of
Boolean functions $f_n \dvtx \{-1,1\}^n \to\{0,1\}$
is said to be (asymptotically) \textit{noise stable} if
\[
\sup_{n\geq0} {\mathbb P}[{f_n(\omega) \neq f_n(\omega^\eps
)}] \mathop
{\longrightarrow}_{\eps\to0} 0.
\]
\end{definition}

\subsection{\texorpdfstring{Structure of the paper.}{Structure of the paper}}
In Section \ref{s.prel}, we will review some necessary background:
Fourier analysis of Boolean functions, the notion of \textit{influence}
and some facts about percolation.
Then three Sections \ref{s.hyper}, \ref{s.algo} and \ref{s.geo}, form
the core of this survey. They present three different approaches, each
of them
enabling to localize with more or less accuracy the ``frequency
domain'' of percolation.

The approach presented in Section \ref{s.hyper} is based on a
technique, \textit{hypercontractivity}, brought from harmonic analysis.
Following \cite{MR2001m60016} and \cite{MR2016607} we apply this
technique to the
sensitivity of percolation as well as to the study of shape fluctuations
in first passage percolation.
In Section \ref{s.algo}, we describe an approach developed by Schramm
and Steif in \cite{SchrammSteif} based on the analysis of randomized
algorithms.
Section \ref{s.geo}, following \cite{GPS}, presents an approach which
considers the ``frequencies'' of percolation as
random sets in the plane; the purpose is then to study the law of these
``frequencies'' and to prove that they behave in some ways like random
Cantor sets.

Finally, in Section \ref{s.dp} we present applications of the detailed
information provided by the last two approaches (\cite{GPS}
and \cite{SchrammSteif})
to the model of dynamical percolation.


The contributions that we have chosen to present reveal personal tastes
of the author. Also, the focus here is mainly on the
applications in statistical mechanics and particularly percolation.
Nevertheless we will try as much as possible, along this survey, to
point toward other contributions Oded made
close to this theme (such as
\cite{MR99i60173,MR2309980,MR2181623,DecisionTrees}). 
See also the very nice survey by Oded \cite{MR2334202}.

\section{\texorpdfstring{Background.}{Background}}\label{s.prel}

In this section, we will give some preliminaries on Boo\-lean functions
which will be used throughout the rest of the present survey.
We will start with an analog of Fourier series for Boolean functions;
then we will define the \textit{influence} of a variable, a notion which will
be particularly relevant in the remainder of the text; we will end the
preliminaries section with a brief review on percolation since most of Oded's
work in noise sensitivity was motivated by applications in percolation theory.

\subsection{\texorpdfstring{Fourier analysis of Boolean functions.}{Fourier analysis of Boolean functions}}\label{ss.fourier}
In this survey, recall that we consider our Boolean functions as
functions from the hypercube $\Omega_n:= \{-1,1\}^n$
into $\{0,1\}$, and $\Omega_n$ will be endowed with the uniform measure
$\PP=\PP^n = (\frac1 2 \delta_{-1} + \frac1 2 \delta_1)^{\otimes n}$.
\begin{remark}
In greater generality, one could consider other natural measures like
$\PP_p = \PP^n_p = ((1-p)\delta_{-1} + p\delta_1)^{\otimes n}$;
these measures are relevant, for example, in the study of sharp
thresholds (where one increases the ``level''~$p$). In the remainder of
the text, it will be
sufficient for us to restrict to the case of the uniform measure $\PP
=\PP_{1/2}$ on $\Omega_n$.
\end{remark}

In order to apply Fourier analysis, the natural setup is to enlarge our
discrete space of Boolean functions and to consider instead the larger space
$L^2(\{-1,1\}^n)$ of real-valued functions on $\Omega_n$ endowed with
the inner product
\begin{eqnarray*}
\langle f,g\rangle:\!&=& \sum_{x_1,\ldots, x_n} 2^{-n} f(x_1,\ldots
,x_n)g(x_1,\ldots,x_n) \\
&=&{\mathbb E}[{f g}] \qquad\mbox{for all } f,g \in
L^2(\Omega_n) ,
\end{eqnarray*}
where $\E$ denotes the expectation with respect to the uniform measure
$\PP$ on~$\Omega_n$.

For any subset $S\subset\{1,2\ldots,n\}$, let $\chi_S$ be the function
on $\{-1,1\}^n$ defined for any $x=(x_1,\ldots,x_n)$ by
%
%
\begin{equation}
\chi_S(x): = \prod_{i\in S} x_i .
\end{equation}
It is straightforward to see that this family of $2^n$ functions forms
an orthonormal
basis of $L^2(\{-1,1\}^n)$. Thus, any function $f$ on $\Omega_n$ (and a
fortiori any Boolean function $f$) can be decomposed as
\[
f = \sum_{S\subset\{1,\ldots,n\} } \hat f(S) \chi_S,
\]
where $\hat f(S)$ are the Fourier coefficients of $f$. They are
sometimes called the \textit{Fourier--Walsh} coefficients
of $f$, and they satisfy
\[
\hat f(S) := \langle f, \chi_S \rangle= {\mathbb E}[{f \chi
_S}] .
\]
Note that $\hat f(\varnothing)$ corresponds to the average ${\mathbb
E}[{f}]$.
As in classical Fourier analysis, if $f$ is some Boolean function, its
Fourier(--Walsh) coefficients provide information on the ``regularity''
of $f$.

Of course one may find many other orthonormal bases for $L^2(\{-1,1\}
^n)$, but there are many situations
for which this particular set of functions $(\chi_S)_{S\subset\{
1,\ldots,n \} }$ arises naturally. First of all there is a well-known theory
of Fourier analysis on groups, a theory which is particularly simple
and elegant on Abelian groups (thus including
our special case of $\{-1, 1\}^n$, but also $\R\slash\Z$, $\R$ and
so on).
For Abelian groups, what turns out to be relevant for doing harmonic analysis
is the set $\hat G$ of \textit{characters} of $G$ (i.e., the group
homomorphisms from $G$ to $\C^*$). In our case of $G=\{-1,1\}^n$,
the characters are precisely our functions $\chi_S$ indexed by
$S\subset\{1,\ldots,n\}$ since they satisfy $\chi_S(x\cdot y) = \chi
_S(x) \chi_S(y)$.

These functions also arise naturally if one performs simple random walk
on the hypercube (equipped with the Hamming graph structure), since they
are the eigenfunctions of the heat kernel on $\{-1,1\}^n$.

Last but not least, the basis $(\chi_S)$ turns out to be particularly
adapted to our study of \textit{noise sensitivity}.
Indeed if $f \dvtx \Omega_n \to\R$ is any real-valued function, then the
correlation between $f(\omega)$ and $f(\omega^\eps)$
is nicely expressed in terms of the Fourier coefficients of $f$ as follows:
%
%
\begin{eqnarray}\label{e.correlationFourier}
{\mathbb E}[{f(\omega) f(\omega^{\eps})}] & =& {\mathbb
E}\biggl[{\biggl( \sum_{S_1} \hat f(S_1) \chi_{S_1}(\omega)\biggr)
\biggl( \sum_{S_2} \hat f(S_2) \chi_{S_2}(\omega^\eps)
\biggr)}\biggr] \nonumber\\
&=& \sum_{S} \hat f(S)^2{\mathbb E}[{\chi_S(\omega) \chi
_S(\omega^\eps)}]
\\
&=& \sum_S \hat f(S)^2 (1-\eps)^{|S|} .\nonumber 
\end{eqnarray}

Therefore, the ``level of sensitivity'' of a Boolean function is
naturally encoded in its Fourier coefficients. More precisely,
for any real-valued function $f \dvtx \Omega_n \to\R$, one can consider
its \textit{energy spectrum} $E_f$ defined by
\[
E_f(m):= \sum_{|S|=m} \hat f(S)^2\qquad \forall m \in\{1,\ldots, n\}
.
\]
Since $\Cov(f(\omega), f(\omega^\eps)) = \sum_{m=1}^n E_f(m)
(1-\eps)^m$, all the information we need is contained
in the energy spectrum of $f$. As argued in the \hyperref[sec1]{Introduction},
a~function of ``high frequency'' will be sensitive to noise
while a function of ``low frequency'' will be stable. This allows us to
give an equivalent definition of \textit{noise sensitivity} (recall
Definition \ref{d.NS}):
\begin{proposition}\label{pr.NSeq}
A sequence of Boolean functions $f_{n} \dvtx \{-1,1\}^{m_n} \to\{0,1\}$
is (asymptotically) \textup{noise sensitive} if and only if, for any
$k\geq1$
\[
\sum_{m=1}^k \sum_{|S|=m} \hat f_n(S)^2 = \sum_{m=1}^k E_{f_n}(m)
\mathop{\longrightarrow}_{n\to\infty} 0 .
\]
\end{proposition}

Before introducing the notion of influence, let us give a simple example:
%
\begin{Example*} Let $\Phi_n$ be the majority function on $n$
variables (a function which is of obvious interest in computer
science). 
$\MAJ_n(x_1,\ldots,x_n):= \sign(\sum x_i)$, where $n$ is an odd integer
here. It is possible in this case to explicitly compute its Fourier
coefficients, and
when $n$ goes to infinity, one ends up with the following asymptotic
formula (see \cite{ODonnellThesis} for a nice overview and
references therein):
\[
E_{\MAJ_n}(m) = \sum_{|S|=m} \widehat{\MAJ_n}(S)^2 = \cases{
\displaystyle \frac{4}{ \pi m 2^m}
\pmatrix{m-1\vspace*{2pt}\cr\dfrac{m-1}{2}}+O(m/n), &\quad if $m$ is odd, \vspace*{5pt}\cr
0, &\quad if $m$ is even.}
\]


Figure \ref{f.MAJ} represents the shape of the energy spectrum of
$\MAJ_n$: its spectrum is concentrated on low frequencies which is typical
%
%
\begin{figure}

\includegraphics{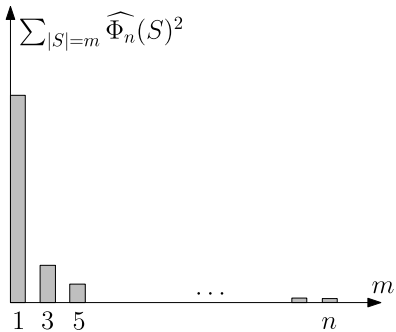}

\caption{Spectrum of the majority function $\Phi_n$.}\label{f.MAJ}
\end{figure}
of stable functions.
\end{Example*}

Henceforth, most of this paper will be concerned with the description
of the Fourier
expansion of Boolean functions (and more specifically of their energy spectrum).

\subsection{\texorpdfstring{Notion of influence.}{Notion of influence}}\label{ss.influence}

%
%
If $f\dvtx\Omega_n\to\R$ is a (real-valued) function,
the influence of a variable $k\in[n]$ is a quantity which measures
by how much (on average) the function $f$ depends on the fixed variable $k$.
For this purpose, we introduce the functions
\[
\nabla_k f \dvtx \cases{
\Omega_n \to \R,\cr
\omega \mapsto f(\omega) -f(\sigma_k\cdot\omega),
}\qquad
\mbox{for all } k\in[n] ,
\]
where $\sigma_k$ acts on $\Omega_n$ by flipping the $k$th bit
(thus $\nabla_k f$ corresponds to a~discrete derivative along the
$k$th bit).

The \textit{influence} $\Inf_k(f)$ of the $k$th variable is defined in
terms of $\nabla_k f$ as follows:
\[
\Inf_k(f) := \| \nabla_k f \|_1= \| f(\omega)- f(\sigma_k \cdot
\omega
)\|_1 .
\]
\begin{remark}
If $f\dvtx\Omega_n \to\{0,1\}$ is a Boolean function corresponding to an
event $A\subset\Omega_n$ (i.e., $f=1_A$), then $\Inf_k(f)=\Inf_k(A)$
is the probability that the $k$th bit is pivotal for $A$ (i.e., $\Inf
_k(f)= {\mathbb P}[{f(\omega)\neq f(\sigma_k \cdot\omega
)}]$).
\end{remark}
\begin{remark}\label{r.mono}
If $f\dvtx\Omega_n \to\R$ is a \textit{monotone} function [i.e., $f(\omega_1)
\le f(\omega_2)$ when $\omega_1\le\omega_2$], then
notice that
\[
\hat f(\{ k\}):= {\mathbb E}\bigl[{f \chi_{\{k\}}}\bigr] = {\mathbb
E}[{f(x_1,\ldots,x_n) x_k}] =
\tfrac1 2 \Inf_k(f)
\]
by monotonicity of $f$. This gives a first hint that influences are
closely related to the Fourier expansion of $f$.
\end{remark}

We define the \textit{total influence} of a (real-valued) function $f$ to
be $\Inf(f):= \sum\Inf_k(f)$. This notion is very relevant in the
study of sharp thresholds (see \mbox{\cite{MR671248,MR1371123}}).
Indeed if $f$ is a monotone function, then by the Margulis--Russo
formula (see, e.g., \cite{GrimmettGraphs})
\[
\frac d {dt} \bigg|_{1/2} \PP_p(f) = \Inf(f) = \sum_k \Inf_k(f) .
\]
This formula easily extends to all $p\in[0,1]$. In particular, for a
monotone event $A$, a ``large'' total influence implies
a ``sharp'' threshold for $p\mapsto\PP_p(A)$.

As it has been recognized for quite some time already (since Ben
Or/Linial), the set of all influences $\Inf_k(f)$, $k\in[n]$ carries
important information
about the function $f$.
Let us then call $\InfV(f):= (\Inf_k(f))_{k\in[n]}$ the \textit{influence vector} of~$f$.
We have already seen that the $L^1$ norm of the influence vector
encodes properties\vspace*{1pt} ``sharp threshold''-type for $f$ since by definition
$\Inf(f)= \| {\InfV}(f)\|_1$.
The $L^2$ norm of this influence vector will turn out to be a key
quantity in the study of noise sensitivity of Boolean functions
$f\dvtx\Omega_n \to\{0,1\}$. We thus define
(following the notations of \cite{MR2001m60016})
\[
\II(f):= \sum_k \Inf_k(f)^2 = \|{\InfV}(f)\|_2^2 .
\]

For Boolean functions $f$ (i.e., with values in $\{0,1\}$), these
notions [$\Inf(f)$ and $\II(f)$]
are intimately related with the above Fourier expansion of $f$. Indeed
we will see in the next section that if
$f \dvtx \Omega_n \to\{0,1\}$, then
\[
\Inf(f) = 4 \sum_S |S| \hat f(S)^2 .
\]

If one assumes furthermore that $f$ is \textit{monotone}, then from Remark
\ref{r.mono}, one has
%
%
\begin{equation}\label{e.H}
\frac1 4 \II(f) = \sum_{k}\hat f(\{k\})^2 = \sum_{|S|=1}\hat f(S)^2
= E_f(1),
\end{equation}
which corresponds to the ``weight'' of the level-one Fourier
coefficients [this property also holds for real-valued
functions, but we will use the quantity~$\II(f)$ only in the Boolean case].

We will conclude by the following general \textit{philosophy}
that we will encounter throughout the rest of the survey (especially in Section
\ref{s.hyper}): if a~function $f \dvtx  \{0,1\}^n \to\R$ is such that each
of its variables has a ``very small'' influence (i.e., $\ll\frac1
{\sqrt{n}}$),
then $f$ should have a behavior very different from a ``Gaussian'' one.
We will see an illustration of this \textit{rule} in the context of
anomalous fluctuations (Lem\-ma~\ref{l.torusfpp}). In the Boolean case,
these functions (such that all their variables have very small
influence) will be noise sensitive (Theorem \ref{th.ns}),
which is not characteristic of Gaussian nor White noise behavior.

\subsection{\texorpdfstring{Fast review on percolation.}{Fast review on percolation}}\label{ss.perco}
We will only briefly recall what the model is, as well as some of its
properties that will be used throughout the text.
For a complete account on percolation see \cite{Grimmettnewbook} and more
specifically in our context the lecture notes \cite{07100856}.

\begin{figure}[t]

\includegraphics{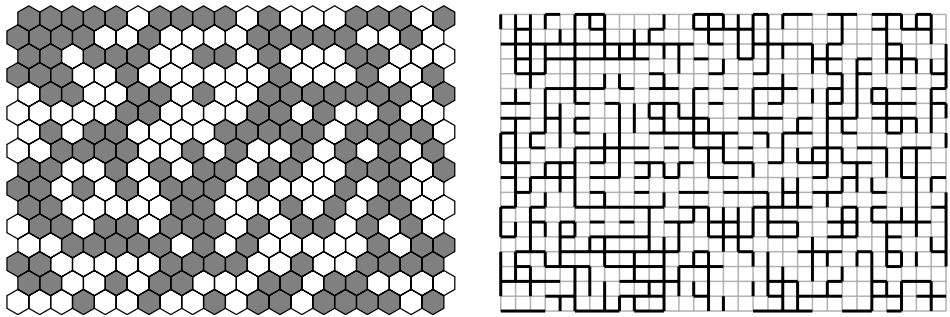}

\caption{Illustrations by Oded representing two critical percolation
configurations, respectively, on $\T$ and on $\Z^2$. The sites of the
triangular grid
are represented by hexagons.}\label{f.perco}
\end{figure}

We will be concerned mainly in two-dimensional percolation, and we will
focus on two lattices, $\Z^2$ and the triangular lattice $\T$ (see Figure
\ref{f.perco}).
All the results stated for $\Z^2$ in this text are also valid for
percolations on ``reasonable'' two-dimensional translation invariant
graphs for which RSW is known to hold.

Let us describe the model on $\Z^2$. Let $\E^2$ denote the set of edges
of the graph~$\Z^2$. For any $p\in[0,1]$ we define a random subgraph of
$\Z^2$ as follows: independently for each edge $e\in\E^2$, we keep this
edge with probability $p$ and remove it with probability $1-p$.
Equivalently, this corresponds to defining a random configuration
$\omega\in\{-1,1\}^{\E^2}$ where, independently for each edge
$e\in\E^2$, we declare the edge to be open [$\omega(e)=1$] with
probability $p$ or closed [$\omega (e)=-1$] with probability $1-p$. The
law of the so-defined random subgraph (or configuration) is denoted
by~$\PP_p$. In percolation theory, one is interested in large-scale
connectivity properties of the random configuration~$\omega$.

In particular as one raises the level $p$, above a certain critical
parameter~$p_c(\Z^2)$, an infinite cluster (almost surely) appears.
This corresponds to the well-known \textit{phase transition} of
percolation. By a famous theorem of Kesten this transition takes place at
$p_c(\Z^2)=\frac1 2$.\vspace*{2pt}

Percolation is defined similarly on the triangular grid $\T$, except
that on this lattice we will rather consider \textit{site}-percolation
(i.e., here we keep each \textit{site} with probability $p$). Also\vspace*{1pt} for
this model, the transition happens at the critical point $p_c(\T
)=\frac
1 2$. 

%
%

\begin{figure}

\includegraphics{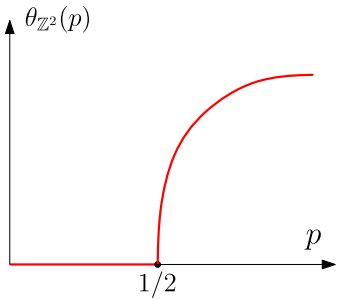}

\caption{An illustration of the density function on $\Z^2$.}\label{fig4}
\vspace*{-2pt}
\end{figure}

The phase transition can be measured with the \textit{density function}
$\theta_{\Z^2}(p): = \PP_p (0\stackrel{\omega}{\longleftrightarrow}
\infty)$
which encodes important properties of the
large-scale connectivities of the random configuration $\omega$: it
corresponds to the density averaged over the space $\Z^2$ of
the (almost surely unique) infinite cluster. The shape of the function
$\theta_{\Z^2}$ is pictured in Figure \ref{fig4} (notice the infinite
derivative at~$p_c$).

Over the last decade, the understanding of the critical regime has
undergone remarkable progress, and Oded himself obviously
had an enormous impact on these developments. The main ingredients of
this productive period were the introduction of the
$\mathrm{SLE}$ processes by Oded (see the survey on $\mathrm{SLEs}$ by
Steffen Rohde in the present volume)
and the proof of conformal invariance on $\T$ by Stanislav Smirnov
\cite{MR1851632}.

\begin{figure}[b]

\includegraphics{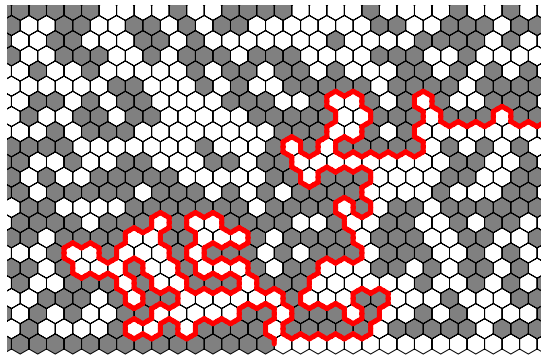}

\caption{A sample of an exploration path on the triangular lattice.}\label{fig5}
\end{figure}

At this point one cannot resist showing another famous (and everywhere
used) illustration by Oded
representing an exploration path on the triangular lattice (see Figure \ref{fig5}); this red
curve which turns right on black hexagons
and left on the white ones, asymptotically converges toward $\mathrm
{SLE}_6$ (as the mesh size goes to 0).

The proof of conformal invariance combined with the detailed
information given by the $\mathrm{SLE}_6$ process
enabled one to obtain very precise information on the critical and \textit{near-critical} behavior of $\T$-percolation. For
instance,
it is known that on the triangular lattice, the density function
$\theta
_{\T}(p)$ has the following behavior near $p_c=1/2$:
\[
\theta(p) = (p-1/2)^{5/36 + o(1)} ,
\]
when $p\to1/2+$ (see \cite{07100856}).

\begin{figure}

\includegraphics{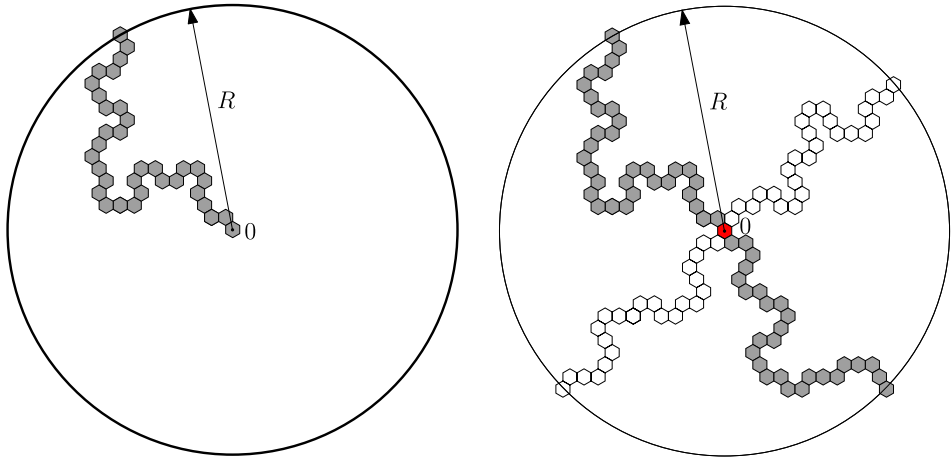}

\caption{A realization of the \textit{one-arm} event is pictured on the
left; the \textit{four-arm} event is pictured on the right.}\label{f.armsevents}
\end{figure}

In the rest of the text, we will often rely on two types of percolation
events: namely
the \textit{one-arm} and \textit{four-arm}
events. They are defined as follows: for any radius $R>1$, let $A_R^1$
be the event that the site 0 is connected to distance~$R$ by some open
path. Also, let $A_R^4$ be the event that there are four ``arms''
of\vadjust{\goodbreak}
alternating color from the site 0 (which can be of either color)
to distance $R$
(i.e., there are four connected paths, two open, two closed from~0 to
radius $R$ and the closed paths lie
between the open paths). See Figure~\ref{f.armsevents} for a
realization of each event.

%
%

It was proved in \cite{MR2002k60204} that the probability of the one-arm
event decays like
\[
{\mathbb P}[{A_R^1}]:= \alpha_1(R) = R^{-5/{48}
+o(1)} .
\]
For the four-arms event, it was proved by Smirnov and Werner in
\cite{MR1879816} that its probability decays like
\[
{\mathbb P}[{A_R^4}]:= \alpha_4(R) = R^{-5/4 +o(1)} .
\]

The three exponents we encountered concerning $\theta_{\T}$, $\alpha_1$
and $\alpha_4$ (i.e., $\frac5 {36}$, $\frac5 {48}$ and $\frac5 4$)
are known as \textit{critical exponents}.

The \textit{four-arm} event will be of particular importance throughout
the rest of this survey. Indeed suppose the four arms
event holds at some site $x\in\T$ up to some large distance $R$. This
means that the site $x$ carries important
information about the large scale connectivities within the euclidean
ball $B(x,R)$. Changing the status of $x$ will drastically
change the ``picture'' in $B(x,R)$. We call such a point a \textit{pivotal
point} up to distance $R$.

Finally it is often convenient to ``divide'' these arm-events into
different scales. For this purpose, we introduce
$\alpha_4(r,R)$ (with $r\le R$) to be the probability that the four-arm
event is realized from radius $r$ to radius $R$
[$\alpha_1(r,R)$ is defined similarly for the one-arm event]. By
independence on disjoint sets, it is clear that if $r_1\le r_2\le r_3$ then
one has $\alpha_4(r_1,r_3) \le\alpha_4(r_1,r_2) \alpha_4(r_2,r_3)$.
A very useful property known as \textit{quasi-multiplicativity}
claims that up to constants, these two expressions are the same (this
makes the division into several scales practical).
This property can be stated as follows.
\begin{proposition}[(Quasi-multiplicativity \cite
{MR88k60174})]\label{pr.quasi}
For any $r_1\le r_2\le r_3$, one has (both for $\Z^2$ and $\T$ percolations)
\[
\alpha_4(r_1,r_3) \asymp\alpha_4(r_1,r_2) \alpha_4(r_2,r_3),
\]
where the constant involved in $\asymp$ are uniform constants.
\end{proposition}

\begin{figure}[t]

\includegraphics{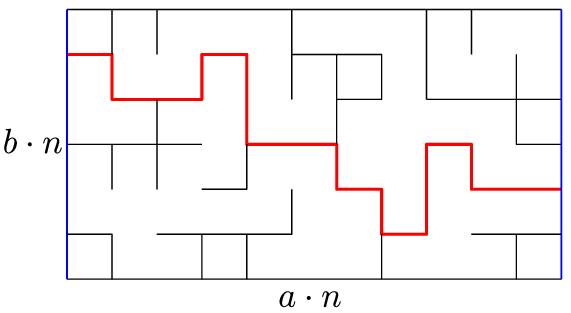}

\caption{Left to right crossing event in a rectangle and RSW
theorem.}\label{fig7}
\end{figure}

See \cite{07100856,NolinKesten} for more details. Note also that the
same property holds for the one-arm event (but is much easier to prove:
it is an easy consequence of the RSW theorem which is stated in Figure
\ref{fig7}).

%

Figure \ref{fig7} represents a configuration for which the
left--right crossing of $[a\cdot n$, $b\cdot n]$ is realized; in this case
we define $f_n(\omega):=1$; otherwise we define $f_n(\omega):= 0$.
There is an extremely useful result known as \textit{RSW theorem} which states
that for any $a,b >0$, asymptotically the probability of crossing the
rectangle ${\mathbb P}[{f_n=1}]$ remains bounded away from
0 and 1.

\begin{figure}[t]

\includegraphics{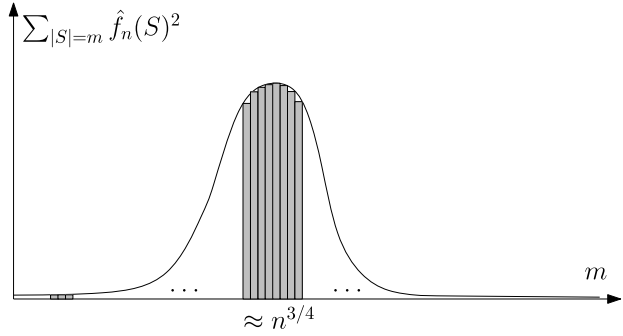}

\caption{Expected shape of the energy spectrum of percolation crossing
events. We will see in Section~\protect\ref{s.geo}
that most of the spectral mass of $f_n$ is indeed localized around
$n^{3/4}$. Compare this shape of the spectrum with the above spectrum
of the majority function $\Phi_n$.}\label{f.CrossingEnergy}
\end{figure}

To end this preliminary section, let us sketch what will be one of the
main goals throughout this survey: if
$f_n$ denotes the left--right crossing event of a large rectangle
$[a\cdot n, b\cdot n ]$ (see Figure \ref{fig7}),
then one can consider these observables as Boolean functions. As such
they admit a Fourier expansion.
Understanding the sensitivity of percolation will correspond to
understanding the energy spectrum of such observables.
We will see by the end of this survey that, as $n\to\infty$, the
energy spectrum of $f_n$ should roughly look as in Figure \ref
{f.CrossingEnergy}.

%
%

\section{\texorpdfstring{The ``hypercontractivity'' approach.}{The ``hypercontractivity'' approach}}\label{s.hyper}
In this section, we will describe the first noise sensitivity estimates
as well as anomalous fluctuations for first passage percolation (FPP).

The notion of influence will be essential throughout this section;
recall that we defined
$\nabla_k f = f -f\circ\sigma_k$ and $\Inf_k(f)= \| \nabla_k f\|_1$.
It is fruitful to consider~$\nabla_k f$ 
as a function from $\Omega_n$ into $\R$. Indeed, its Fourier
decomposition is directly related to the Fourier decomposition
of $f$ itself in the following way: it is straightforward to check that
for any $S\subset[n]$
\[
\widehat{\nabla_k f} (S) = \cases{
2 \hat f (S), &\quad if $k\in S$, \cr
0, &\quad else.}
\]

If one assumes furthermore that $f$ is Boolean (i.e., $f\dvtx  \Omega_n \to
\{
0,1\}$), then our discrete derivatives $\nabla_k f$ take their values
in $\{-1,0,1\}$; this implies
in particular that for any $k\in[n]$
\[
\Inf_k(f) = \|\nabla_k f \|_1 = \| \nabla_k f\|_2^2 .
\]
Using Parseval, this enables us to write the total influence of $f$ in
terms of its Fourier coefficients as follows:
%
%
\begin{eqnarray}\label{e.inf1}
\Inf(f) &=& \sum_k \Inf_k(f) = \sum_k \| \nabla_k f\|_2^2
\nonumber\\
&=& \sum_k \| \widehat{\nabla_k f}\|_2^2 \\
&=& 4 \sum_S |S| \hat f(S)^2,\nonumber
\end{eqnarray}
since each ``frequency'' $S\subset[n]$ appears $|S|$ times. 
Before the appearance of the notion of ``noise sensitivity'' there have
been several works whose focus were to get lower bounds on
influences. For example (see \cite{MR1371123}), if one wants to
provide general criteria for sharp thresholds, one needs to obtain
lower bounds
on the total influence $\Inf(f)$. The strategy to do this in great
generality was developed in \cite{KKL}, where they proved that for any
balanced Boolean function $f$ (balanced meaning here ${\mathbb P}
[{f=1}] =1/2)$,
there always exists a variable whose influence is greater than $O(\frac
{\log n }{n})$.

Before treating in detail the case of noise sensitivity, let us
describe in an informal way what was
the ingenious approach from \cite{KKL} to obtain lower bounds on
influences. Let us consider some arbitrary
balanced Boolean function $f$.
We want to show that at least one of its variables has large influence
$\geq c \frac{\log n} n$.
Suppose all its influences $\Inf_k(f)$ are ``small'' (this would need to
be made quantitative); this means that all the functions $\nabla_k f$
have small $L^1$ norm.
Now if $f$ is Boolean (into $\{0,1\}$), then as we have noticed above,
$\nabla_k f$ is almost Boolean (its values are in $\{-1,0,1\}$);
hence $\| \nabla_k f\|_1$ being small implies that $\nabla_k f$ has a
small support.
Using the intuition coming from Weyl--Heisenberg uncertainty, $\widehat
{\nabla_k f}$ should then be quite spread; in
particular, most of its spectral mass should be concentrated on high
frequencies.

This intuition (which is still vague at this point) somehow says that
having small influences pushes the spectrum of $\nabla_k f$
toward high frequencies. Now summing up
as we did in (\ref{e.inf1}), but only restricting ourselves to
fequencies~$S$ of size smaller than some large (well-chosen) $1\ll M
\ll n$, one obtains
%
%
\begin{eqnarray}\label{e.estimate}
\sum_{0<|S|< M} \hat f(S)^2 & \le& 4 \sum_{0<|S| < M} |S| \hat f(S)^2
\nonumber\\
&=& \sum_k \sum_{0<|S|<M} \widehat{\nabla_k f}(S)^2
\nonumber\\[-8pt]\\[-8pt]
& \mbox{``}\ll\mbox{''} & \sum_k \| \widehat{\nabla_k f} \|_2^2 \nonumber\\
&=& \Inf(f) ,\nonumber
\end{eqnarray}
where in the third line, we used the informal statement that $\widehat
{\nabla_k f}$ should be supported on high frequencies if
$f$ has small influences. Now recall that we assumed $f$ to be
balanced, hence
\[
\sum_{|S|>0} \hat f(S)^2 = 1/4 .
\]
Therefore, in the above equation (\ref{e.estimate}), if we are in the
case where
a positive fraction of the Fourier mass of $f$ is concentrated below
$M$, then (\ref{e.estimate})
says that $\Inf(f)$ is much larger than one. In particular, at least
one of the influences has to be ``large.''
If, on the other hand, we are in the case where most of the spectral
mass of $f$ is supported on frequencies of
size higher than $M$, then we also obtain that $\Inf(f)$ is large by
the previous formula
\[
\Inf(f)= 4 \sum_S |S| \hat f(S)^2 .
\]

In \cite{KKL}, this intuition is converted into a proof. The main
difficulty here is to formalize, or, rather, to implement,
the above argument, that is, to obtain
spectral information on functions with values in $\{-1,0,1\}$ knowing
that they have small support.
This is done \cite{KKL} using techniques brought from harmonic
analysis, namely \textit{hypercontractivity}.

\subsection{\texorpdfstring{About hypercontractivity.}{About hypercontractivity}}
First, let us state what hypercontractivity corresponds to.
Let $(K_t)_{t\ge0}$ be the heat kernel on $\R^n$.
Hypercontractivity
is a statement which\vadjust{\goodbreak} quantifies how functions are regularized under
the heat flow. The statement, which goes back to Nelson and Gross, can
be simply stated as follows:
\begin{lemma}[(Hypercontractivity)]\label{l.hyp}
If $1<q<2$, there is some $t=t(q)>0$ (which does not depend on the
dimension $n$) such that for any $f\in L^q(\R^n)$,
\[
\| K_t \ast f \|_2 \le\| f \|_q .
\]
\end{lemma}

The dependence $t=t(q)$ is explicit but will not concern us in the
Gaussian case. Hypercontractivity is thus a regularization statement:
if one starts with some initial ``rough'' $L^q$ function $f$ outside of
$L^2$ and waits long enough [$t(q)$] under the heat flow, we end up
being in
$L^2$ with a good control on its $L^2$ norm.

This concept has an interesting history that is nicely explained in
O'Don\-nell's lectures notes (see \cite{OdonnellBlog}). It was
originally invented by Nelson in \cite{MR0210416}
when he needed regularization estimates on Free Fields (which are the
building blocks of quantum field theory) in order to apply these in
``constructive field theories.''
It was then generalized by Gross in his elaboration of Logarithmic
Sobolev Inequalities~\cite{MR0420249},
which are an important tool in analysis. Hypercontractivity is
intimately related to these Log--Sobolev inequalities (they are
somewhat equivalent concepts)
and thus has many applications in the theory of semi-groups, mixing of
Markov chains and so on.

We now state the result in the case which concerns us, the hypercube.
For any $\rho\in[0,1]$, let $T_\rho$ be the following ``noise
operator'' on the functions of the hypercube:
recall from the preliminary section that if $\omega\in\Omega_n$, we
denote by $\omega^\eps$ an $\eps$-noised configuration of $\omega$.
For any $f \dvtx  \Omega_n \to\R$, we define $T_{\rho} f\dvtx  \omega\mapsto
{\mathbb E}[{f(\omega^{1-\rho}) \mid \omega}]$. This
noise operator
acts in a very simple way on the Fourier coefficients
\[
T_\rho\dvtx  f=\sum_S \hat f(S) \chi_S \mapsto\sum_S \rho^{|S|} \hat
f(S) \chi_S .\vadjust{\goodbreak}
\]

We have the following analog of Lemma \ref{l.hyp}:
\begin{lemma}[(Bonami--Gross--Beckner)]
For any $f \dvtx  \Omega_n \to\R$,
\[
\| T_\rho f \|_{2} \le\| f\|_{1+\rho^2} .
\]
\end{lemma}

The analogy with the classical Lemma \ref{l.hyp} is clear: the Heat
flow is replaced here by the random walk
on the hypercube.

Before applying hypercontractivity to noise sensitivity, let us sketch
how this functional inequality helps to implement
the above idea from \cite{KKL}. If a~Boolean function $f$ has small
influences, its discrete derivatives $\nabla_k f$ have small support.
Now these functions have values in $\{-1,0,1\}$; thus for any $1<q<2$
we have that $\|\nabla_k f\|_q = (\| \nabla_k f\| _2)^{2/q} \ll\|
\nabla_k f\|_2$
(because of the small support of $\nabla_k f$). Now applying
hypercontractivity (with $q=1+\rho^2$), we obtain that $\|T_\rho
\nabla
_k f\|_2 \ll\| \nabla_k f\|_2$.
Written on the Fourier side this means that
\[
\sum_S \rho^{2 |S|} \widehat{\nabla_k f}(S)^2 \ll\sum_S \widehat
{\nabla
_k f}(S)^2 ,
\]
and this happens only if most of the spectral mass of $\nabla_k f$ is
supported on high frequencies.
It remains to make the above heuristics precise in the case which
interests us here.

\subsection{\texorpdfstring{Applications to noise sensitivity.}{Applications to noise sensitivity}}
Let us now see hypercontractivity in action. As in the \hyperref[sec1]{Introduction} we
are interested in the noise sensitivity of a sequence of Boolean functions
$f_n \dvtx  \Omega_{m_n} \to\{0 , 1\}$. A deep theorem from
\cite{MR2001m60016} can be stated as follows:
\begin{theorem}[\cite{MR2001m60016}]\label{th.ns}
Let $(f_n)_n$ be a sequence of Boolean functions. If $\II(f_n) \to
0$, when $n\to\infty$, then $(f_n)_n$ is \textup{noise sensitive};
that is, for any $\eps>0$, the correlation ${\mathbb E}
[{f_n(\omega) f_n(\omega^\eps)}] - {\mathbb E}[{f}
]^2$ converges to 0.
\end{theorem}

The theorem is true independently of the speed of convergence of $\II
(f_n)$. Nevertheless, if one assumes that there is some
exponent $\delta>0$, such that $\II(f_n) \le(m_n)^{-\delta}$, then the
proof is quite simple as was pointed out to us by Jeff Steif, and furthermore
one obtains some ``logarithmic bounds'' on the sensitivity of $f_{n}$.
We will restrict ourselves to this stronger assumption since it will be
sufficient for our application to Percolation.
\begin{remark}
$\!\!\!$If the Boolean functions $(f_n)_n$ are assumed to
be~\mbox{monotone}, it is interesting to note that as we observed in (\ref{e.H}),
$\II(f_n) = \sum_{|S|=1} \hat f_n(S)^2$. So $\II(f_n)$ corresponds here
to the level-1 Fourier weights. Thus in the monotone case,
Theorem \ref{th.ns} says that if asymptotically there is no weight on
the level one coefficients, then there is no weight on any finite level
Fourier\vadjust{\goodbreak} coefficient (of course the Boolean hypothesis on $f_n$ is also
essential here).
In particular, in the monotone case, the condition $\II(f_n)\to0$ is
equivalent to noise sensitivity.
\end{remark}
\begin{pf*}{Proof of Theorem \ref{th.ns} \textup{[under the stronger assumption $\II
(f_n) \le(m_n)^{-\delta}$ for some $\delta>0$]}}
The spirit of the proof is similar to the one carried out in~\cite{KKL},
but the target is different here. Indeed in \cite{KKL},
the goal
was to obtain
good lower bounds on the total influence in great generality; for
example, the (easy) sharp threshold encountered by
the majority function around $p=1/2$ fits into their framework.
Majority is certainly not a sensitive function, so Theorem \ref{th.ns}
requires more assumption
than \cite{KKL} results. Nevertheless the strategy will be similar: we
still use the ``spectral decomposition'' of $f$ with respect to its
``partial derivatives'' $\nabla_k f$
\[
\Inf(f) = \sum_k \| \nabla_k f\|_2^2 = \sum_k \|\widehat{\nabla_k
f} \|
_2^2 = 4 \sum_S |S| \hat f(S)^2 .
\]
But now we want to use the fact that $\II(f) = \sum_k \| \nabla_k f
\|
_2^4$ is very small. If~$f$ is Boolean, this implies that the (almost) Boolean
$\nabla_k f$ have very small support (and this can be made
quantitative). Now, again we expect $\widehat{\nabla_k f}$ to be
spread, but this time,
we need more: not only $\nabla_k f$ have high frequencies, but in some
sense ``all their frequencies'' are high leaving no mass (after summing
up over the $m_n$ variables)
to the finite level Fourier coefficients. This is implemented using
hypercontractivity, following similar lines as in~\cite{KKL}.

Let us then consider a sequence of Boolean functions $f_n \dvtx  \Omega
_{m_n} \to\{0,1\}$ satisfying $\II(f_n) \le(m_n)^{-\delta}$ for some
exponent $\delta>0$.
We want to show that there is some constant $M=M(\delta)$, such that
\[
\sum_{0< |S| < M\log(m_n)} \hat f_n(S)^2 \to0 ,
\]
which gives a quantitative (logarithmic) noise sensitivity statement
\begin{eqnarray}
&&\sum_{0<|S| < M \log(m_n)} \hat f_n(S)^2 \nonumber\\
&&\qquad \le 4 \sum_{0< |S| < M
\log
(m_n)} |S| \hat f_n(S)^2 \nonumber\\
&&\qquad=  \sum_k \sum_{ 0< |S| < M \log(m_n)} \widehat{\nabla_k
f}(S)^2\nonumber\\
&&\qquad\le \sum_k \biggl(\frac1 {\rho^2}\biggr)^{M \log(m_n)} \| T_\rho\nabla_k
f\|
_2^2 \nonumber\\
&&\qquad \le \sum_k \biggl(\frac1 {\rho^2}\biggr)^{M \log(m_n)} \| \nabla_k f\|
_{1+\rho
^2}^2
\nonumber\\
\eqntext{\mbox{(by hypercontractivity)}.}
\end{eqnarray}
Now, since $f$ is Boolean, one has $\|\nabla_k f\|_{1+\rho^2} = \|
\nabla
_k f \|_2^{2/(1+\rho^2)}$, hence
\begin{eqnarray*}
&&\sum_{0<|S| < M \log(m_n)} \hat f_n(S)^2 \\
&&\qquad \le \rho^{-2 M \log(m_n)
}\sum_k \| \nabla_k f \|_2^{4/(1+\rho^2)}\\
&&\qquad= \rho^{-2 M \log(m_n) }\sum_k \Inf_k(f)^{2/(1+\rho
^2)}\\
&&\qquad \le \rho^{-2 M \log(m_n)} (m_n)^{\rho^2/(1+\rho^2)} \biggl(
\sum_k
\Inf_k(f)^2 \biggr)^{1/({1+\rho^2})}\qquad\mbox{(By H\"older)}
\\
&&\qquad = \rho^{-2 M \log(m_n)} (m_n)^{\rho^2/(1+\rho^2)} \II
(f_n)^{
1/ ({1+\rho^2})} \\
&&\qquad\le \rho^{-2 M \log(m_n)} (m_n)^{({\rho^2 - \delta})/({1+\rho
^2})}
.
\end{eqnarray*}

Now by choosing $\rho\in(0,1)$ close enough to 0, and then by choosing
$M=M(\delta)$ small enough, we obtain the desired
logarithmic noise sensitivity.\vspace*{8pt}

\textit{Application to percolation}.
We want to apply the above result to the case of percolation crossings.
Let $D$ be some smooth domain of the plane (not necessarily simply connected),
and let $\p_1$, $\p_2$ be two smooth arcs on the boundary $\p D$. For
all $n\ge1$, let $D_n \subset\frac1 n \Z^2$ be a domain
approximating $D$, and call
$f_n$, the indicator function of a left to right crossing in $D_n$
from~$\p_1^n$ to~$\p_2^n$ (see Figure \ref{f.revealment} in Section \ref
{s.algo}). Noise sensitivity of percolation means that the sequence of
events $(f_n)_{n\ge1}$
is noise sensitive. Using Theorem \ref{th.ns}, it would be enough (and
necessary) to show that $\II(f_n) \to0$. But if we want a
self-contained proof here, we would
prefer to have at our disposal the stronger claim $\II(f_n) \le
n^{-\delta}$ for some $\delta>0$ (here $m_n\asymp n^2$).

There are several ways to see why this stronger claim holds.
The most natural one is to get good estimates on the probability for an
edge $e$ to be pivotal. Indeed, recall that
in the monotone case, $\II(f_n):= \sum_{\mathrm{edges}\ e\in D_n }
{\mathbb P}[e \mbox{ is}\break  \mbox{Pivotal}]^2$ (see Section
\ref{ss.influence}).
This probability, without considering boundary effects, is believed to
be of order $n^{-5/4}$, which indeed makes $\II(f_n) \approx n^2\cdot
(n^{-5/4})^2 \approx n^{-1/2}$
decrease to 0 polynomially fast. This behavior is now known in the case
of the triangular grid thanks to Schramm's $\mathrm{SLEs}$ and
Smirnov's proof of
conformal invariance (see \cite{MR1879816} where the relevant
critical exponent is computed).

At the time of \cite{MR2001m60016}, of course, critical exponents were not
available (and anyway, they still remain unavailable today on $\Z^2$),
but Kesten had estimates which implied that for some $\eps>0$,
${\mathbb P}[{e \mbox{ is Pivotal}}] \le n^{-1-\eps}$
[which is enough to obtain $\II(f_n) \le n^{-\delta}$].

Furthermore, an ingenious alternative way to obtain the polynomial
convergence of $\II(f_n)$ was developed in \cite{MR2001m60016} which
did not
need Kesten's results.
This alternative way, on which we will say a few words in Section \ref
{s.algo}, in some sense prefigured the randomized algorithm approach
that we will describe in the next section.\vspace*{8pt}


\textit{Some words on the general Theorem} \ref{th.ns} \textit{and its proof}.
The proof of the general result is a bit more involved than the one we
outlined here.
The main lemma is as follows:
\begin{lemma}
There exist absolute constants $C_k$ for all $k\geq1$, such that for
any monotone Boolean function $f$ one has
\[
\sum_{|S|= k} \hat f(S)^2 \le C_k \II(f) (-\log\II(f) )^{k-1} .
\]
\end{lemma}

This lemma ``mimics'' a result from Talagrand \cite{MR1401897}.
Indeed Proposition~2.2 in \cite{MR1401897} can be translated as
follows: for any monotone
Boolean function $f$, its level-$2$ Fourier weight [i.e., $\sum_{|S|=2}
\hat f(S)^2$] is bounded by $O(1) \II(f) \log(1/\break\II(f))$.
It obviously implies Theorem \ref{th.ns} in the monotone case; the
general case being deduced from it by a monotonization procedure.
Hypercontractivity is used in the proof of this lemma.


\subsection{\texorpdfstring{Anomalous fluctuations, or chaos.}{Anomalous fluctuations, or chaos}}

In this section, we will outline how hypercontractivity was used in
\cite{MR2016607} in order to prove that shape fluctuations in the
model of first passage percolation are sub-Gaussian.

\subsubsection{\texorpdfstring{The model of first passage percolation
(FPP).}{The model of first passage percolation (FPP)}}
Let us start with the model and then state the theorem proved in
\cite{MR2016607}.
First passage percolation can be seen as a model of a \textit{random
metric} on $\Z^d$; it is defined simply as follows:
independently for each edge $e\in\E^d$, fix the length of $e$ to be $1$
with probability $1/2$, $2$ else.
In greater generality, the lengths of the edges are i.i.d. nonnegative
random variables, but here, following \cite{MR2016607},
we will restrict ourselves to the above uniform distribution on $\{1,2\}
$ to simplify the exposition (see \cite{MR2451057} for an
extension to more general laws).
In fact, in \cite{MR2016607}, they handle the slightly more general case
of a uniform distribution on $\{a,b\}$ with $0 <a < b$ but we decided
here to stick to the case
$a=1$ and $b=2$ since it makes the analogy with the previous results on
influences simpler.

For any $\omega\in\{ 1,2\}^{\E^d}$, this defines a (random) metric,
$\mathrm{dist}_\omega$, on $\Z^d$ satisfying for any $x,y\in\Z^d$,
\[
\mathrm{dist}_\omega(x,y):= \inf_{\gamma\dvtx \mathrm{path}\ \mathrm{from}\ x\
\mathrm{to}\ y}
l_\omega(\gamma),
\]
where $l_\omega(\gamma)$ is the length of $\gamma$.

Using sub-additivity, it is known that the renormalized ball $\frac1 n
B_\omega(0,n)$ converges toward a deterministic shape (which
can be in certain cases computed explicitly).

\subsubsection{\texorpdfstring{Fluctuations around the limiting
shape.}{Fluctuations around the limiting shape}}

The fluctuations around the asymptotic limiting shape have received
tremendous interest over the last 15 years.
In the two-dimensional case, using very beautiful combinatorial
bijections with random matrices,
certain cases of \textit{directed} last passage percolation (where the law
on the edges is taken to be geometric or exponential) have been understood
very deeply. For example, it is known \cite{MR1737991}
that the fluctuations of the ball of radius $n$\vadjust{\goodbreak}
(i.e., the points whose last passage time are below $n$)
around $n$ times its asymptotic deterministic shape, are of order
$n^{1/3}$, and the law of these fluctuations properly renormalized
follow the Tracy--Widom distribution
(as do the fluctuations of the largest eigenvalue of GUE ensembles).

``\textit{Universality}'' is believed to hold for these models in the
sense that the behavior of the fluctuations around the deterministic shape
should not depend on the ``microscopic'' particularities of the model
(e.g., the law on the edges). The shape itself does depend of course.
In particular in the above model
of (nondirected) first passage percolation in dimension $d=2$,
fluctuations are widely believed to be also of order $n^{1/3}$
following as well the Tracy--Widom law.
Still, the current state of understanding of this model is far from
this conjecture.

Kesten first proved that the fluctuations of the ball of radius $n$
were at most $\sqrt{n}$ (which did not exclude yet Gaussian behavior).
Benjamini, Kalai and Schramm then strengthened this result by showing
that the fluctuations were sub-Gaussian. This does not yet reach the
conjectured $n^{1/3}$-fluctuations, but their approach has the great
advantage to be very general; in particular their result holds in any
dimension $d\geq2$.

Let us now state their main theorem concerning the fluctuations of the
metric $\dist_\omega$.
%
\begin{theorem}[\cite{MR2016607}]\label{th.fpp}
For all $d$, there exists an absolute constant $C=C(d)$ such that in
$\Z^d$
\[
\operatorname{var}(\mathrm{dist}_\omega(0,v)) \le C\frac{|v|}{\log|v|}
\]
for any point $v\in\Z^2$, $|v|\geq2$.
\end{theorem}

\subsubsection{\texorpdfstring{Link with ``noise sensitivity.''}{Link with ``noise sensitivity''}}

This result about sub-Gaussian fluctuations might seem at first
disjoint from our initial study of noise sensitivity, but they turn out
to be intimately
related. First of all the methods to understand one or the other, as we
will see, follow very similar lines.
But also, as is very nicely explained in \cite{ChatterjeeChaos}, the
phenomenon of ``anomalous fluctuations'' is in some sense
equivalent to a certain ``noise sensitivity'' of the geodesics of
first-passage-percolation. More precisely, the variance of the first
passage time
is of order ${\mathbb E}[{|\gamma(\omega(0))\cap\gamma(\omega
(\tau))|}]$, where
$\tau\sim\mathcal{E}(1)$ is an exponential variable.
Thus we see that if the metric, or rather the geodesic, is highly
perturbed when
the configuration is noised; then the distances happen to be very
concentrated. Chatterjee calls this phenomenon \textit{chaos}. Of course,
our short description
here was informal since in our present setup there might be many
different geodesics between two fixed points. The above link between
concentration
and sensitivity discovered by Chatterjee works very nicely in the
context of Maxima of Gaussian processes (which in that case
arise a.s. at a single point, or a single ``geodesic'' in a geometrical
context) (see \cite{ChatterjeeChaos} for more details).

\subsubsection{\texorpdfstring{The simpler case of the torus.}{The simpler case of the torus}}
Following the approach of \cite{MR2016607}, we will first consider the
case of the torus. The reason for this is that it is a~much simpler case.
Indeed, in the torus, for the least-passage time that we will consider,
any edge will have up to constant the same influence, while in the case
of $\mathrm{dist}_\omega(0,v)$,
edges near the endpoints $0$ or $v$ have a high influence on the
outcome (in some sense there is more symmetry and invariance to play with
in the case of the torus).

Let $\T_m^d$ be the $d$-dimensional torus $(\Z/m\Z)^d$.
As in the above (lattice) model, independently for each edge of $\T
_m^d$, we choose its length to be either~$1$ or $2$. We are interested here
in the smallest (random) length among closed paths $\gamma$ ``turning''
around the torus along the first coordinate $\Z/ m\Z$ (i.e., these
paths $\gamma$, once
projected onto the first cycle, have winding number one). In
\cite{MR2016607}, this is called the shortest \textit{circumference}.
For any configuration $\omega\in\{1,2\}^{E(\T_m^d)}$, call $\Circ
_m(\omega)$ this shortest circumference.
\begin{theorem}[\cite{MR2001m60016}]\label{th.torus}
There is a constant $C>0$ (which does not depend on the dimension $d$),
such that
\[
\operatorname{var}(\Circ_m(\omega)) \leq C \frac m {\log m}.
\]
\end{theorem}
\begin{remark}
A similar analysis as the one carried out below works in greater
generality: if $G=(V,E)$ is some finite connected graph endowed with a
random metric $d_\omega$ with $\omega\in\{1,2\}^{\otimes E}$,
then one can obtain bounds on the fluctuation of the random diameter
$D=D_\omega$ of $(G,d_\omega)$. See \cite{MR2016607}, Theorem 2, for a
precise statement in this more general context.
\end{remark}
\begin{pf*}{Sketch of proof of Theorem \ref{th.torus}}
In order to highlight the similarities with the above case of noise
sensitivity of percolation, we will not follow
exactly~\cite{MR2016607}; it will be more ``hands-on'' with the
disadvantage of being less general (we take $a=1$ and $b=2$).

As before, for any edge $e$, let us consider the gradient along the
edge $e$, $\nabla_e \Circ_m$; these gradient functions have values in
$\{-1, 0, 1 \}$, since changing the length of $e$ can only have this
effect on the circumference. Note here that if the lengths of edges
were in $\{a,b\}$ for any fixed choice of
$0<a<b$, then it would not always be the case anymore that $\nabla_e
\Circ_m \in\{ -(b-a), 0, b-a \}$. Even though this is not crucial
here, this is why we stick to the case $a=1$ and $b=2$. See \cite{MR2001m60016}
for a way
to overcome this lack of ``Boolean behavior.''

Since our gradient functions have values in $\{-1,0,1\}$, we end up
being in the same setup as in our previous study;
influences are defined in the same way and so on.
We sill see that our gradient functions (which are ``almost Boolean'')
have small support, and hypercontractivity will imply the desired bounds.

Let us work in the general case of a function $f \dvtx  \{-1,1\}^n \to\R$,
such that for any variable $k$,
$\nabla_k f \in\{-1,0,1\}$. We are interested in $\operatorname
{var}(f)$ [and we
want to show that if ``influences are small'' then $\operatorname
{var}(f)$ is small].
It is easy to check that the variance can be written
\[
\operatorname{var}(f) = \frac1 4 \sum_k \sum_{\varnothing\neq S
\subset[n]} \frac1
{|S|} \widehat{\nabla_k f} (S)^2 .
\]
We see on this expression, that if variables have very small influence,
then as previously, the almost Boolean $\nabla_k f$ will be
of high frequency. Heuristically, this should then imply that
\begin{eqnarray*}
\operatorname{var}(f) & \ll& \sum_k \sum_{S\neq\varnothing}
\widehat{\nabla_k f}
(S)^2 \\
& = & \sum_k \Inf_k(f) .
\end{eqnarray*}

We prove the following lemma on the link between the fluctuations of
a~real-valued function $f$ on $\Omega_n$ and its influence
vector. 
%
\begin{lemma}\label{l.torusfpp}
Let $f \dvtx  \Omega_n \to\R$ be a (real-valued) function such that each of
its discrete derivative $\nabla_k f, k\in[n]$ have their values in
$\{-1,0,1\}$.
If we assume that the influences of $f$ are small in the following sense:
there exists some $\alpha>0$ such that for any $k\in\{1, \ldots, n\}
$, $\Inf_k(f) \le n^{-\alpha}$, then there is some constant
$C=C(\alpha)$,
such that
\[
\operatorname{var}(f) \le\frac{C} {\log n} \sum_k \Inf_k(f) .
\]
\end{lemma}

Before proving the lemma, let us see that in our special case of first
passage percolation, the assumption on small
influences is indeed verified.
Since the edges' length is in $\{1,2\}$, the smallest contour $\Circ
_m(\omega)$ in $\T_m^d$ around the first coordinate
lies somewhere in $[m, 2m]$. Hence, if $\gamma$ is a geodesic (a path
in the torus) satisfying $l(\gamma)=\Circ_m(\omega)$,
then $\gamma$ uses at most $2 m$ edges. There might be several
different geodesics minimizing the circumference.
Let us choose randomly one of these in an ``invariant'' way and call it
$\tilde\gamma$. For any edge $e\in E(\T_m^d)$, if by changing the length
of $e$, the circumference increases, then $e$ has to be contained in
any geodesic $\gamma$, and in particular in $\tilde\gamma$. This
implies that
${\mathbb P}[{\nabla_e \Circ_m(\omega) > 0}] \le
{\mathbb P}[{e \in\tilde\gamma}]$. By
symmetry we obtain that
\[
{\mathbb P}[{\nabla_e \Circ_m(\omega) \neq0}] \le2
{\mathbb P}[{e \in\tilde\gamma}] .
\]

As we have seen above, $\nabla_e \Circ_m \in\{-1,0,1\}$; therefore
$\Inf_e(\Circ_m) \le O(1)\times{\mathbb P}[{e \in\tilde\gamma
}]$.
Now using the symmetries both of the torus $\T_m^d$ and of our
observable $\Circ_m$, if $\tilde\gamma$ is chosen in an appropriate
invariant way
(uniformly among all geodesics would work), then it is clear that all
the vertical edges (the edges which, once projected on the first cycle,
project on a single vertex) have the same probability to lie in $\tilde
\gamma$; the same goes for horizontal edges. In particular,
\[
\sum_{``\mathrm{vertical}\ \mathrm{edges}{\fontsize{8.36pt}{9.36pt}\selectfont{\mbox{''}}}\ e} {\mathbb P}[{e \in\tilde
\gamma}] \le{\mathbb E}[{|\tilde\gamma|}] \le2
m .
\]
Since there are $O(1)m^d$ vertical edges, the influence of these is
bounded by $O(1)m^{1-d}$; the same goes for horizontal edges.
All together this gives the desired assumption needed in Lemma \ref
{l.torusfpp}. Applying this lemma, we indeed obtain that
\[
\operatorname{var}(\Circ_m(\omega)) \le O(1) \frac m {\log m} ,
\]
where the constant does not depend on the dimension $d$ (since the
dimension helps us here).
\begin{pf*}{Proof of the Lemma \ref{l.torusfpp}}
As for noise sensitivity, the proof relies on implementing
hypercontractivity in the right way.
\begin{eqnarray*}
\operatorname{var}(f) & = & \frac1 4 \sum_k \sum_{S\neq\varnothing}
\frac1 {|S|}
\widehat{\nabla_k f} (S)^2 \\
&\le& \frac1 4 \sum_k \sum_{0 < |S| < c \log n} \widehat{\nabla_k f}
(S)^2 + \frac{O(1)} {\log n} \sum_k \Inf_k(f) .
\end{eqnarray*}

Hence it is enough to bound the contribution of small frequencies,
$0<|S|< c \log n$, for some constant $c$ which will be chosen
later. As previously we have for any $\rho\in(0,1)$ and using
hypercontractivity,
%
%
\begin{eqnarray}\qquad
\sum_k \sum_{0 < |S| < c \log n} \widehat{\nabla_k f} (S)^2 &\le&
\rho
^{-2 c \log n} \sum_k \| T_\rho\nabla_k f \|_2^2 \nonumber\\
&\le& \rho^{-2 c \log n} \sum_k \| \nabla_k f\|_{1+\rho^2} ^2
\nonumber\\
&=& \rho^{-2 c \log n} \sum_k \Inf_k(f)^{2/(1+\rho^2)}
\nonumber\\[-8pt]\\[-8pt]
&\le& \rho^{-2 c \log n} \Bigl(\sup_{k} \Inf_k(f) \Bigr)^{
({1-\rho
^2})/({1+\rho^2})} \sum_{k} \Inf_k(f) \nonumber\\
&\le& \rho^{-2 c \log n} n^{-\alpha({1-\rho^2})/({1+\rho^2}) }
\sum
_k \Inf_k(f) \nonumber\\
&&\eqntext{\mbox{(by our assumption)}.}
\end{eqnarray}

Now fixing $\rho\in(0,1)$, and then choosing the constant $c$ depending
on $\rho$ and~$\alpha$, the lemma follows
(by optimizating on the choice of $\rho$, one could get better
constants).\vspace*{-2pt}
\end{pf*}

\subsubsection{\texorpdfstring{Some hints for the proof of Theorem
\protect\ref{th.fpp}.}{Some hints for the proof of Theorem \protect\ref{th.fpp}}}

The main difficulty here is that the quantity of interest, $f(\omega
):=\dist_\omega(0,v)$, is not anymore invariant under a large class of
graph automorphisms. This lack of symmetry makes the study of
influences more difficult. (e.g., as was noticed above, edges
near the endpoints 0 or $v$ will have high influence). To gain some
more symmetry, the authors in \cite{MR2016607} rely on a nice
``averaging'' procedure. The idea is as follows: instead of looking at
the (random) distance form 0 to $v$, they first pick a
point $x$ randomly in the mesoscopic box $[-|v|^{1/4}, |v|^{1/4}]^d$
around the origin and then consider the distance from this
point $x$ toward $v+x$. Let~$\tilde f$ denote this function [$\dist
_\omega(x,v+x)$]. $\tilde f$ uses extra randomness compared to~$f$, but
it is clear
that ${\mathbb E}[{f}]= {\mathbb E}[{\tilde f}
]$, and it is not hard to see that when $|v|$
is large, $\operatorname{var}(f) \asymp\operatorname{var}(\tilde
f)$. Therefore it is enough to
study the fluctuations
of the more symmetric $\tilde f$. (We already see here that thanks to
this averaging procedure, the endpoints 0 and $v$ no longer have a high
influence.)
In some sense, along geodesics, this procedure ``\textit{spreads}'' the
influence on the $|v|^{1/4}$-neighborhood of the geodesics. More precisely,
if $e$ is some edge, the influence of this edge is bounded by $2
{\mathbb P}[{e \in x + \gamma}]$, where $\gamma$ is
chosen among geodesics from 0 to $v$.
Now, as we have seen in the case of the torus, geodesics are
essentially one-dimensional [of length less than $O(1)|v|$]; this is
still true on the mesoscopic
scale: for any box $Q$ of radius $m:=|v|^{1/4}$, $|\gamma\cap Q| \le
O(1) m$. Now by considering the mesoscopic box around $e$, it is like
moving a ``line''
in a box of dimension $d$; the probability for an edge to be hit by
that ``line'' is of order $m^{1-d}$. Therefore the influence of any
edge $e$ for the
``spread'' function $\tilde f$ is bounded by $O(1) |v|^{(1-d)/4} \le
O(1)|v|^{-1/4}$.
This implies the needed assumption in Lemma \ref{l.torusfpp} and hence
concludes the sketch of proof of Theorem \ref{th.torus}. See
\cite{MR2016607} for a more detailed
proof.\vspace*{-2pt}
\begin{remark}
In this survey, we relied on Lemma \ref{l.torusfpp}, since its proof is
very similar to the noise sensitivity proof. In \cite{MR2016607}, the authors
use (and reprove with better constants) an inequality from Talagrand
\cite{MR1303654} which states that there is some universal
constant $C>0$ such that for any
$f\dvtx \{0,1\}^n \to\R$,
\[
\operatorname{var}(f) \le C \sum_k \frac{\| \nabla_k f\|_2^2} { 1+
\log(\| \nabla_k
f \|_2 / \| \nabla_k f\|_1)} .\vspace*{-2pt}
\]
\end{remark}

Conclusion: Benjamini, Kalai and Schramm \cite{MR2001m60016,MR2016607}
developed
multiple and very interesting techniques. The results of \cite{MR2016607}
have since
been extended to more general laws \cite{MR2451057}, but
essentially, their control of the variance in $\Omega(n/\log n)$ is to
this day still the best.
The paper \cite{MR2001m60016} had a profound impact on the field.
As we will see, some of the ideas present in \cite{MR2001m60016} already
announced some ideas of the next section.\vadjust{\goodbreak}

\section{\texorpdfstring{The randomized algorithm approach.}{The randomized algorithm approach}}\label{s.algo}
In this part, we will describe the quantitative estimates on noise
sensitivity obtained in \cite{SchrammSteif}.
Their applications to the model of dynamical percolation will be
described in the last section of this survey.
But before we turn to the remarkable paper \cite{SchrammSteif}, where
Schramm and Steif introduced deep techniques to control Fourier
spectrums of
general functions,
let us first mention and explain that the idea of using randomized
algorithms was already present in essence in \cite{MR2001m60016},
where they
used an algorithm in order to prove that $\II(f)$ converges quickly
(polynomially) toward 0.

\subsection{\texorpdfstring{First appearance of randomized algorithm
ideas.}{First appearance of randomized algorithm ideas}}
In \cite{MR2001m60016}, as we have seen above, in order to prove that
percolation crossings are asymptotically noise sensitive, the authors needed
the fact that $\II(f_n)= \sum_k \Inf_k(f_n)^2 \to0$ (see Theorem
\ref{th.ns});
if furthermore this $L^2$ quantity converges to zero more quickly than
a polynomial of the number of
variables, $(m_n)^{-\delta}$ for some $\delta>0$, then the proof of
Theorem \ref{th.ns} is relatively simple as we outlined above.
This fast convergence to zero of $\II(f_n)$ was guaranteed by the work
of Kesten (in particular his work \cite{MR88k60174} on hyperscaling
from which follows the
fact that the probability for a point to be pivotal until distance $m$
is less than $O(1) m^{-1-\alpha}$ for some $\alpha>0$).

Independently of Kesten's approach, the authors provided in \cite{MR2001m60016}
a different way of looking at this problem (an approach
more in the spirit of noise sensitivity). They noticed the remarkable
property that if a \textit{monotone} Boolean function $f$ happens to be
correlated very little with majority functions (for all subsets of the
bits); then $\II(f)$ has to be very small, and hence the function has
to be
sensitive to noise. They obtained a quantitative version of this
statement that we briefly state here.

Let $f\dvtx \{-1,1\}^n \to\{0,1\}$ be a Boolean function. We want to use
its correlations with majority functions. Let us define these:
for all $K \subset[n]$, define the majority function on the subset
$K$ by $\mathbf{M}_K(x_1,\ldots,x_n):= \sign\sum_K x_i$ (where
$\sign0 := 0$ here).
The correlation of the Boolean function $f$ with these majority
functions is measured by
\[
\Lambda(f):= {\max_{K \subset[n]}} | {\mathbb E}[{f
\mathbf{M}_K}] | .
\]

Being correlated very little with majority functions corresponds to
$\Lambda(f)$ being very small.
The following quantitative theorem about correlation with majority is
proved in \cite{MR2001m60016}.
\begin{theorem}\label{th.cm}
There exists a universal constant $C>0$ such that for any $f \dvtx \Omega_n
\to\{0,1\}$ \textup{monotone}
\[
\II(f) \le C \Lambda(f)^2 \bigl(1-\log\Lambda(f)\bigr) \log n
\]
(the result remains valid if $f$ has values in $[0,1]$
instead).\vadjust{\goodbreak}
\end{theorem}

With this result at their disposal, in order to obtain fast
convergence of $\II(f_n)$ to zero
in the context of percolation crossings, the authors of \cite{MR2001m60016}
investigated the correlations
of percolation crossings $f_n$ with majority on subsets $K\subset
[n]$. They showed that there exist $C,\alpha>0$ universal constants,
so that for any subset of the lattice $K$, $|{\mathbb E}[{f_n
\mathbf{M}_K}]
| \le C n^{-\alpha}$. For this purpose, they used
a nice appropriate randomized algorithm. We will not detail this
algorithm used in \cite{MR2001m60016}, since it was considerably
strengthened in \cite{SchrammSteif}. We will now describe the approach
of \cite{SchrammSteif} and then return to ``correlation with
majority'' using
the stronger algorithm from \cite{SchrammSteif}.

\subsection{\texorpdfstring{The Schramm/Steif approach.}{The Schramm/Steif approach}} 

The authors in \cite{SchrammSteif} introduced the following beautiful
and very general idea: suppose a real-valued function, $f\dvtx\Omega_n \to
\R$
can be \textit{exactly} computed with a randomized algorithm $A$, so that
every fixed variable is used by the algorithm $A$ only with small probability;
then this function $f$ has to be of ``high frequency'' with quantitative
bounds which depend on how
unlikely it is for any variable to be used by the algorithm.

\subsubsection{\texorpdfstring{Randomized algorithms, revealment and
examples.}{Randomized algorithms, revealment and examples}}
Let us now define more precisely what types of randomized algorithms
are allowed here.
Take a function $f\dvtx\{-1,1\}^n \to\R$. We are looking for algorithms
which compute the output of $f$ by examining some of the bits (or variables)
\textit{one by one}, where the choice of the next bit may depend on the
set of bits discovered so far, plus if needed, some additional randomness.
We will call an algorithm satisfying this property a \textit{Markovian}
(\textit{randomized}) \textit{algorithm}. Following \cite{SchrammSteif}, if $A$ is a
Markovian algorithm
computing the function $f$, we will denote by $J\subset[n]$ the
(random) set of bits examined by the algorithm.

In order to quantify the property that variables are unlikely to be
used by an algorithm $A$, we define the \textit{revealment} $\reveal=
\reveal_A$ of the
algorithm $A$ to be the supremum over all variables $i \in[n]$ of the
probability that $i$ is examined by $A$. In other words,
\[
\reveal=\reveal_A = \sup_{i \in[n]} {\mathbb P}[{i \in
J}] .
\]

We can now state one of the main theorems from \cite{SchrammSteif} (we
will sketch its proof in the next subsection).
\begin{theorem}\label{th.algo}
Let $f \dvtx \{-1,1\}^n \to\R$ be a function. Let $A$ be a Markovian
randomized algorithm for $f$ having revealment $\reveal=\reveal_A$.
Then for every $k=1,2,\ldots.$ The ``level k''-Fourier coefficients of
$f$ satisfy
\[
\sum_{S\subset[n], |S|=k} \hat f(S)^2 \le\reveal_A k \| f \|
_2^2 .
\]
\end{theorem}
\begin{remark}
If one is looking for a Markovian algorithm computing the output of the
majority function on $n$ bits, then it is clear that
the only way to proceed is to examine variables one at a time (the
choice of the next variable being irrelevant since they all play the
same role).
The output will not be known until at least half of the bits are
examined; hence the revealment for majority is at least $1/2$.
\end{remark}

In the case of percolation crossings, as opposed to the above case of
majority, one has to exploit the ``richness'' of the percolation
picture in order to find algorithms which detect crossings while
examining very few bits. A~natural idea for a left-to-right crossing event
in a large rectangle is to use an \textit{exploration path}. The idea of
an exploration path, which was highly influential in the \hyperref[sec1]{Introduction}
by Schramm
of the $\mathrm{SLE}$ processes, was pictured in Section \ref
{ss.perco} in the case of the triangular lattice.

More precisely, for any $n\geq1$, let $D_n$ be a domain consisting of
hexagons of mesh $1/n$ approximating the square $[0,1]^2$, or more generally
any smooth ``quad'' $\Omega$ with two prescribed arcs $\p_1, \p_2
\subset\Omega$ (see Figure \ref{f.revealment}).
We are interested in the left-to-right crossing event (in the general
setting, we look at the crossing event from $\p_1$ to $\p_2$ in $D_n$).
Let $f_n$ be the corresponding Boolean function and call $\gamma_n$ the
``exploration path'' as in Figure \ref{f.revealment} (which starts at
the upper left corner $a$). We run this exploration path until it
reaches either the bottom side (in which case $f_n=0$) or the
right-hand side (corresponding to $f_n=1$).

%
%
\begin{figure}[b]

\includegraphics{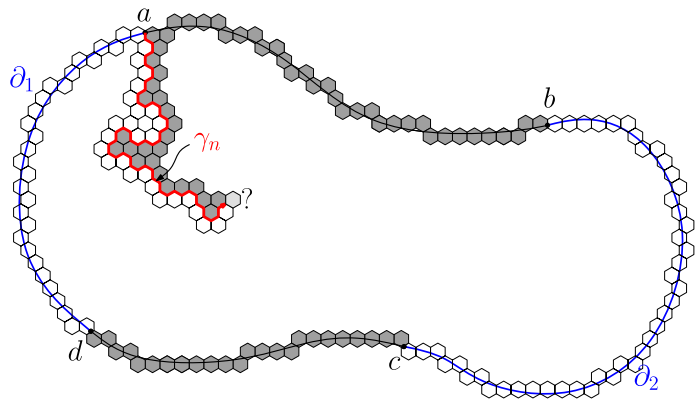}

\caption{The random interface $\gamma_n$ is discovered one site at a
time (which makes our algorithm Markovian). The exploration stops
when $\gamma_n$ reaches either $(bc)$ (in which case $f_n=1$) or $(cd)$
($f_n=0$).}\label{f.revealment}
\end{figure}

This thus provides us with a Markovian algorithm to compute $f_n$ where
the set $J=J_n$ of bits examined by the algorithm is the set of
``hexagons'' touching $\gamma_n$ on either side.
The nice property of both the exploration path~$\gamma_n$ and its
$1/n$-neighborhood $J_n$, is that they have a scaling limit when the
mesh $1/n$ goes to zero, this scaling limit being the well-known $\SLE
_6$. (This scaling limit of the exploration path was as we mentioned
above one of the main motivations of Schramm to introduce these $\SLE$
processes.) This scaling limit is a.s. a random fractal curve in
$[0,1]^2$ (or $\Omega$) of Hausdorff dimension~$7/4$. This means that
asymptotically, the exploration path uses a~very small amount of all
the bits included in this picture. With some more work (see~\cite{MR1879816,07100856}),
we can see that inside the domain (not near the corner or the sides),
the probability for a point to be on $\gamma_n$ is of order $n^{-1/4 +
o(1)}$, where $o(1)$
goes to zero as the mesh goes to zero.

One therefore expects the revealment $\reveal_n$ of this algorithm to
be of order $n^{-1/4+o(1)}$. But the corner/boundaries
have a nontrivial contribution here: for example, in this setup, the
single hexagon on the upper-left corner of the domain (where the
interface starts) will be used systematically by the algorithm making
the revealment equal to one! There is an easy way to handle this
problem: the idea in \cite{SchrammSteif} is to use some additional
randomness and to start the exploration path from a random point $x_n$
on the left-hand side of $[0,1]^2$. Doing so, this ``smoothes'' the
singularity of the departure
along the left boundary. There is a small counterpart to this: with
this setup, one interface might not suffice to detect the existence of
a left-to-right crossing, and a second interface starting from the same
random point~$x_n$ might be needed (see~\cite{SchrammSteif} for more
details). Using \textit{arms exponents} from~\cite{MR1879816} and ``quasimultiplicativity'' of arms events
\cite{MR88k60174,SchrammSteif}, it can be checked that indeed
the revealment of this modified
algorithm is $n^{-1/4+o(1)}$.
\begin{theorem}
Let $(f_n)_{n\geq1}$ be the left-to-right crossing events in domains
$(D_n)_{n\geq1}$ approximating the unit square (or more generally a
smooth domain $\Omega$).
Then there exists a sequence of Markovian algorithms, whose revealments
$\reveal_n$ satisfy that for any $\eps>0$,
\[
\reveal_n \le C n^{-1/4+\eps} ,
\]
where $C=C(\eps)$ depends only on $\eps$.
\end{theorem}


Therefore, applying Theorem \ref{th.algo}, one obtains.
\begin{corollary}\label{cor.algo}
Let $(f_n)_n$ be the above sequence of crossing events. Let~$\Spec
_{f_n}$ denote the spectral sample of these Boolean functions.
Since $\|f_n \|_2 \le1$, we obtain that for any sequence $(L_n)_{n\geq1}$,
%
%
\begin{equation}\label{e.lowertail}
{\mathbb P}[{0 < |\Spec_{f_n}| < L_n}] \le\reveal_n
L_n^2 .
\end{equation}

In particular, this implies that for any $\eps>0$, ${\mathbb P}
[{0< |\Spec_{f_n}| < n^{1/8-\eps}}] \to0$.
\end{corollary}

This result gives precise lower bound information about the ``spectrum
of percolation'' (or its ``energy spectrum''). It implies a\vadjust{\goodbreak}
``polynomial sensitivity'' of crossing events in the sense
that for any level of noise $\eps_n \gg n^{-1/8}$, we have that
${\mathbb E}[{f_n(\omega) f_n(\omega^{\eps_n})}] -
{\mathbb E}[{f_n}]^2 \to0$.\vadjust{\goodbreak}
%
%
\begin{remark}
\begin{itemize}
\item
Note that equation (\ref{e.lowertail}) gives a good control on the
lower tail of the spectral distribution, and as we will see in the last
section, these
lower-tail estimates are essential in the study of dynamical percolation.
\item
Similar results are obtained by the authors in \cite{SchrammSteif} in
the case of the $\Z^2$-lattice, except that for this lattice, conformal
invariance
and convergence toward $\SLE_6$ are not known; therefore critical
exponents such as $1/4$ are not available; still \cite{SchrammSteif}
obtains polynomial controls (thus strengthening~\cite{MR2001m60016})
but with
small exponents (their value coming from RSW estimates).
\end{itemize}
\end{remark}

\subsubsection{\texorpdfstring{Link with correlation with majority.}{Link with correlation with majority}}
Before proving Theorem \ref{th.algo}, let us briefly return to the
original motivation of these types of algorithms.
Suppose we have at our disposal the above Markovian Algorithms for the
left-to-right crossings $f_n$ with small revealments $\reveal_n=
n^{-1/4+o(1)}$;
then it follows easily that the events $f_n$ are mostly not correlated with
majority functions; indeed let $n\geq1$ and $K\subset D_n \approx
[0,1]^2$ some fixed subset of the bits.
By definition of the revealment, we have that ${\mathbb E}[{|K
\cap J_n|}] \le
|K| \reveal_n$. This means that on average, the algorithm visits very
few variables
belonging to~$K$. Since
\[
|{\mathbb E}[{f_n \mathbf{M}_K}]| = |{\mathbb E}
[{f_n {\mathbb E}[{\mathbf{M}_K \mid J_n }]}
]| \le{\mathbb E}[{| {\mathbb E}[{\mathbf{M}_K \mid
J_m}]|}] ,
\]
and using the fact that on average, $|K\cap J_m|$ is small compared to
$|K|$, it is easy to deduce that there is some $\alpha>0$
such that
\[
\Lambda(f_n) = {\max_{K}} | {\mathbb E}[{f_n
\mathbf{M}_K}] | \le
n^{-\alpha
} .
\]

This, together with Theorem \ref{th.ns} and Theorem \ref{th.cm} implies
a logarithmic noise sensitivity for percolation.
In \cite{MR2001m60016}, they rely on another algorithm which instead of
following interfaces, in some sense ``invades'' clusters attached to
the left-hand side.
Since clusters are of fractal dimension $91/48$, intuitively their
algorithm, if boundary issues can be properly taken care of, would give
a bigger revealment
of order $n^{-5/48 + o(1)}$ (the notion of revealment only appeared in
\cite{SchrammSteif}). So the major breakthrough in \cite{SchrammSteif}
is that they simplified
tremendously the role played by the algorithm by introducing the notion
of revealment and they noticed a more direct link with the Fourier transform.
Using their correspondence between algorithm and spectrum, they greatly
improved the control on the Fourier spectrum (polynomial v.s. logarithmic).
Furthermore, they improved the randomized algorithm.

\subsubsection{\texorpdfstring{Proof of Theorem
\protect\ref{th.algo}.}{Proof of Theorem \protect\ref{th.algo}}}

Let $f \dvtx \{-1,1\}^n \to\R$ be some real-valued function, and consider
$A$, a Markovian algorithm associated to it with\vadjust{\goodbreak}
revealment $\reveal=\reveal_A$. Let $k\in\{1,\ldots, n\}$;
we want to bound from above the size of the level-$k$ Fourier
coefficients of $f$ [i.e., $\sum_{|S|=k} \hat f(S)^2$].
For this purpose, one introduces the function $g= g^{(k)}=\sum_{|S|=k}
\hat f(S) \chi_S$, which is the projection of $f$ onto the
subspace of level-$k$ functions. By definition, one has $\sum_{|S|=k}
\hat f(S)^2 = \|g \|_2^2$.

Very roughly speaking, the intuition goes as follows: if the
revealment~$\reveal$ is small, then for low level $k$
there are few ``frequencies'' in $g^{(k)}$ which will be ``seen'' by
the algorithm. More precisely, for any fixed ``frequency''~$S$, if $|S|=k$
is small, then
with high probability none of the bits in $S$ will be visited by the algorithm.
This means that ${\mathbb E}[{g^{(k)}\mid J}]$ (recall
$J$ denotes the set of
bits examined by
the algorithm) should be of small
$L^2$ norm compared to $g^{(k)}$. Now since $f= {\mathbb E}
[{f\mid J}] = \sum
_k {\mathbb E}[{g^{(k)} \mid J}]$, most of the Fourier transform
should be supported on high frequencies. There is some difficulty in
implementing this intuition, since the conditional expectations
${\mathbb E}[{g^{(k)} \mid J}]$
are not orthogonal.


Following \cite{SchrammSteif} very closely, one way to implement this
idea goes as follows:\looseness=1
%
%
\begin{equation}\label{e.fg}
{\mathbb E}[{g^2}] = {\mathbb E}[{f g}] =
{\mathbb E}[{f {\mathbb E}[{g \mid J}]}] \le
\| f\|_2 \| {\mathbb E}[{g\mid J}] \|_2 .
\end{equation}\looseness=0
As hinted above, the goal is to control the $L^2$ norm of ${\mathbb
E}[{g \mid J}]$. In order to achieve this, it will be
helpful to interpret
${\mathbb E}[{g\mid J}]$ as the expectation of a~random
function $g_J$ whose
construction is explained below.

Recall that $J$ is the set of bits examined by the
randomized algorithm $A$. Since the randomized algorithm depends on the
configuration $\omega\in\{-1,1\}^n$ and possibly some additional randomness,
one can view $J$ as a random variable on some extended probability
space $\tilde\Omega$, whose elements can be represented
as $\tilde\omega= (\omega, \tau)$ ($\tau$~corresponding here to the
additional randomness).

For any function $h \dvtx \{-1,1\}^n \to\R$, one defines the \textit{random}
function $h_J$ which corresponds to the function $h$ where bits in $J$
are fixed to match with what was examined by the algorithm. More
exactly, if $J(\tilde\omega)=J(\omega,\tau)$ is the random set of bits
examined,
then the random function $h_J = h_{J(\tilde\omega)}$ is the function
on $\{-1,1\}^n$ defined by $h_{J(\omega,\tau)} (\omega') := h(\omega'')$,
where $\omega''=\omega$ on $J(\omega,\tau)$ and $\omega''=\omega
'$ on
$\{1,\ldots,n\} \setminus J(\omega,\tau)$.
Now, if the algorithm has examined the set of bits $J=J(\tilde\omega
)$, then with the above definition it is clear that the conditional
expectation ${\mathbb E}[{h \mid J}]$ [which is a
measurable function of
$J=J(\tilde\omega)$] corresponds to averaging
$h_J$ over configurations $\omega'$ [in other words we average on the
variables outside of $J(\tilde\omega)$]; this can be written as
\[
{\mathbb E}[{h \mid J}] = \int h_J:= \hat h_J (\varnothing),
\]
where the integration $\int$ is taken with respect to the uniform
measure on $\omega'\in\Omega_n$.
In particular ${\mathbb E}[{h}]= {\mathbb E}[{\int
h_J}] ={\mathbb E}[{\hat h_J(\varnothing)}]$. Since
$(h^2)_J = (h_J)^2$, it follows
that
%
%
\begin{equation}\label{e.hJ}
\|h\|_2^2 = {\mathbb E}[{h^2}] = {\mathbb E}\biggl[{\int
h_J^2}\biggr] = {\mathbb E}[{\|h_J \|_2^2}] .
\end{equation}

Recall from (\ref{e.fg}) that it only remains to control $\| {\mathbb
E}[{g\mid J}] \|_2^2 = {\mathbb E}[{\hat
g_J(\varnothing)^2 }]$.
For this purpose,\vadjust{\goodbreak} we apply Parseval to the (random) function $g_J$:
this gives (for any $\tilde\omega\in\tilde\Omega$),
\[
\hat g_J(\varnothing)^2 = \|g_J\|_2^2 - \sum_{|S|>0} \hat g_J(S)^2 .
\]
Taking the expectation over $\tilde\omega=(\omega,\tau)\in\tilde
\Omega$, this leads to:
\begin{eqnarray*}
{\mathbb E}[{\hat g_J(\varnothing)^2}] & = & {\mathbb
E}[{\| g_J \|_2^2}] - \sum_{|S|>0}
{\mathbb E}[{\hat g_J(S)^2}]  \\
&=& \|g \|_2^2 - \sum_{|S|>0} {\mathbb E}[{\hat g_J(S)^2}]
\qquad\mbox{[by (\ref
{e.hJ})]}  \\
&=& \sum_{|S|=k} \hat g(S)^2 - \sum_{|S|>0} {\mathbb E}[{\hat
g_J(S)^2}]
\qquad\cases{
\mbox{since $g$ is supported}\cr
\mbox{on level-$k$ coefficients}}
\\
&\le& \sum_{|S|=k} {\mathbb E}[{\hat g(S)^2 - \hat
g_J(S)^2}] \qquad\cases{
\mbox{by restricting on} \cr
\mbox{level-$k$ coefficients}.}
\end{eqnarray*}

Now, since $g_J$ is built randomly from $g$ by fixing the variables in
$J=J(\tilde\omega)$, and since $g$ by definition
does not have frequencies larger than $k$, it is clear that
\[
\hat g_J(S) = \cases{
\hat g(S) = \hat f(S), &\quad if $S\cap J(\tilde
\omega) =\varnothing$, \cr
0, &\quad else.}
\]
Therefore, we obtain
\[
\| {\mathbb E}[{g\mid J}] \|_2^2 = {\mathbb E}
[{\hat g_J(\varnothing)^2}] \le\sum_{|S|=k}
\hat g(S)^2 {\mathbb P}[{S\cap J \neq\varnothing}] \le\|
g\|_2^2 k
\reveal
.
\]
Combining with (\ref{e.fg}), this proves Theorem \ref{th.algo}.

\subsection{\texorpdfstring{Is there any better algorithm\textup{?}}{Is there any better algorithm}}

One might wonder whether there exist better algorithms which detect
left-to-right crossings (better in the sense of \textit{smaller revealment}).
The existence of such algorithms would immediately imply sharper
concentration results for the ``Fourier spectrum of percolation'' than in
Corollary~\ref{cor.algo}.

This question of the
``most effective'' algorithm is very natural and has already been
addressed in another paper of Oded et al. \cite{MR2309980}, where
they study random turn hex.
Roughly speaking the ``best'' algorithm might be close to the
following: assume $j$ bits (forming the set $S_j$) have been explored
so far; then
choose for the next bit to be examined, the bit having the highest
chance to be \textit{pivotal} for the left--right crossing
conditionally on what was previously examined ($S_j$). This algorithm
stated in this way does not require additional randomness. Hence one
would need to randomize it
in some way in order to have a chance\vadjust{\goodbreak} to obtain a small revealment. It
is clear that this algorithm (following pivotal locations) in some
sense is
``smarter'' than the one above, but on the other hand
its analysis is highly nontrivial. It is not clear to guess what the
revealment for that algorithm would be.


Before turning to the next section, let us observe that even if we had
at our disposal the most effective algorithm (in terms of
revealment),
it would not yet necessarily imply the expected concentration behavior
of the spectral measure of $f_n$ around its mean (which for an $n\times
n$ box on $\T$ turns out to be of order $n^{3/4}$).
Indeed, in an $n\times n$ box, the algorithm will stop once it has
found a crossing from left to right OR a dual crossing from top to bottom.
In either case, the lattice or dual path is at least of length~$n$;
therefore we always have $|J|=|J_n| \ge n$. But if $\reveal_n$ is the
revealment
of any algorithm computing~$f_n$, then by definition of the revealment,
one has ${\mathbb E}[{|J_n|}] \le O(1) n^2 \reveal_n$
[there are $O(1)n^2$
variables]; since
\mbox{$|J_n|\geq n$}, this implies that $\reveal_n$ is necessarily bigger than
$O(1) n^{-1}$. Now by Corollary \ref{cor.algo}, one has that
\[
\sum_{0<|S|< n^{\alpha}} \hat f_n(S)^2 \le n^{2 \alpha} \reveal
_n .
\]
Therefore, the restriction that $\reveal_n$ has to be bigger than
$O(1)n^{-1}$ implies that with the above algorithmic approach,
one cannot hope to control the Fourier tail of percolation above
$n^{1/2+o(1)}$ (while as we will see in the next section,
the spectral measure of $f_n$ is concentrated around $n^{3/4}$).

The above discussion raises the natural question of finding \textit{good
lower bounds} on the ``revealment of percolation crossings.''
It turns out that 
one can obtain much better lower bounds on the smallest possible
revealment than the straightforward one obtained above.
In our present case (percolation crossings), the best known lower bound
on the revealment follows from the following theorem by O'Donnell and
Servedio.
\begin{theorem}[\cite{MR2341918}]\label{th.OdS}
Let $f \dvtx \Omega_n \to\{-1,1\}$ be a \textup{monotone} Boolean function.
Any randomized algorithm $A$ which computes $f$ satisfies the following:
\[
\reveal_A \geq\frac{\Inf(f)^2}{n} = \frac{(\sum_k \Inf_k(f))^2}
{n} .
\]
\end{theorem}

In our case, $f_n$ depends on $O(n^2)$ variables, and if we are on the
triangular grid~$\T$, the total influence is known to be of order
$n^{3/4+o(1)}$.\vadjust{\goodbreak}
Hence the above theorem implies that
\[
\reveal_n \geq n^{-1/2+o(1)} .
\]
Note that one could have recovered this lower bound also from the case
$k=1$ in Theorem \ref{th.algo}.
Now, using Corollary \ref{cor.algo}, this lower bound shows that
one cannot hope to obtain
concentration results on $\hat f_n$ above level $n^{1/4+o(1)}$. Note
that this is still far from the expected $n^{3/4+o(1)}$.

In fact, Oded (and his coauthors) obtained several results which can be
used to give lower bounds on revealments.
Since these results are related (but slightly \textit{tangential}) to this
survey, we list some of them in this last subsection.

\subsection{\texorpdfstring{Related works of Oded an randomized
algorithms.}{Related works of Oded an randomized algorithms}}

The first related result was obtained in \cite{DecisionTrees}.
Their main theorem can be stated as follows:
\begin{theorem}[\cite{DecisionTrees}]\label{th.DT}
For any Boolean function $f\dvtx \Omega_n \to\{-1,1\}$ and any Markovian
randomized algorithm $A$ computing $f$, one has
\[
\Var[f] \le\sum_{i=1}^n \reveal_i \Inf_i(f) ,
\]
where for each $i\in[n]$, $\reveal_i$ is the probability that the
variable $i$ is examined by~$A$. (In particular, with this notation,
$\reveal_A:= \sup_i \reveal_i$.)
\end{theorem}

This beautiful result can be seen as a strengthening of Poincar\'e's
inequality which states that $\Var[f] \le\sum_i \Inf_i(f)$.
The proof of the latter inequality is straightforward and well known,
and in some sense the proof of the above theorem pays attention
to what is ``lost'' when one derives Poincar\'e's inequality.

This result has deep applications in \textit{complexity theory}. It has
also clear applications in our context since it provides lower bounds
on revealments. For example, together with Theorem \ref{th.OdS}, it
implies the second related result of Oded we wish to mention, the
following theorem from \cite{MR2181623}:
%
\begin{theorem}[\cite{MR2181623}]
Let $f\dvtx \Omega_n \to\{-1,1\}$ be any \textup{monotone} Boolean function; then
any Markovian randomized algorithm $A$ computing $f$ satisfies the following:
\[
\reveal_A \geq\frac{\operatorname{var}[f]^{2/3}} {n^{1/3}} .
\]
\end{theorem}

In \cite{MR2181623} (where other results of this kind are proved),
it is also shown that this theorem is \textit{sharp} up to logarithmic terms.

Since the derivation of this estimate is very simple,
let us see how it follows from Theorems \ref{th.OdS} and \ref{th.DT}.
\begin{eqnarray*}
\Var[f] &\le& \sum_i \reveal_i \Inf_i(f) \qquad\mbox{(note that
monotonicity is not needed here)} \\
&\le& \reveal_A \sum_i \Inf_i(f) \\
&\le& \reveal_A \sqrt{n \reveal_A} \qquad\mbox{(by Theorem \ref{th.OdS})}
\\
&=& \reveal_A^{3/2} \sqrt{n} ,
\end{eqnarray*}
which concludes the proof.
\end{pf*}

Note that for our example of percolation crossings, this implies the
following lower bound on the revealment:
\[
\reveal_n \geq C n^{-2/3} .
\]
This lower bound in this case is not quite as good as the one given by
Theorem~\ref{th.OdS}, but is much better than the
easy lower bound of $\Omega(1/n)$.




\section{\texorpdfstring{The ``geometric'' approach.}{The ``geometric'' approach}}\label{s.geo}

As we explained in the previous section, the randomized algorithm
approach cannot lead to the expected sharp behavior
of the spectrum of percolation. In \cite{GPS}, Pete, Schramm and the
author of this paper
obtain, using a new approach which will be described in this section, a
sharp control of the Fourier spectrum of percolation. This approach
implies among other things
exact bounds on the sensitivity of percolation (the applications of
this approach to dynamical percolation will be explained
in the next section).

\subsection{\texorpdfstring{Rough description of the approach and \textup{``}spectral
measures.\textup{''}}{Rough description of the approach and ``spectral measures''}}

The general idea is to study the ``geometry'' of the frequencies $S$.
For any\vspace*{1pt} Boolean function $f \dvtx \{-1,1\}^n \to\{0,1\}$,
with Fourier expansion $f= \sum_{S\subset[n]} \hat f(S) \chi_S$, one
can consider the different frequencies $S\subset[n]$ as ``random''
subsets of the bits $[n]$; however, they are ``Random'' according to
which measure? Since we are interested in quantities like correlations
\[
{\mathbb E}[{f(\omega) f(\omega^\eps)}] = \sum
_{S\subset[n]} \hat f(S)^2
(1-\eps)^{|S|} ,
\]
it is very natural to introduce the \textit{spectral measure} $\hat\Q=
\hat\Q_f$ on subsets of $[n]$ defined by
%
%
\begin{equation}\label{e.specmeasdef}
\hat\Q[{S}]:= \hat f(S)^2,\qquad S\subset[n] .
\end{equation}
By Parseval, the total mass of the so-defined spectral measure is
\[
\sum_{S\subset[n]} \hat f(S)^2 = \| f \|_2^2 .
\]
One can therefore define a \textit{spectral probability measure} $\hat\PP
= \hat\PP_f$ on the subsets of $[n]$ by
\[
\hat\PP:= \frac1 {\| f \|^2} \hat\Q.\vadjust{\goodbreak}
\]
The random variable under this probability corresponds to the random
frequency and will be denoted by $\Spec=\Spec_f$ (i.e.,
$\hat\PP[{\Spec= S}]:= \hat f(S)^2 / \|f \|^2$). We will
call $\Spec_f$ the
\textit{spectral sample} of $f$.
\begin{remark}
Note that one does not need the Boolean hypothesis here: $\hat\PP_f$
is defined similarly for any real-valued $f \dvtx \Omega_n \to\R$.
\end{remark}

In the remainder of this section, it will be convenient for simplicity
to consider our Boolean functions from $\{-1,1\}^n$
into $\{-1,1\}$ (rather than $\{0,1\}$) thus making $\| f\| ^2 =1$ and
$\hat\Q= \hat\PP$.

Back to our context of percolation, in the following, for all $n\geq
1$, $f_n$ will denote the Boolean function with values in $\{-1,1\}$
corresponding to the left--right crossing in a domain $D_n \subset\Z^2$
(or $\T$) approximating
the square $[0,n]^2$ (or more generally $n\cdot\Omega$ where $\Omega
\subset\C$ is some smooth quad with prescribed boundary arcs).
To these Boolean functions, one associates their spectral samples
$\Spec_{f_n}$.

In \cite{GPS}, we show that most of the spectral mass of $f_n$ is
concentrated around
$n^2 \alpha_4(n)$ (which in the triangular grid $\T$ is known to be of
order~$n^{3/4}$). More exactly
we obtain the following result:
\begin{theorem}[\cite{GPS}]\label{th.GPS}
If $f_n, n\geq1$ are the above indicator functions (in $\{-1,1\}$)
of the left-to-right crossings, then
%
%
\begin{equation}
\limsup_{ n \to\infty} \sum_{0< |S| < \eps n^2 \alpha_4(n)}
\hat
f_n(S)^2 \mathop{\longrightarrow}_{\eps\to0} 0 ,
\end{equation}
or equivalently in terms of the spectral probability measures
%
%
\begin{equation}
\limsup_{n\to\infty} \hat\PP[{0 < |\Spec_{f_n}| < \eps
n^2 \alpha_4(n)}] \mathop{\longrightarrow}_{\eps\to0} 0 .
\end{equation}
\end{theorem}

This result is optimal in localizing where the spectral mass is, since
as we will soon see, it is easy to show that the upper tail
satisfies
\[
\limsup_{n\to\infty} \hat\PP\biggl[{|\Spec_{f_n}| > \frac1
{\eps} n^2 \alpha_4(n)}\biggr] \mathop{\longrightarrow}_{\eps\to
0} 0 .
\]
Hence as hinted long ago
in Figure \ref{f.CrossingEnergy}, we indeed obtain that most of the
spectral mass is localized around $n^2 \alpha_4(n)$ (\mbox{$\approx$}$n^{3/4}$
on $\T$).
It is worth comparing this result with the bound given by the
algorithmic approach in the case of $\T$ which lead to
a spectral mass localized above \mbox{$\approx$}$n^{1/8}$.
We will later give sharp bounds on the behavior of the lower tail
[i.e., at what ``speed'' does it decay to zero below $n^2 \alpha_4(n)$].

Even though we are only interested in the size $|\Spec_{f_n}|$ of the
spectral sample,
the proof of this result will go through a detailed study of the
``geometry'' of~$\Spec_{f_n}$.
The study of this random object for its own sake was suggested by Gil
Kalai and was also considered by Boris Tsirelson
in his study of \textit{black noises}. The general belief, for quite some
time already, was that asymptotically~$\frac1 n \Spec_{f_n}$
should behave like a random Cantor set in $[0,1]^2$.

In some particular cases of Boolean functions $f$, the geometry of
$\Spec_{f}$ can be exactly analyzed: for example, in \cite
{mathPr9903068}, Tsirelson
considers coalescing random walks and follows the trajectory under the
coalescing flow of some particle. The position
of the tagged particle at time $n$ can be represented by a real-valued
function $g_n$ of the array of bits which define the coalescing flow.
Tsirelson shows that the projection
of $\Spec_{g_n}$ on the time axis has the same law as the zero-set of a
random walk. Asymptotically (at least once projected
on the $x$-axis), the spectral sample is indeed looking like a random
Cantor set, and the zeros of a random walk have
a simple structure which enables one to prove a sharp localization of
the spectral mass of $\Spec_{g_n}$ around~$n^{1/2}$.

In the case we are interested in, (i.e., $\Spec_{f_n}$), we do not have
at our disposal such a simple description of $\Spec_{f_n}$.
Until recently, there was some hope that~$\Spec_{f_n}$ would
asymptotically be very similar to the set of pivotal points of~$f_n$,
$\Piv_{f_n}$,
and that one could study the behavior of $|\Spec_{f_n}|$ by focusing on
the geometry of $\Piv_{f_n}$.
We will see that indeed $\Spec_{f_n}$ and $\Piv_{f_n}$ share many common
properties; nevertheless they are conjectured to be asymptotically
\textit{singular} with respect to each other.
This is somewhat ``unfortunate'' since the study of $\Spec_{f_n}$ turns
out to be far less tractable than the study of $\Piv_{f_n}$.

Therefore, since the law of $\Spec_{f_n}$ cannot be ``exactly''
analyzed, the strategy to prove Theorem \ref{th.GPS} will
be to understand some useful properties of the law of $\Spec_{f_n}$
which will then be sufficient to imply the desired
concentration results. Imagine one could show that $\frac1 n \Spec
_{f_n}$ behaves like a Cantor set with a~lot of ``\textit{independence}''
built in (something like a supercritical Galton--Watson tree
represented in $[0,1]^2$ as in Mandelbrot fractal percolation). Then
our desired result would easily follow. In some sense, we will show
that~$\frac1 n \Spec_{f_n}$ indeed behaves like a Cantor set, but we will
only obtain a~very weak independence statement. (It is a posteriori
surprising that such a~weak control on the dependencies
within $\Spec_{f_n}$ could lead to Theorem~\ref{th.GPS}.)

In the next subsections, we will list a few simple observations on
$\Spec_{f_n}$, some of them
illustrating why our spectral sample $\frac1 n \Spec_{f_n}$ should
asymptotically look like a ``random Cantor set''
of $[0,1]^2$.

\subsection{\texorpdfstring{First and second moments on the size of the spectral
sample.}{First and second moments on the size of the spectral sample}}

Let $f \dvtx \{-1,1\}^n \to\R$ be any real-valued function. We have the
following lemma on the
the spectral measure $\hat\Q=\hat\Q_f$ (we state the lemma in
greater generality than for
functions into $\{-1 ,1\}$ since we will need it later).
\begin{lemma}\label{l.A}
For any subset of the bits $A\subset[n]$, one has
%
%
\begin{equation}
\hat\Q[{\Spec_f \subset A}] = {\mathbb E}
[{{\mathbb E}[{f \mid A}]^2}]
= \| {\mathbb E}[{f \mid A}]
\|^2 ,
\end{equation}
where ${\mathbb E}[{f \mid A}]$ is the conditional
expectation of $f$ knowing
the bits in $A$.
\end{lemma}

The proof of this lemma is straightforward and follows from the
definition of $\hat\Q_f$
\begin{eqnarray*}
\hat\Q[{\Spec_f \subset A }] :\!& = & \sum_{S \subset A}
\hat f(S)^2
\\
& = & \biggl\| \sum_{S \subset A} \hat f(S) \chi_S \biggr\|^2 \\
& = & \| {\mathbb E}[{f \mid A}] \|^2 .
\end{eqnarray*}

Now let us return to our particular case of a Boolean function $f \dvtx \{
-1,1\}^n \to\{-1 ,1\}$. Recall that its set of pivotal points
$\Piv=\Piv_f(\omega)$ is the random set of variables $k$ for which
$f(\omega)\neq f(\sigma_k\cdot\omega)$.
\begin{lemma}[(Observation of Gil Kalai)] \label{l.Kalai}
Let $f \dvtx \{-1,1\}^n \to\{-1,1\}$ be a~Boolean function, then
%
%
\begin{equation}
\hat\E[{|\Spec_f|}] = {\mathbb E}[{|\Piv_f|}
] \quad\mbox{and}\quad \hat\E[{|\Spec_f|^2}]
= {\mathbb E}[{|\Piv_f|^2}] .
\end{equation}
\end{lemma}

This lemma is a very useful tool. Its poof is simple, let us sketch how
to derive the first moment; by
definition one has $\hat\E[{|\Spec_f|}] = \sum_{k\in
[n]} \hat\PP[{k \in\Spec_f}]$. Now for $k\in[n]$,
\begin{eqnarray*}
\hat\PP[{k \in\Spec_f}] & = & 1- \hat\PP\bigl[{\Spec
_f \subset[n]\setminus k}\bigr]
\\
& = & \| f \|^2 - \bigl\| {\mathbb E}\bigl[{f \mid [n]\setminus k}
\bigr] \bigr\|^2 \\
&=& \| f - {\mathbb E}[{f \mid \{k\}^c}] \|^2 = \| f
\cdot1 _{ k \in\Piv
_f} \|
^2 \\
&=& {\mathbb P}[{k \in\Piv_f}] .
\end{eqnarray*}
With a similar computation, one can see that for any $k, l \in[n]$,
$\hat\PP[{k,l \in\Spec_f}] = {\mathbb P}[{k,l
\in\Piv_f}]$
thus proving the second moment relation.
\begin{remark}
We mentioned already that $\frac1 n \Spec_{f_n}$ and $\frac1 n \Piv
_{f_n}$ are believed to be asymptotically
singular. This lemma shows that in terms of size, they have the same
first and second moments, but it is not hard to check
that their higher moments differ.\vadjust{\goodbreak}
\end{remark}

\textit{Consequence for the spectral measures of percolation
crossings $f_n$.}
It is a~standard fact in critical percolation (see \cite{GPS}) that
${\mathbb E}[{|\Piv_{f_n}|^2}] \le O(1) {\mathbb E}
[{|\Piv_{f_n}|}]^2$
(this follows from the quasimultiplicativity property). It then follows
from Lemma \ref{l.Kalai} that
\[
\hat\E[{| \Spec_{f_n}| ^2}] \le C \hat\E
[{|\Spec{f_n}|}]^2 .
\]
Hence, using the Paley--Zygmund inequality, there is a universal
constant $c>0$ such that
\[
\hat\PP\bigl[{| \Spec_{f_n}| > c \hat\E[{| \Spec_{f_n}|
}]}\bigr] > c .
\]

Since $\hat\E[{|\Spec_{f_n}|}] = {\mathbb E}[{|
\Piv_{f_n}|}] \asymp n^2 \alpha
_4(n)$, only using simple observations, one already
obtains that at least a ``positive fraction'' of the spectral mass is
localized around $n^2 \alpha_4(n)$.
This property was known for some time already; but in order to be
useful, a kind of 0--1 \textit{law} needs to be proved here.
This turns out to be much more difficult and in some sense this was the
task of~\cite{GPS}.

Notice also that $\hat\E[{|\Spec_{f_n}|}] \asymp n^2
\alpha_4(n) $ easily
implies, by Markov's inequality,
the estimate on the upper tail behavior mentioned above, that is,
\[
\limsup_{n\to\infty} \hat\PP\biggl[{|\Spec_{f_n}| > \frac1
{\eps} n^2 \alpha_4(n)}\biggr] \mathop{\longrightarrow}_{\eps\to
0} 0 .
\]

\subsection{\texorpdfstring{``Scale invariance'' properties of
$\Spec_{f_n}$.}{``Scale invariance'' properties of $\Spec_{f_n}$}}
The goal of this subsection is to identify
some easy properties of $\Spec_{f_n}$ that reveal a certain scale
invariance behavior. This is a first step
toward what we mean by describing the ``geometry'' of $\Spec_{f_n}$.

If one compares with the first two approaches (i.e., hypercontractivity
and algorithmic), both of them gave estimates on the size of
$\Spec_{f_n}$ but did not give any information on the typical shape of
$\Spec_{f_n}$. With these approaches, one cannot make the distinction
between $\Spec_{f_n}$ being typically ``localized'' (all the points
within a Ball of small radius) or on the other hand a spectrum~%
$\Spec_{f_n}$ being ``self-similar.'' We will see, with simple
observations, that at least a~positive fraction of the spectral measure
$\hat\PP_{f_n}$
is supported on self-similar sets.

Before stating these, let us consider a simple analog: let $n:=2^k$ and
let us consider a fractal percolation on
$[0,n]^2$ with parameter $p$. This is defined iteratively in a very
simple way, as a Galton--Watson process: we start at scale $n=2^k$;
assume we are left with a certain set of squares of side-length~$2^l$,
$l\geq1$. If $Q$ is one of these squares of sidelength $2^l$, then
divide this square into four dyadic subsquares and keep these
independently with probability $p$.
Let $M$ denote the random discrete set that remains at the end of the
procedure. By definition on has ${\mathbb E}[{|M|}] = p^k
n^2 = n^{2+\log_2 p}$.
Now, let $1 \ll r \ll n$ be some \textit{mesoscopic} scale; assume $r=
2^l, 1< l < k$. If one is interested in the random set $M$ only
up to scale $r$ (i.e., we do not keep the detailed information below
scale $r$); this corresponds to the set $M_{(r)}$ where we stop the
above induction
once we reach squares of sidelength $r$ (hence $M\subset M_{(r)}$). By
definition one has that ${\mathbb E}[{|M_{(r)}|}] =
(n/r)^{2+\log_2 p}$ (where
$|M_{(r)}|$ denotes
the number of $r$-squares in $M_{(r)}$).

In many ways, $\Spec_{f_n}$ should behave in a similar way as these
fractal percolation sets (which are by construction self-similar).
In order to test the self-similarity of~$\Spec_{f_n}$, let us introduce
some mesoscopic scale $1\ll r \ll n$; one would like to study
the spectral sample $\Spec_{f_n}$ only at the level of detail of the
scale $r$ and then claim that this $r$-``smoothed'' picture
of $\Spec_{f_n}$ should look like $\Spec_{f_{n/r}}$.

Let us precisely define an $r$-smoothed picture of $\Spec_{f_n}$.
Divide the plane
into a lattice of $r\times r$ squares (i.e., $r \Z^2$), and for any
subset $S$ of the window $[0,n]^2$,
we define
\[
S_{(r)}:= \{ \mbox{those $r\times r$ boxes that intersect $S$}\}.
\]

We would like to see now that if $\Spec\sim\hat\PP_{f_n}$, then
$\Spec_{(r)}$ is similar to $\Spec_{f_{n/r}}$.
We have the following estimate on the size of $\Spec_{(r)}$:
\begin{lemma}
\[
{\hat\E}_{f_n} \bigl[ \bigl|\Spec_{(r)}\bigr| \bigr] \asymp\frac{n^2}
{r^2} \alpha_4(r,n)
\asymp
\hat\E[{| \Spec_{f_{n/r}}|}] \approx\biggl( \frac n
r\biggr)^{3/4} ,
\]
where the last approximate equality is known only on $\T$.
\end{lemma}

The proof\vspace*{1pt} is the same as the one for Lemma \ref{l.Kalai}, except that
instead of summing over points, we sum over
squares of radii $r$ in the window $[0,n]^2$. There are $(n/r)^2$ such
squares, and for any such square $Q=Q_r$, it remains to
compute
\[
{\hat\PP}_{f_n} \bigl[ Q\in\Spec_{(r)} \bigr] = {\hat\PP}_{f_n}
[ \Spec\cap Q \neq\varnothing] .
\]
As in Lemma \ref{l.Kalai}, one gets ${\hat\PP}_{f_n} [ \Spec
\cap Q \neq\varnothing] = 1- {\hat\PP}_{f_n} [ \Spec
\subset Q^c ] = \| f - {\mathbb E}[{f \mid Q^c}]
\|^2$.
This expression is obviously bounded from above by the probability that
$Q$ is a pivotal square for $f$ (i.e., that without the information in
$Q$, one cannot
yet predict $f$). This happens only if there is a four-arm event around
$Q$. Neglecting the boundary issues, this is of probability~$\alpha_4(r,n)$.
It can be shown that not only ${\hat\PP}_{f_n} [ \Spec\cap Q
\neq\varnothing]
\le\alpha_4(r,n)$ but in fact ${\hat\PP}_{f_n} [ \Spec\cap Q
\neq\varnothing]
\asymp
\alpha_4(r,n)$, thus concluding the (sketch of) proof.
\end{pf*}

Let us now state two observations about the ``local'' scale invariance
of~$\Spec_{f_n}$.
First, if we consider one of the squares $Q$ of the above lattice of
$r\times r$ squares:
how does $\Spec_{f_n}$ look \textit{inside} $Q$ when one conditions on
$Q\cap\Spec_{f_n} \neq\varnothing$?
This is answered (in the averaged sense) by the following lemma:
\begin{lemma}\label{l.insideQ}
Let $Q$ be an $r\times r$ square in $[0,n]^2$, then
%
%
\begin{equation}\label{e.firstmomS}
{\hat\E}_{f_n} [ | \Spec\cap Q| \mid Q\cap\Spec\neq
\varnothing]
\asymp r^2 \alpha_4(r)\approx r^{3/4} .
\end{equation}
\end{lemma}

The proof works as follows: there are about $r^2$ sites (or edges for
$\Z^2$) in~$Q$, and if $x$ if one of these we have
\begin{eqnarray*}
{\hat\PP}_{f_n} [ x \in\Spec\mid \Spec\cap Q \neq
\varnothing] &=& \frac{
{\hat\PP}_{f_n} [ x\in\Spec] }{ {\hat\PP}_{f_n}
[ \Spec\cap Q \neq\varnothing] } =
\frac
{ {\mathbb P}[{x \in\Piv_{f_n} }] } { {\hat\PP}_{f_n}
[ \Spec\cap Q
\neq\varnothing]
}\\
&\asymp& \frac{\alpha_4(n)} {\alpha_4(r,n)} ,
\end{eqnarray*}
neglecting as usual boundary issues (and using the proof of the
previous lemma for estimating $ {\hat\PP}_{f_n} [ \Spec\cap Q
\neq\varnothing]$). By quasi-multiplicativity, this is indeed of
order $\alpha_4(r)$.

Now, if one now conditions on the whole $\Spec_{f_n}$ to be contained
in one of these mesoscopic box $Q$ (of side-length $r$),
then under this strong conditioning, one expects that inside $Q$,
$\Spec
_{f_n}$ should look like $\Spec_{f_r}$ in an $r\times r$ window.
Using the same techniques one can obtain (the conditioning here
requires some more work).
\begin{lemma}\label{l.strongQ}
Let $Q$ be an $r\times r$-square inside $[0,n]^2$, then
\[
{\hat\E}_{f_n} [ |\Spec| \mid \varnothing\neq\Spec
\subset Q ] \asymp r^2
\alpha_4(r) .
\]
\end{lemma}

All of these observations are relatively easy, since they are just \textit{averaged} estimates instead of almost sure properties.
Notice that one can easily extend these lemmas to their second moment
analog (but this does not yet lead to almost sure properties).

Before turning to the strategy of the proof of Theorem \ref{th.GPS},
let us state a last result in a different vain, which illustrates
the asymptotic Cantor set behavior of $\Spec_{f_n}$.
\begin{theorem}[\cite{GPS}]\label{th.SLS}
In the case of the triangular lattice $\T$, as $n\to\infty$, $\frac1
n \Spec_{f_n}$ converges in law toward a random subset $\Spec$ of $[0,1]^2$
which, if not empty, is almost surely a perfect compact set (i.e.,
without isolated points and in particular infinite).
\end{theorem}

The proof of this result combines ideas from Tsirelson as well as a
``noise''-representation of percolation by Schramm and Smirnov
\cite{SSblacknoise}.
The proof of Theorem \ref{th.GPS} also works in the ``continuous''
setting and implies that the scaling limit of $\frac1 n \Spec_{f_n}$
is (if not empty) an almost sure compact set of Hausdorff dimension $3/4$.

\subsection{\texorpdfstring{Strategy assuming a good control on the dependency
structure.}{Strategy assuming a good control on the dependency structure}}
In this subsection, we will give the basic strategy for proving Theorem
\ref{th.GPS} but without being very careful
concerning the dependency structure within the random set $\Spec
_{f_n}$. In the next subsection, not being able
in \cite{GPS} to obtain a~sharp understanding of this dependency
structure, we will modify the basic strategy accordingly.
Along the explanation of the strategy we will be in good shape to state
a more precise theorem (than Theorem \ref{th.GPS}) on the
lower tail behavior of $|\Spec_{f_n}|$.

The proof of Theorem \ref{th.GPS} will heavily rely on the
self-similarity of $\Spec_{f_n}$.
Before explaining the strategy which will be used for the spectral
sample, let us explain the same strategy applied to
a much simpler example.

Consider a simple random walk on $\Z$, $(X_k)$ starting uniformly on
$\{ -\lfloor{\sqrt{n}}\rfloor,\break \ldots,$ $\lfloor{\sqrt
{n}}\rfloor\}$.
Let $Z_n:= \{ j \in[0,n], X_j=0 \}$ be the zero-set of $(X_k)$ until
time $n$. (We do not start the walk from 0, since as in the case of the
spectral sample in general quads $\Omega$, with positive probability
one can have $Z_n=\varnothing$.) It is clear that ${\mathbb E}[{|
Z_n|}] \asymp n^{1/2}$
as well as ${\mathbb E}[{|Z_n|^2}] \asymp n$. One wishes to
prove a statement for
$|Z_n|$ similar to Theorem \ref{th.GPS} and possibly
with detailed information on the lower tail of the distribution of
$|Z_n|$ (all of this of course is classical and very well known but
this example
we believe helps as an analogy).

Let us consider as above a mesoscopic scale $1 \ll r \ll n$. If $Z=Z_n$
is the zero-set of the random walk, we define as above
$Z_{(r)}$ by dividing the line into intervals of length $r$; $Z_{(r)}$
being the set of these intervals that intersect $Z=Z_n$.
If $J=J_r$ is one of these intervals, then it is clear that
%
%
\begin{equation}\label{e.firstmomZ}
{\mathbb E}[{|Z_n \cap J | \mid Z_n \cap J \neq\varnothing
}] \asymp r^{1/2} .
\end{equation}
This is the analog of one of the (easy) above observations on $\Spec_{f_n}$.
More importantly, the simple structure of the zero-set of the random
walk enables one to obtain very good control
on the dependency structure of $Z_n$ in the following sense: for any
interval $J=J_r$ of length $r$, conditioned on the
whole zero-set $Z_n$ outside of~$J$, one still has good control on the
size of $Z_n$ inside~$J$. More precisely,
there is a universal constant $c>0$ such that
%
%
\begin{equation}\label{e.strongcond}
{\mathbb P}\bigl[{|Z_n \cap J| \geq r^{1/2} \mid J\cap Z_n \neq
\varnothing, Z_n \cap[n]\setminus J}\bigr] \geq c .
\end{equation}
This type of control is much stronger than the estimate \ref
{e.firstmomZ} (or its analog~\ref{e.firstmomS}).
What does it say about the distribution of $|Z_n|$? Each time $Z_n$
intersects an interval $J$ of length $r$, independently
of everything outside, it has a positive chance to be of size $r^{1/2}$
inside. Therefore, if one is interested in the lower tail
estimate ${\mathbb P}[{|Z_n| < r^{1/2}}]$, it seems that
under the conditioning
$\{ |Z_n| < r^{1/2} \}$, the set $Z_n$ will ``typically'' touch
few intervals of length $r$; in other words $|Z_{(r)}|$ will have to be
small. Now, how does the set~$Z_{(r)}$ look when it is conditioned to
be of small
size? Intuitively,
it is clear that conditioned on $Z_{(r)}$ to be very small, say
$|Z_{(r)}|=2$, the two intervals~$J_1$ and $J_2$ intersecting $Z_n$ will
(with high conditional probability) be very close to each other. Indeed
else, one would would have to pay twice
the unlikely event to touch the line and leave it right away. So for
$k$ ``small,'' one expects that
${\mathbb P}[{|Z_{(r)}| = k}] \approx{\mathbb P}
[{|Z_{(r)}| =1}] \asymp(r/n)^{1/2}$
[there are $\frac n r$ $r$-intervals, and each of these has
probability $\asymp(r/n)^{3/2}$ to be
the only interval crossed by $Z_n$, the boundary issues being treated
easily]. Summarizing the discussion, ${\mathbb P}[{|Z_n| <
r^{1/2}}]$ should be
of order
${\mathbb P}[{|Z_{(r)}| \mbox{ is ``small''} }]$ (since
the different intervals
touched by $Z_n$ behave more or less independently), which should be of
order ${\mathbb P}[{|Z_{(r)}|=1}]$;
hence one expects
%
%
\begin{equation}\label{e.LTZ}
{\mathbb P}\bigl[{|Z_n| < \sqrt{r} }\bigr] \leq O(1) \sqrt{\frac
r n} .
\end{equation}

In order to make these heuristics precise, the strong control on the
dependency structure (\ref{e.strongcond}) implies
%
%
\begin{equation}\label{e.expodecayZ}
{\mathbb P}\bigl[{0 < |Z_n| < r^{1/2} \mid \bigl|Z_{(r)}\bigr| = k}\bigr] <
(1-c)^k,
\end{equation}
from which follows
\[
{\mathbb P}[{0 < |Z_n| < r^{1/2} }] \leq\sum_{k\geq
1} {\mathbb P}\bigl[{\bigl|Z_{(r)}\bigr| = k }\bigr]
(1-c)^k .
\]
Thanks to the exponential decay of $(1-c)^k$, it is enough to obtain a
precise understanding of ${\mathbb P}[{|Z_{(r)}| = k}]$
for small values of $k$ and to make sure that this estimate does not
explode exponentially fast for largest values of $k$.
More precisely we need to show that for any $k\geq1$,
%
%
\begin{equation}\label{e.roughZ}
{\mathbb P}\bigl[{\bigl| Z_{(r)}\bigr| = k }\bigr] \leq g(k) {\mathbb
P}\bigl[{\bigl| Z_{(r)}\bigr| = 1}\bigr] \asymp g(k)
\sqrt{\frac r n} ,
\end{equation}
where $g(k)$ is a sub-exponentially fast growing function (which does
not depend on $r$ or $n$). Plugging this estimate into (\ref
{e.expodecayZ}) leads to the precise lower tail
estimate (\ref{e.LTZ}).

The advantage of this strategy is that by introducing the mesoscopic
scale~$r$, and assuming a good dependency structure at that scale,
precise lower tail estimates boil down to rather rough estimates
as (\ref{e.roughZ}) which need to be sharp only for the very bottom
part of the tail-distribution of~$|Z_{(r)}|$.

Needless to say, there are simpler ways to prove the lower tail
estima\-te~(\ref{e.roughZ}) for the zero-set of random walks,
but back to our spectral sample~$\Spec_{f_n}$, this is the way we will
prove Theorem \ref{th.GPS}.

Recall from Lemma \ref{l.insideQ} that if $Q$ is an $r\times r$ square
in $[0,n]^2$, then
%
%
\begin{equation}
{\hat\E}_{f_n} [ | \Spec\cap Q| \mid Q\cap\Spec\neq
\varnothing]
\asymp r^2 \alpha_4(r)\approx r^{3/4} .
\end{equation}
One can prove a second moment analog of this expression which leads to
\[
{\hat\PP}_{f_n} [ |\Spec\cap Q| \geq a r^2 \alpha_4(r)
\mid Q \cap\Spec\neq\varnothing] \geq a
\]
for some absolute constant $a>0$. Since we believe that $\Spec_{f_n}$
behaves asymptotically like a ``nice'' random
Cantor set, it is natural to expect a~nice dependency structure for
$\Spec_{f_n}$ similar to (\ref{e.strongcond}) for $Z_n$.
It turns out that such a precise understanding is very hard to achieve
(and is presently unknown), but let's assume for the remainder
of this subsection that we have such a~nice independence at our disposal.
Hence one ends up with
%
%
\begin{equation}\label{e.expodecayS}
\hat\PP[{0< |\Spec_{f_n}| < a r^2 \alpha_4(r)}] \leq
\sum_{k \geq1}
\hat\PP\bigl[{\bigl|\Spec_{(r)}\bigr| = k}\bigr] (1-c)^k .
\end{equation}
As for the zero-set of random walk, $\Spec_{(r)}$ typically likes to be
spread; hence, if one conditions on
$|\Spec_{(r)}|$ to be small, it should be, with high conditional
probability, very concentrated in space. This hints that for $k$ ``small''
one should have
\[
\hat\PP\bigl[{\bigl|\Spec_{(r)}\bigr| = k}\bigr] \approx\hat\PP\bigl[{\bigl|
\Spec_{(r)}\bigr| = 1}\bigr] \asymp
\frac
{n^2} {r^2} \alpha_4(r,n)^2 .
\]
This estimate follows from the fact that there are $(n/r)^2$ squares
$Q$ of side-length~$r$, and for each of these,
the probability ${\hat\PP}_{f_n} [ \varnothing\neq\Spec\subset
Q ]$ can be
computed using Lemma \ref{l.A}. It is up to constants
the square of the probability that $Q$ is pivotal for $f_n$ [i.e., away
from boundaries it is $\alpha_4(r,n)^2$].

Now to complete the proof (assuming a nice dependency structure), one
needs the following control on the
probability that $|\Spec_{(r)}|$ is very small.
\begin{proposition}[(\cite{GPS}, Section 4)]\label{p.rough}
There is a sub-exponentially fast growing function $g(k), k\geq1$,
such that for any
$1\leq r \leq n $,
%
%
\begin{equation}
{\hat\PP}_{f_n} \bigl[ \bigl| \Spec_{(r)}\bigr | = k \bigr] \leq g(k)
\frac{n^2}{r^2}
\alpha_4(r,n)^2.
\end{equation}
\end{proposition}

The proof of this proposition (which constitutes a nonnegligeable part
of~\cite{GPS}) is of a
``combinatorial'' flavor. It involves an induction over scales which
roughly speaking uses the fact that if $\Spec_{(r)}$
is both spread and of small size, then one can detect several large
empty annuli. The main idea is to ``classify''
the set of all possible $\Spec_{(r)} , |\Spec_{(r)}|=k$ into broader
families, each family corresponding to a certain empty annuli
structure. This classification into classes is easier to handle, since
the families are built in such a~way that they keep
only the meaningful geometric information (i.e., here the large empty
annuli) and in some sense ``average'' over the remaining
``microscopic information.'' The inductive proof reveals that the
families which contribute the most (under $\hat\PP_{f_n}$) are the ones
with fewer empty annuli. (This corresponds to the fact that very small
spectral sets tend to be localized.)
\begin{remark}
Notice that since $g(k) \ll\exp(k)$, then summing (\ref{e.roughS})
with $r:=1$ from $k=1$ to $k=\log n$,
yields to
\[
{\hat\PP}_{f_n} [ 0 < |\Spec| < \log n ] \le\log(n)
g(\log n) n^2
\alpha
_4(n)^2 \ll\log n\, n^{2+\eps} \alpha_4(n)^2\vadjust{\goodbreak}
\]
for any exponent $\eps>0$. Since $n^2 \alpha_4(n)^2 \approx n^{-1/2}$
on $\T$ (and on $\Z^2$, it is known to be $\le n^{-\alpha}$, for some
$\alpha>0$),
Proposition \ref{p.rough} by itself reproves the results from
\cite{MR2001m60016} with better quantitative bounds. This thus gives a
proof more combinatorial
than analytic (hypercontractivity).
\end{remark}

Now Proposition \ref{p.rough} applied to the expected (\ref
{e.expodecayS}) (where we assumed a good dependency knowledge)
leads to following result on the lower tail behavior of $\Spec_{f_n}$.
\begin{theorem}[\cite{GPS}]\label{th.GPSprecise}
Both on the triangular lattice $\T$ and on $\Z^2$, one has
\[
\hat\PP[{0 < | \Spec_{f_n}| < r^2\alpha_4(r) }] \asymp
\frac{n^2} {r^2}
\alpha_4(r,n)^2 ,
\]
where the constants involved in $\asymp$ are absolute constants.
\end{theorem}

This theorem is sharp (up to constants), and it obviously implies the
weaker Theorem \ref{th.GPS}.

In the case of the triangular lattice, where critical
exponents are known, this can be written in a more classical form
for lower tail estimates
\begin{proposition}[\cite{GPS}]
For every $\lambda\in(0,1]$, one has
\[
\limsup_{n\to\infty} \hat\PP\bigl[{0< |\Spec_{f_n}| \leq\lambda
\hat\E[{|\Spec_{f_n}|}] }\bigr] \asymp\lambda
^{2/3} .
\]
\end{proposition}

\subsection{\texorpdfstring{The weak dependency control we achieved and how to modify
the strategy.}{The weak dependency control we achieved and how to modify
the strategy}}

In the previous subsection, we gave a rough sketch of the proof of
Theorem \ref{th.GPSprecise} (and thus of Theorem \ref{th.GPS})
assuming a good knowledge on the dependency structure for $\Spec
_{f_n}$. Indeed it is natural to believe that if
$Q$ is some square of side-length $r$ in $[0,n]^2$ and if $W$ is some
subset of $[0,n]^2$ not too close to $Q$ [e.g., $d(Q,W)\geq r$ is enough],
then for any subsets $\varnothing\neq S \subset Q$
, the following strong independence statement should hold:
\[
{\hat\PP}_{f_n} [ \Spec\cap Q = S \mid \Spec\cap Q \neq
\varnothing]
\asymp
{\hat\PP}_{f_n} [ \Spec\cap Q = S \mid \Spec\cap Q \neq
\varnothing, \Spec\cap W ] .
\]
It is natural to believe such a statement is true since, for example,
it is known to hold for the random set of pivotals points $\Piv_{f_n}$.
Since our goal is to understand the size $|\Spec_{f_n}|$, and since we
expect the spectral sample $\Spec$ to be typically of size $r^2 \alpha_4(r)$
in a square $Q$ satisfying $\Spec\cap Q\neq\varnothing$ (see
Lemma~\ref{l.insideQ}), it would be enough to prove the following slightly weaker
statement:
there is some constant $c>0$ such that 
%
%
\begin{equation}\label{e.expected}
{\hat\PP}_{f_n} [ |\Spec\cap Q| \geq r^2 \alpha_4(r) \mid
\Spec\cap Q \neq\varnothing, \Spec\cap W ] \geq c ,%
\end{equation}
which is the exact analog of (\ref{e.strongcond}). Unfortunately, such
a statement is out of reach at the moment (due to a lack of knowledge
and control on the
law of $\Spec_{f_n}$), but in \cite{GPS} one could prove the following
weaker independence statement.
\begin{proposition}[\cite{GPS}]\label{pr.W}
There is a universal constant $c>0$, such that for any square $Q$ of
side-length $r$ in $[0,n]^2$ and any $W\subset[0,n]^2$ only satisfying
$d(Q,W)\geq r$,
one has
\[
{\hat\PP}_{f_n} [ |\Spec\cap Q| \geq r^2 \alpha_4(r) \mid
\Spec\cap Q \neq\varnothing, \Spec\cap W = \varnothing] \geq
c .
\]
\end{proposition}

This result strengthens the easy estimate from Lemma \ref{l.insideQ}
and constitutes another nonnegligible part of \cite{GPS}
(we will not comment on its proof here).
It says that as far as we looked for $\Spec$ in some region $W$ but did
not find anything there, one can still control the size
of the spectrum $\Spec$ in any disjoint square $Q$ (knowing $\Spec
\cap
Q \neq\varnothing$).

This control on the dependency structure is much weaker than the
expected (\ref{e.expected}), but it turns out that one can modify
the strategy explained in the above subsection in order to deal with
this lack of independence. Let us briefly summarize what is the idea
behind the modified approach.

Let $1\leq r \leq n$, we want to estimate the lower-tail quantity $\PP
[0< |\Spec_{f_n}| < r^2 \alpha_4(r)]$.
As in the definition of $\Spec_{(r)}$, we divide the window $[0,n]^2$
into a~lattice of
$r\times r$ squares. Now, the \textit{naive} idea is to ``scan'' the
window $[0,n]^2$ with $r$-squares $Q$ scanned \textit{one at a time}
[there are $\asymp(n/r)^2$ such squares to scan]. We hope to show that
with high probability, if $\Spec\neq\varnothing$, then at some point
we will find some ``good'' square $Q$, with $|\Spec\cap Q| > r^2\alpha
_4(r)$. The problem being that, if by the time one encounters a good
square, we
discovered some points but only few of them in the region $W$ that was
scanned so far, then due to our weak
dependency control (one keeps a good control only as far as $\Spec\cap
W = \varnothing$), we could not push the scanning procedure any longer.

One way to handle this problem is to scan the window in a \textit{dilute
manner}; that is, at any time, the region $W$ we discovered so far
should be some well-chosen ``sparse'' set. It is not hard to figure out
what should be the right ``intensity'' of this set:
assume we have discovered the region $W$ so far and that we did not
find any spectrum there (i.e., $\Spec\cap W =\varnothing$). We want to continue
the exploration (or scanning) inside some unexplored square $Q$ such
that $d(Q,W)\geq r$. Assume
furthermore that we know $\Spec\cap Q \neq\varnothing$. Then, we are
exactly under a conditional law on $|\Spec\cap Q|$ that we are
able to handle using our Proposition \ref{pr.W}. At that point we
might scan the whole square $Q$ and by Proposition \ref{pr.W}, we have
a positive probability
to succeed [i.e., to find more than $r^2 \alpha_4(r)$ points in the
spectral sample]; but in case we only find few points (say $\log r$ points),
then, as mentioned above, we will not be able to continue the
procedure. Thus, a natural idea is to explore each site (or edge) $x\in
Q$ only
with probability $p_r:= (r^2 \alpha_4(r))^{-1}$ (independently of the
other sites of the square and independently of $\Spec$).
As such one only scans a random subset $\mathcal{Z}$ of $Q$. 
The advantage of this particular choice is that in
the situation where $\varnothing\neq\Spec\cap Q$ is of small
cardinality [compared to $r^2 \alpha_4(r)$], then there is a small chance
that $\mathcal{Z} \cap\Spec\neq\varnothing$, and we will be able to
continue the scanning procedure in a some new square $Q'$
(with $W' = W \cup\mathcal{Z}$). On the other hand, if one discovers
$\Spec\cap\mathcal{Z} \neq\varnothing$,
then with high probability one can guess that $| \Spec\cap Q| $ is of
order $r^2 \alpha_4(r)$, which is what we are looking for.

There are two difficulties one has to face with this modified strategy:

\begin{itemize}
\item The first technicality is that in the above
description, we needed to choose new squares $Q$ satisfying $d(Q,W)\geq r$,
as such we might not be able to scan the whole window $[0,n]^2$, and
we might miss where the spectrum actually lies. There is an
easy way out here: if $Q$ is some new square in the scanning procedure
only satisfying $Q\cap W =\varnothing$, then one still has a
good control on the size of the spectral sample under\vspace*{1pt} the above
conditioning but now restricted
to the concentric square $\bar Q \subset Q$ of side-length $r/2$. More
precisely we have
\[
{\hat\PP}_{f_n} [ |\Spec\cap\bar Q| \geq r^2 \alpha_4(r)
\mid \Spec\cap Q \neq\varnothing, \Spec\cap W = \varnothing
] \geq c .
\]
Hence, we can now explore all the squares $Q$ of the $r\times r$
subgrid of $[0,n]^2$ one at a time and look for the spectrum only within
some dilute sets $\mathcal{Z}$ inside their concentric square $\bar Q$.

\item  The second difficulty is harder to handle. It has
to do with the fact that if one applies the above ``naive'' strategy,
that is,
scanning the squares one by one, then one needs to keep track of the
conditional probabilities ${\hat\PP}_{f_n} [ \Spec\cap Q \neq
\varnothing\mid \Spec\cap W = \varnothing]$ as the explored
set $W$ is growing along
the procedure (indeed in Proposition \ref{pr.W}, we need to condition
on $\Spec\cap Q\neq\varnothing$).
These conditional probabilities turn out to be hard to control and
because of this, exploring squares one at a time is not doable.
Instead, in~\cite{GPS}, we rely on a more abstract or ``averaged'' scanning
procedure where all squares $Q$ are considered in some sense simultaneously
(see \cite{GPS}, ``A~large deviation lemma,'' Section 6).
\end{itemize}

To summarize, the proof of Theorem \ref{th.GPSprecise} (which implies
Theorem \ref{th.GPS}), relies on three main steps:
\begin{longlist}[(3)]
\item[(1)] Obtaining a sharp control on the very beginning of the lower tail
of the spectrum $\Spec_{(r)}$ for all mesoscopic scales
$1 \le r \le n$: this is given by Proposition~\ref{p.rough} (this step
gives an independent proof of Theorem \ref{th.ns}).

\item[(2)] Proving a sufficient control on the dependency structure within
$\Spec_{f_n}$: Proposition \ref{pr.W}.

\item[(3)] Adapting the ``naive'' strategy of scanning mesoscopic squares
$Q$ of size $r$ ``one at a time'' into
an averaged procedure which avoids keeping track of conditional
probabilities like ${\hat\PP}_{f_n} [ \Spec\cap Q \neq
\varnothing\mid \Spec\cap W =\varnothing]$
(Section~6 in~\cite{GPS}).
\end{longlist}

Explaining in detail all these steps would take us too far afield in
this survey, but we hope that this explanation as well as the analogy with
the much simpler case of the zeros of random walks gives some intuition
on how the proof works.

\section{\texorpdfstring{Applications to dynamical percolation.}{Applications to dynamical percolation}}\label{s.dp}

We already mentioned applications to first-passage percolation in the
hypercontractivity section. Let us now
present a very natural model where percolation undergoes a
time-evolution, the model of \textit{dynamical percolation} described below.
The study of the ``dynamical'' behavior of percolation as opposed to
its ``static'' behavior
turns out to be very rich: interesting phenomena arise especially at
the point where the phase transition takes place. We will see that in
some sense,
dynamical \textit{planar} percolation at criticality is a very unstable
(or chaotic) process.
In order to understand this instability, sensitivity of percolation
(and so its Fourier analysis) will play a key role.

\subsection{\texorpdfstring{The model of dynamical percolation.}{The model of dynamical percolation}} 
As mentioned earlier in the \hyperref[sec1]{Introduction} of this survey, \textit{dynamical
percolation} was introduced in \cite{MR1465800}
as a~natural dynamic illustrating \textit{thermal agitation} (similar
to
\textit{Glauber dynamics} for the Ising model).

The dynamic introduced in \cite{MR1465800} (also invented
independently by
Itai Benjamini) is essentially the unique dynamic which preserves the
i.i.d. law of percolation and which is Markovian.

Due to the i.i.d. structure of the percolation model, the dynamic is
very simple: on $\Z^d$, at parameter $p\in[0,1]$,
each edge $e\in\E^d$ is updated at rate one independently of the other
edges. This means that at rate one, independently of everything else,
one keeps the edge with probability $p$, and removes it with
probability $1-p$. See \cite{MR1465800,SchrammSteif} and especially
\cite{SurveySteif}
for background on this model.

In \cite{MR1465800}, it is shown that no surprises arise outside of the
critical point~$p_c(\Z^d)$ in dimension $d\geq2$.
If $(\omega^p_t)_{t\geq0}$ denotes a trajectory
of percolation configurations of parameter $p$ on $\Z^d$ evolving as
described above, they prove that if $p>p_c$, then almost surely,
there is an infinite cluster in ALL configurations~$\omega^p_t$,
$t\geq0$
(as well, in the subcritical regime $p<p_c$, a.s. clusters remain
finite along the dynamic).

What happens at the critical point $p_c(\Z^d)$ requires more care: in
\cite{MR1465800}, H{\"a}ggstr{\"o}m, Peres and Steif prove that if the
dimension is high enough
($d\geq19$), then as in the subcritical regime, the clusters of
$\omega
_t^{p_c}$ remain finite as the time $t$ goes on.
Their proof relies essentially on the ``mean-field'' behavior of
percolation rigorously obtained by Hara and Slade \cite{MR1043524}
in dimension $d\geq19$
(conjectured to hold for $d>6$ and $d=6$ with ``logarithmic'' corrections).

In lower dimensions, the phase transition is somehow more ``abrupt.''
This can be seen by the fact that the density of the infinite cluster
represented by $p\mapsto\theta(p)$ has an infinite derivative at the
critical point (see \cite{MR893131} in $d=2$). Because of
this phenomenon, one cannot easily
rule out the sudden appearance of infinite clusters along the dynamic;
consequently the study of dynamical percolation needs a more detailed study.
We will restrict ourselves to the case of dimension $d=2$, since in
dimension $d=3$ say, very little is known
(even the static behavior of $3d$-percolation at the critical point is
unknown, which makes its dynamical study at this point hopeless).


In the planar case $d=2$, where there is a very rich literature on the
critical behavior, it was asked back in \cite{MR1465800} whether at the
critical point
(say on $\Z^2$ at $p_c=1/2$), such exceptional infinite clusters appear
along the dynamic or not.
More exactly, is it the case that there exists almost surely a (random)
subset $\varnothing\neq\Exc\subset\R$ of \textit{exceptional times},
such that for any
$t\in\Exc$, there is an infinite cluster in $\omega_t^{p_c}$?
We will see in the present section that there indeed exists such
exceptional times.
This illustrates the above mentioned unstable character of dynamical
percolation at criticality.
Note that from the continuity of the phase transition in $d=2$, it
easily follows
that the exceptional set $\Exc$ is almost surely of Lebesgue measure zero.

This appearance of exceptional infinite clusters along the dynamic will
be guaranteed by the fact\vspace*{1pt} that large-scale connectivity properties
should decorrelate very fast. Roughly speaking, if the large-scale
geometry of the configurations $\omega_t^{p_c}$ ``changes'' very
quickly as
time goes on, then it will have better chances to ``catch'' infinite
clusters. As we argued earlier,
the initial study of noise sensitivity in \cite{MR2001m60016} was
partially motivated
by this problem since the rapid decorrelation needed to prove the
existence of exceptional times corresponds exactly to their notion
of noise sensitivity. The results obtained in \cite{MR2001m60016} already
illustrated the fact that at large scales, connectivity properties of
$(\omega_t)_{t\geq0}$
evolve very quickly. Nonetheless, their ``logarithmic'' control of noise
sensitivity was not sufficient to imply the appearance of infinite clusters.
The first proof of the above conjecture (about the existence of
exceptional times) in the context of a two-dimensional lattice was achieved
by Schramm and Steif in \cite{SchrammSteif} for the triangular lattice
$\T$. Their proof relies on their algorithmic approach of the ``Fourier
spectrum''
of percolation explained in Section~\ref{s.algo}. Before explaining in
more detail the results from~\cite{SchrammSteif} and \cite{GPS}, let
us present the strategy to prove the existence of exceptional times as
well as a variant model of dynamical percolation in $d=2$ introduced by
Benjamini and Schramm where the decorrelation structure is somehow
easier to handle (does not require any Fourier analysis).

\subsection{\texorpdfstring{Second moment analysis and the ``exceptional planes of
percolation.''}{Second moment analysis and the ``exceptional planes of
percolation''}}

Suppose we are considering some model of percolation at criticality,
for example, $\Z^2$ or the triangular lattice $\T$ at $p_c=1/2$,
but it could also be Boolean percolation (a Poisson point process of
balls in the plane) at its critical intensity.
Let's consider a trajectory of percolation configurations $(\omega
_t)_{t\geq0}$ for the model under consideration evolving as above
(in the case of the Boolean model, one could think of many natural
dynamics, for example, each ball moving independently according to a Brownian
motion). One would like to prove the existence of \textit{exceptional
times} $t$ for which there is an infinite cluster in $\omega_t$.
It turns out to be more convenient to show the existence of exceptional times
where the origin itself is connected to infinity via some open path;
furthermore it is enough to show that with positive probability
there is some time $t\in[0,1]$ where $0 \stackrel{\omega
_t}{\longleftrightarrow}\infty$.

For any large radius $R\gg1$, let us then introduce $Q_R$ to be the
set of times where
the origin $0$ is connected to distance $R$
\[
Q_R:= \bigl\{t \in[0,1] \dvtx 0\stackrel{\omega_t}{\longleftrightarrow}R\bigr\}.
\]

Proving the existence of exceptional times boils down to proving that
with positive probability $\bigcap_{R>0} Q_R \neq\varnothing$.
Even though the sets $Q_R$ are not closed, with some additional
technicality (see \cite{SchrammSteif}),
it is enough to prove that $\inf_{R>1} {\mathbb P}[{Q_R \neq
\varnothing}] > 0$.

In order to achieve this, we introduce the random variable $X_R$
corresponding to the amount of time where 0 is connected to distance $R$
\[
X_R:= \int_0^1 1_{0\stackrel{\omega_t}{\longleftrightarrow}R} \,dt.
\]
Our goal is therefore to prove that there is some universal constant
$c>0$, so that ${\mathbb P}[{X_R>0}]>c$ for any radius $R$.
As usual, this is achieved by the second moment method,
\[
{\mathbb P}[{Q_R \neq\varnothing}] = {\mathbb P}[{X_R
>0}] \geq\frac{{\mathbb E}[{X_R}]^2}{{\mathbb
E}[{X_R^2}]},
\]
and it remains to prove the existence of some
constant $C>0$ such that for all $R>1$, ${\mathbb E}[{X_R^2}
] < C {\mathbb E}[{X_R}]^2$.
Now, notice that the second moment can be written
%
%
\begin{eqnarray}\label{e.XR2}
{\mathbb E}[{X_R^2}]& = &\mathop{\iint_{0\leq s \leq
1}}_{\hspace*{12pt}0 \leq t \leq1} {\mathbb P}\bigl[{0\stackrel{\omega
_s}{\longleftrightarrow}R, 0\stackrel{\omega
_t}{\longleftrightarrow}R}\bigr] \,ds\,dt\nonumber\\[-8pt]\\[-8pt]
& \leq& 2 \int_0^1 {\mathbb P}\bigl[{0\stackrel{\omega
_0}{\longleftrightarrow}R, 0\stackrel{\omega_t}{\longleftrightarrow
}R}\bigr] \,dt .\nonumber
\end{eqnarray}

One can read from this expression that a fast decorrelation over time
of the event $0\stackrel{\omega}{\longleftrightarrow}R$
will imply the existence of exceptional times. The study of the
correlation structure of these events usually goes through Fourier analysis
unless there is a more direct way to get a hand on these (we will give
such an example of dynamical percolation below).
Fourier analysis helps as follows: if $f_R$ denotes the indicator
function of the event $\{ 0\stackrel{\omega}{\longleftrightarrow}R \}$
then, as a~Boolean function (of the bits inside the disk of radius $R$
into $\{0,1\}$), it has a~Fourier--Walsh expansion.
Now, as was explained in Section~\ref{ss.fourier}, the desired
correlations are simply expressed in terms of the Fourier coefficients
of $f_R$\looseness=-1
\begin{eqnarray*}
&&{\mathbb P}\bigl[{0\stackrel{\omega_0}{\longleftrightarrow}R,
0\stackrel{\omega_t}{\longleftrightarrow}R}\bigr]\\
&&\qquad = {\mathbb
E}[{f_R(\omega_0) f_R(\omega_t)}]
\\
&&\qquad= {\mathbb E}\biggl[{\biggl(\sum_{S\subset B(0,R)} \hat{f_R}(S)
\chi_S(\omega_0) \biggr) \biggl( \sum_{S\subset B(0,R)} \hat
{f_R}(S) \chi_S(\omega_t) \biggr)}\biggr]
\\
&&\qquad= {\mathbb E}[{f_R}]^2 + \sum_{\varnothing\not\subset S
\subset B(0,R)} \hat
{f_R}(S)^2 \exp(-t|S|) . 
\end{eqnarray*}\looseness=0
%

We see in this expression, that for small times $t$, the frequencies
contributing in the correlation between $\{ 0\stackrel{\omega
_0}{\longleftrightarrow}R \}$
and $\{ 0\stackrel{\omega_t}{\longleftrightarrow}R \}$ are of ``small''
size $|S| \lesssim1/t$. In particular, in order to detect the
existence of exceptional times,
one needs to achieve good control on the \textit{lower tail} of the \textit{Fourier spectrum} of $f_R$.

We will give some more details on this setup in the
coming subsection, but before going more deeply into the Fourier side,
let us briefly present the first planar percolation model for which
this phenomenon
of exceptional infinite clusters was proved: \textit{exceptional planes of
percolation} \cite{MR99i60173}.

The model can be described as follows: sample balls of radius $r>0$ in
$\R^4$ according to a Poisson point process with Lebesgue intensity
1 on $\R^4$. Let~$\mathcal{U}_r$ denote the union of these balls. Now,
for each plane $P\subset\R^4$,
let $\omega_P^r= \omega_P$ be the percolation configuration
corresponding to the set $\mathcal{U}_r \cap P$. For a fixed plane $P$,
$\omega_P^r$
follows the law of
a Boolean percolation model with random i.i.d. radii (in $[0,r]$).

This provides us with a natural planar dynamical percolation where the
dynamic is here parametrized by planes $P\subset\R^4$
instead of time.
It is\break known~\cite{MR1409145} that there is a critical parameter
$r_c>0$ so that the corresponding planar model is critical (with a.s.
no infinite cluster at the critical radius $r_c$). This means that at
criticality, for ``almost all'' planes $P\subset\R^4$, there is no
infinite cluster in $\omega_P$. Whence the natural question:
are there are exceptional planes of percolation or not?
\begin{theorem}[(Benjamini and Schramm \cite{MR99i60173})]
At the critical radius $r_c$, there exist exceptional planes $P\in\R
^4$ for which
$\omega_P:= \mathcal{U}_{r_c} \cap P$ has an infinite cluster.
\end{theorem}

The proof of this result follows the above setup. To prove the
existence of such exceptional planes, it will be enough to restrict
ourselves to the subclass of planes $P$ which include the origin, and
we will ask that the origin itself is connected to infinity within
$\omega_P$. If $P_0$ is some plane of reference,
by the second moment method, the key quantity to estimate is the
correlation ${\mathbb P}[{0\stackrel{\omega_{P_0}}
{\longleftrightarrow}R, 0\stackrel{\omega_{P}}
{\longleftrightarrow} R }]$, this correlation
being then integrated uniformly over planes $P$ containing the origin.

The big advantage of this model is that if $P\neq P_0$, since the
radius $r_c$ of the balls in $\R^4$ is finite,
percolation configurations $\omega_P$ and $\omega_{P_0}$ will
``interact'' together
only within a finite radius $R_{P_0,P}$ which depends on the angle~%
$\theta_{P_0,P}$ formed by the two planes
(the closer they are, the further away they interact). We can therefore
bound the correlation as follows (using similar notations as in
Section \ref{ss.perco}):
\begin{eqnarray*}
{\mathbb P}\bigl[{0\stackrel{\omega_{P_0}} {\longleftrightarrow
}R, 0\stackrel{\omega_{P}} {\longleftrightarrow} R }\bigr] &\le&
{\mathbb P}\bigl[{0\stackrel{\omega_{P_0}} {\longleftrightarrow
}R_{P_0,P}, R_{P_0, P}\stackrel{\omega_{P_0}} {\longleftrightarrow
} R, R_{P_0, P}\stackrel{\omega_{P}} {\longleftrightarrow} R}
\bigr] \\
&\le& \alpha_1(R_{P_0,P}) \alpha_1(R_{P_0,P}, R)^2
\end{eqnarray*}
by independence and since configurations do not interact after
$R_{P_0,P}$. 
Now, if we had at our disposal a RSW theorem similar as what is known
for $\Z^2$ percolation, this would imply a quasimultiplicativity result
on the one-arm event
\[
\alpha_1(r_1,r_3) \asymp\alpha_1(r_1,r_2) \alpha_1(r_2,r_3)
\]
for all $r_1\le r_2 \le r_3$. It would then follow that the correlation
is bounded by $O(1) \alpha_1(R_{P_0,P})^{-1} \alpha_1(R)^2$
(see \cite{MR99i60173} where they relied on a more complicated argument
since RSW was not available).

In the second moment argument, we need to bound the second moment by
the squared of the first one [here $\alpha_1(R)^2$]. Hence
it remains to show that $\int_P \alpha_1(R_{P_0,P})^{-1} \,d\lambda(P)
< \infty$. This can be done using the expression of $R_{P_0,P}$ in terms
of the angle $\theta_{P_0,P}$ as well as some control on the decay of
$\alpha_1(u)$ as $u\to\infty$. See \cite{MR99i60173} where the computation
is carried out.

\subsection{\texorpdfstring{Dynamical percolation, later
contributions.}{Dynamical percolation, later contributions}}
In this subsection, we present the results from \cite{SchrammSteif} and
\cite{GPS} on dynamical percolation.

\subsubsection{\texorpdfstring{Consequences of the algorithmic approach
\cite{SchrammSteif}.}{Consequences of the algorithmic approach [40]}}
As we described in Section \ref{s.algo}, Schramm and Steif introduced a
powerful technique in \cite{SchrammSteif} to control
the Fourier coefficients of general Boolean functions $f \dvtx \{-1,1\}^n
\to\{0,1\}$. We have already seen in Corollary \ref{cor.algo}
that their approach enabled one to obtain a polynomial control on the
sensitivity of percolation crossings, thus strengthening previous results
from \cite{MR2001m60016}. But more importantly, using their technique (mainly
Theorem \ref{th.algo}) they obtained deep results
on dynamical percolation on the triangular grid $\T$.
Their main result can be summarized in the following theorem (see
\cite{SchrammSteif} for their other contributions).
\begin{theorem}[(Schramm and Steif \cite{SchrammSteif})]\label{th.SS}
For dynamical percolation on the triangular lattice $\T$ at the
critical point $p_c=1/2$,
one has:
\begin{itemize}
\item Almost surely, there exist exceptional times $t\in[0,\infty]$
such that $\omega_t$ has an infinite cluster.

\item Furthermore, if $\Exc\subset[0,\infty]$ denotes the (random)
set of these exceptional times, then the Hausdorff dimension of $\Exc$
is an almost sure constant in $[1/6, 31/36]$.
\end{itemize}
\end{theorem}

They conjectured that the dimension of these exceptional times is a.s.
$31/36$. Let us now explain how they obtained such a result.

As we outlined above, in order to prove the existence of exceptional
times, one needs to prove fast decorrelation of
$f_R(\omega):= 1_{\{ 0\stackrel{\omega} {\longleftrightarrow}R\}} $.
Recall that the correlation between the configuration at time 0,
$\omega_0$
and the configuration at time $t$, $\omega_t$ can be written
%
%
\begin{eqnarray}\label{e.CorrelQ}
{\mathbb P}\bigl[{0\stackrel{\omega_0}{\longleftrightarrow}R,
0\stackrel{\omega_t}{\longleftrightarrow}R}\bigr]& = &{\mathbb
E}[{f_R(\omega_0) f_R(\omega_t)}]
\nonumber\\
&=& {\mathbb E}[{f_R}]^2 + \sum_{\varnothing\not\subset S
\subset B(0,R)} \hat
{f_R}(S)^2 \exp(-t|S|) \\
&\lesssim& {\mathbb E}[{f_R}]^2 + O(1) \sum_{k=1}^{1/t}
\sum_{|S|=k} \hat
f_R(S)^2 .\nonumber
\end{eqnarray}
The correlation for small times $t$ is thus controlled by the lower
tail of the Fourier spectrum of $f_R$.
The great advantage of the breakthrough approach in \cite{SchrammSteif}
is that their technique is particularly well
suited to the study of the lower tail of $\hat f_R$.
Indeed Theorem \ref{th.algo} says that if one can find a Markovian
algorithm computing $f_R$ with small
revealment $\reveal$, then for any $l\geq1$, one has
\[
\sum_{k=1}^l \sum_{|S|=k} \hat f_R(S)^2 \le l^2 \reveal\| f_R \|
_2^2 = l^2 \reveal \alpha_1(R) ,
\]
since by definition of $f_R$, $\| f_R \|_2^2 = {\mathbb P}[{0
\stackrel{\omega}{\longleftrightarrow} R}] =\alpha_1(R)$
(which in the triangular
lattice is of order $R^{-5/48}$).
What remains to be done is to find an algorithm minimizing the
revealment as much as possible. But there is a difficulty here, similar
to the one
encountered in the proof of the sub-Gaussian fluctuations of FPP in $\Z
^2$: any algorithm computing $f_R$ will have to examine sites
around the origin $0$ making its revealment close to one. This is not
only a technical difficulty; the deeper reason
comes from the following interpretation of the correlation ${\mathbb
P}[{0\stackrel{\omega_0}{\longleftrightarrow}R, 0\stackrel
{\omega_t}{\longleftrightarrow}R }]$:
this correlation can be seen as ${\mathbb P}[{0\stackrel{\omega
_0}{\longleftrightarrow}R }] {\mathbb P}[{0\stackrel
{\omega_t}{\longleftrightarrow}R \mid 0\stackrel{\omega
_0}{\longleftrightarrow}R }]$. Now, assume $R$
is very large and $t$ very small; if one conditions on the event
$\{ 0\stackrel{\omega_0}{\longleftrightarrow}R \}$, since few sites are
updated, the open path in $\omega_0$ from $0$ to distance $R$
will still be preserved in~$\omega_t$ at least up to some distance
$L(t)$ while further away, large scale connections start being ``noise
sensitive.''
In some sense the geometry associated to the event $\{ 0\stackrel
{\omega
}{\longleftrightarrow}R \}$ is ``frozen'' on a certain
scale between time $0$ and time $t$. Therefore, it is natural to divide
our correlation analysis into two scales: the ball of radius $r=r(t)$
and the annulus
from $r(t)$ to $R$ [we might need to take $r(t) \gg L(t)$ since the
control on the Fourier tail given by the algorithmic approach is not sharp].
Obviously the ``frozen radius'' $r=r(t)$ increases as $t\to0$. As in
the exceptional planes case,
let us then bound our correlation by
%
%
\begin{eqnarray}
{\mathbb P}\bigl[{0\stackrel{\omega_0}{\longleftrightarrow}R,
0\stackrel{\omega_t}{\longleftrightarrow}R}\bigr]& \le& {\mathbb
P}\bigl[{0\stackrel{\omega_0}{\longleftrightarrow} r}\bigr]
{\mathbb P}\bigl[{r \stackrel{\omega_0}{\longleftrightarrow}R, r
\stackrel{\omega_t}{\longleftrightarrow}R}\bigr]\nonumber\\[-8pt]\\[-8pt]
&\le& \alpha_1(r) {\mathbb E}[{f_{r,R}(\omega_0) f_{r,R}
(\omega_t)}],\nonumber
\end{eqnarray}
where $f_{r,R}$ is the indicator function of the event $r \stackrel
{\omega}{\longleftrightarrow}R$. Now, as above
\begin{eqnarray*}
{\mathbb E}[{f_{r,R}(\omega_0) f_{r,R} (\omega_t)}]
&\lesssim& {\mathbb E}[{f_{r,R}}]^2 +
O(1)\sum_{k=1}^{1/t} \sum_{|S|=k} \hat f_{r,R}(S)^2 \\
& \lesssim& \alpha_1(r,R)^2 + O(1) t^{-2} \reveal_{f_{r,R}}
\alpha_1(r,R) .
\end{eqnarray*}

The Boolean function $f_{r,R}$ somehow avoids the singularity at the
origin, and it is possible to find algorithms for this function
with small revealments. It is an interesting exercise to adapt the
exploration algorithm which was used above in the left--right crossing case
to our current radial case. We will sketch at the end of this
subsection, why using a natural exploration-type algorithm, one can obtain
a revealment for $f_{r,R}$ of order $\reveal\asymp\alpha_2(r) \alpha
_1(r,R)$. Assuming this, one ends up with
\begin{eqnarray*}
{\mathbb P}\bigl[{0\stackrel{\omega_0}{\longleftrightarrow}R,
0\stackrel{\omega_t}{\longleftrightarrow}R}\bigr]& \lesssim&
\alpha_1(r) \bigl( \alpha
_1(r,R)^2 + t^{-2} \alpha_2(r)\alpha_1(r,R)^2 \bigr) \\
&\lesssim& \biggl( \alpha_1(r)^{-1} + t^{-2} \frac{\alpha
_2(r)}{\alpha
_1(r)} \biggr) \alpha_1(R)^2 ,
\end{eqnarray*}
using the quasi-multiplicativity property \ref{pr.quasi}. Now using the
knowledge of the critical exponents on $\T$ without
being rigorous (these exponents are known only up to logarithmic
corrections), one gets
\[
{\mathbb P}\bigl[{0\stackrel{\omega_0}{\longleftrightarrow}R,
0\stackrel{\omega_t}{\longleftrightarrow}R}\bigr] \lesssim(1 +
t^{-2} r^{-1/4}) r^{5/48}
\alpha_1(R)^2 .
\]

Optimizing in $r=r(t)$, one chooses $r(t):= t^{-8}$ which implies the
following bound on the correlation:
%
%
\begin{equation}\label{e.correlalgo}
{\mathbb P}\bigl[{0\stackrel{\omega_0}{\longleftrightarrow}R,
0\stackrel{\omega_t}{\longleftrightarrow}R}\bigr] \lesssim
t^{-5/6} \alpha_1(R)^2 .
\end{equation}
Now, since $\int_0^1 t^{-5/6} \,dt <\infty$, by integrating the above
correlation over the unit interval, one has from (\ref{e.XR2})
that ${\mathbb E}[{X_R^2}] \le C {\mathbb E}
[{X_R}]^2$ for some universal $C>0$ (since, by
definition of $X_R$, ${\mathbb E}[{X_R}]=\alpha_1(R)$).
This thus proves the existence of exceptional times.

If one had obtained a much weaker control on the correlation than that
in (\ref{e.correlalgo}), for example, a bound of
$t^{-1}\log(t)^{-2} \alpha_1(R)^2$, this would have still implied the
existence of exceptional times: one can thus exploit
more the good estimate provided by (\ref{e.correlalgo}). This
estimate in fact easily implies the second part of Theorem
\ref{th.SS}, that is, that almost surely the Hausdorff dimension of the
set of exceptional times is greater than $1/6$. (The upper bound of
$31/36$ is rather easy to obtain; it is a first moment analysis and
follows from the behavior of the density function $\theta$
near $p_c=1/2$ (see \cite{SchrammSteif} and \cite{SurveySteif}).) The
proof of the lower bound on the Hausdorff dimension which is based on
the estimate (\ref{e.correlalgo}),
is classical and essentially consists into defining a (random) Frostman
measure on the set $\Exc$ of exceptional times; one concludes
by bounding its expected $\alpha$-energies for all $\alpha<1/6$. 

See \cite{SchrammSteif} for rigorous proofs (taking care of the
logarithmic corrections etc.).

Let us conclude this subsection by briefly explaining why one can
achieve a revealment of order $\alpha_2(r) \alpha_1(r,R)$
for the Boolean function $f_{r,R}(\omega):=1_{\{ r \stackrel{\omega
}{\longleftrightarrow}R\}}$. We use an algorithm that mimics the
one we used in the ``chordal'' case except the present setup is
``radial.'' As in the chordal case, we randomize the starting point of
our exploration
process: let's start from a site taken uniformly on the ``circle'' of
radius $R$. Then, let's explore the picture with an exploration path
$\gamma$ directed
toward the origin; this means that as in the chordal case, when the
interface encounters an open (resp., closed) site, it turns say on the
right (resp., left),
the only difference being that when the exploration path closes a loop
around the origin, it continues its exploration inside the connected
component of the
origin (see \cite{07100856} for more details on the radial exploration
path). It is known that this discrete curve converges toward \textit{radial} $\mathrm{SLE}_6$
on $\T$, when the mesh goes to zero. It turns out that the so-defined
exploration path gives all the information we need. Indeed, if the
exploration path closes
a clockwise loop around the origin, this means that there is a closed
circuit around the origin making $f_{r,R}$ equal to zero. On the other hand,
if the exploration path does not close any clockwise loop until it
reaches radius $r$, it means that $f_{r,R}=1$. Hence, we run the
exploration path until
either it closes a clockwise loop or it reaches radius $r$. Now, what
is the revealment for this algorithm? Neglecting boundary issues
(points near
radius $r$ or $R$), if $x$ is a point at distance $u$ from 0, with $2r
< u < R/2$, in order for $x$ to be examined by the algorithm, it is
needed that the exploration
path did not close any clockwise loop before radius $u$; hence there is
an open path from $u$ to $R$, and this already costs a factor $\alpha_1(u,R)$.
Now at the level of the scale $u$, what is the probability that $x$
lies on the interface? The interface is asymptotically of dimension
$7/4$, so in the
annulus $A(u/2, 2u)$, about $u^{7/4+o(1)}$ points lie on $\gamma$ and,
the probability that the point $x$ is in~$\gamma$ is of order
$u^{-1/4}$. In fact it is exactly
up to constant $\alpha_2(u)$. Hence, the probability ${\mathbb P}
[{x\in J}]$, for
$|x|=u$ is of order $\alpha_2(u) \alpha_1(u,R)$. Now recall that the
revealment is the supremum
over all sites of this probability; it is easy to see that the smaller
the scale ($r$), the higher this probability is; therefore neglecting
boundary issues, one obtains
that $\reveal\approx\alpha_2(r) \alpha_1(r,R)$.

\vspace*{-2pt}\subsubsection{\texorpdfstring{Consequences of the Geometric approach
\cite{GPS}.}{Consequences of the Geometric approach [12]}}
We now pres\-ent the applications to dynamical percolation of the sharp
estimates on the Fourier spectrum obtained in \cite{GPS}.
Since the results in \cite{GPS} are sharp on the triangular lattice
$\T
$ as well as on $\Z^2$, the first notable consequence
is the following result:\vspace*{-2pt}

\begin{theorem}[\cite{GPS}]
If $(\omega_t)_{t\geq0}$ denotes some trajectory of dynamical
percolation on the square grid $\Z^2$ at $p_c=1/2$,
then almost surely, there exist exceptional times $t\in[0,\infty)$,
such that $\omega_t$ has an infinite cluster.

Furthermore there is some $\alpha>0$, such that if $\Exc\subset
[0,\infty)$ denotes the random set of these exceptional times, then
almost surely, the Hausdorff dimension of $\Exc$ is greater than
$\alpha$.\vspace*{-2pt}
\end{theorem}
\begin{remark}
The paper \cite{SchrammSteif} came quite close to proving that
exceptional times exist
for dynamical critical bond percolation on $\Z^2$, but there is still a
missing gap if one wants to use
randomized algorithm techniques in this case.\vspace*{-2pt}
\end{remark}

On the triangular lattice $\T$, since the critical exponents are known,
the sharp control on $\hat f_R$ in \cite{GPS} and especially on its
\textit{lower tail} enables one to obtain detailed information on the
structure of the (random) set of exceptional times $\Exc$. One can
prove, for instance.\vspace*{-2pt}
\begin{theorem}[\cite{GPS}]\label{th.31}
If $\Exc\subset[0,\infty)$ denotes the (random) set of exceptional
times of dynamical percolation on $\T$, then almost surely the
Hausdorff dimension
of $\Exc$ equals $\frac{31}{36}$.\vspace*{-2pt}
\end{theorem}

This theorem strengthens Theorem \ref{th.SS} from \cite{SchrammSteif}.
Finally, one obtains the following interesting phenomenon:\vspace*{-2pt}
\begin{theorem}[\cite{GPS}]
Almost surely on $\T$, there exist exceptional times ``of the second
kind'' $t\in[0,\infty)$, for which an infinite
(open) cluster \textup{coexists} in $\omega_t$ with an infinite dual
(closed) cluster.

Furthermore, if $\Exc^{(2)}$ denotes the set of these exceptional
times; then almost surely, the Hausdorff dimension
of $\Exc^{(2)}$ is greater than $1/9$.\vspace*{-2pt}
\end{theorem}

We conjectured in \cite{GPS} that $\Exc^{(2)}$ should be almost surely
of dimension $2/3$;
but unfortunately, the methods in \cite{GPS}\vadjust{\goodbreak} do not apply well to
nonmonotone functions, in particular in this case, to
the indicator function of a two-arm event.

Let us now focus on Theorem \ref{th.31}. As in in the previous
subsection, in order to have detailed information on the dimension
of the exceptional set $\Exc$, one needs to obtain a sharp control on
the correlations ${\mathbb P}[{0\stackrel{\omega
_0}{\longleftrightarrow}R, 0\stackrel{\omega
_t}{\longleftrightarrow}R}]$. More exactly, the almost
sure dimension $31/36$ would follow from the following estimate:
%
%
\begin{equation}\label{e.correlgeom}
{\mathbb P}\bigl[{0\stackrel{\omega_0}{\longleftrightarrow}R,
0\stackrel{\omega_t}{\longleftrightarrow}R}\bigr] \lesssim
t^{-5/36} \alpha_1(R)^2 ,
\end{equation}
which gives a more precise control on the correlations than the one
given by (\ref{e.correlalgo}) obtained in \cite{SchrammSteif} (which
implied a lower bound on the dimension of~%
$\Exc$ equal to $1/6$). As previously, proving such an estimate will be
achieved through a sharp understanding of the Fourier spectrum of
$f_R(\omega):= 1_{0\stackrel{\omega} {\longleftrightarrow} R}$ (but
this time, in order to get a sharp control, we will rely on the
``geometric approach'' explained
in Section~\ref{s.geo}). Indeed, recall from (\ref{e.CorrelQ}) that
%
%
\begin{equation}\label{e.CorrelSQb}
{\mathbb P}\bigl[{0\stackrel{\omega_0}{\longleftrightarrow}R,
0\stackrel{\omega_t}{\longleftrightarrow}R}\bigr] \lesssim\hat\Q
[{0< | \Spec_{f_R}|
< 1/t}] .
\end{equation}

Before considering the precise behavior of the lower tail of $\hat\Q
_{f_R}$, let us detail a heuristic argument
which explains why one should expect the above control on the
correlations (estimate \ref{e.correlgeom}).
Recall $\omega_t$ can be seen as a noised configuration of~$\omega_0$.
What we learned from Section \ref{s.geo} is that in the left--right case,
once the noise $\eps$ is high enough so that ``many pivotal points''
are touched,
the system starts being \textit{noise sensitive}. This means that the
left-to-right crossing event under consideration for $\omega^\eps$
starts being independent
of the whole configuration $\omega$. On the other hand (and this side
is much easier), if the noise is such that pivotal points are
(with high probability) not flipped, then the system is \textit{stable}.
Applied to our radial setting, if $R$ is some large radius and $r$ some
intermediate scale then, conditioned on the
event $\{ 0\stackrel{\omega}{\longleftrightarrow}R \}$, it is easy to
check that on average, there are about $r^{3/4}$ pivotal points at
scale $r$ [say, in the annulus $A(r,2r)$]. Hence,
if the level of noise $t$ is such that $t r^{3/4} \ll1$, then the
radial event is ``stable'' below scale~$r$. On the other hand if $t
r^{3/4} \gg1$, then many
pivotal points will be touched and the radial event should be sensitive
above scale $r$. As such, there is some limiting scale $L(t)$ separating
a frozen phase (where radial crossings are stable) from a sensitive
phase. The above argument gives $L(t) \approx t^{-4/3}$; this
limiting scale is often called the \textit{characteristic length}: below
this scale, one does not feel the effects of the dynamic from $\omega
_0$ to
$\omega_t$, while above this scale, connectivity properties start being
uncorrelated between $\omega_0$ and $\omega_t$.
\begin{remark}
Note that in our context of radial events, one has to be careful with
the notion of being \textit{noise sensitive} here. Indeed our sequence of
Boolean functions $(f_R)_{R>1}$ satisfy $\Var(f_R) \to0$; therefore,
the sequence trivially satisfies our initial definition of being noise
sensitive (Definition~\ref{d.NS}).
For such a sequence of Boolean functions corresponding to events of
probability going to zero, what we mean by being noise sensitive here is
rather that the correlations satisfy ${\mathbb E}[{f(\omega)
f(\omega^\eps)}]
\le
O(1) {\mathbb E}[{f}]^2$. This changes the intuition by
quite a lot: for example,
previously, for
a sequence of functions with $L^2$-norm one, analyzing the sensitivity
was equivalent to localizing where most of the Fourier mass was; while
if $\|f_R\|_2$ goes to zero, and if one wishes to understand for which
levels of noise $\eps_R>0$, one has ${\mathbb E}[{f_R(\omega)
f_R(\omega^{\eps_R})}] \le O(1) {\mathbb E}[{f_R}
]^2$, then it is the size of the lower tail of
$\hat f_R$ which is now relevant.
\end{remark}

To summarize and conclude our heuristical argument, as in the above
subsection, one can bound our correlations as follows:
\[
{\mathbb P}\bigl[{0\stackrel{\omega_0}{\longleftrightarrow}R,
0\stackrel{\omega_t}{\longleftrightarrow}R}\bigr] \le{\mathbb
P}\bigl[{0 \stackrel{\omega_0}{\longleftrightarrow} L(t) }\bigr]
{\mathbb P}\bigl[{L(t) \stackrel{\omega_0}{\longleftrightarrow}R,
L(t) \stackrel{\omega_t}{\longleftrightarrow}R}\bigr] .
\]
Not much is lost in this bound since connections are stable below the
correlation length $L(t)$. Now, above the correlation length, the
system starts being noise sensitive, hence one expects
\begin{eqnarray*}
{\mathbb P}\bigl[{L(t) \stackrel{\omega_0}{\longleftrightarrow}R,
L(t) \stackrel{\omega_t}{\longleftrightarrow}R}\bigr] &\le& O(1)
{\mathbb P}\bigl[{L(t) \stackrel{\omega}{\longleftrightarrow
}R}\bigr]^2 \\
&=& O(1) \alpha_1(L(t), R)^2 .
\end{eqnarray*}
Using quasi-multiplicativity and the values of the critical exponents,
one ends up with
\begin{eqnarray*}
{\mathbb P}\bigl[{0\stackrel{\omega_0}{\longleftrightarrow}R,
0\stackrel{\omega_t}{\longleftrightarrow}R}\bigr] & \le& O(1)
\alpha_1 (L(t))^{-1} \alpha
_1(R)^2 \\
& \lesssim& t^{-5/36} \alpha_1(R)^2
\end{eqnarray*}
as desired.

We wish to briefly explain how to adapt the geometric
approach of Section~\ref{s.geo} to our present radial setting.
Recall that the goal is to prove an analog of Theorem~\ref
{th.GPSprecise} for the radial crossing event $f_R(\omega) =
1_{0\stackrel{\omega}{\longleftrightarrow}R}$.\vspace*{1pt}
Before stating such an analogous statement (which would seem at first
artificial), 
we need to understand some properties of $\Spec_{f_R}$. 

First of all, similarly to Section \ref{s.geo}, it is easy to check
that $\hat\PP[{|\Spec_{f_R}|}] \asymp R^2 \alpha_4(R)$
(recall that $\hat\PP_{f_R}$ is the renormalized spectral measure
$\frac1 {\| f_R\|^2} \hat\Q_{f_R} = \frac1 {\alpha_1(R)} \hat\Q
_{f_R}$). The second moment
is also easy to compute. One can conclude from these estimates that
with positive probability (under $\hat\PP_{f_R}$), $|\Spec_{f_R}|$ is
of order $R^2 \alpha_4(R)$ (\mbox{$\approx$}$R^{3/4}$ on~$\T$).

Now, following the strategy explained in Section \ref{s.geo}, we divide
the ball of radius $R$ into a grid of mesoscopic squares of radii $r$
(with $1< r < R$).
It is easy to check that for any such $r$-square $Q$, one has
\[
{\hat\PP}_{f_R} [ |\Spec\cap Q| \mid \Spec\cap Q \neq
\varnothing]
\asymp
r^2 \alpha_4(r) .\vadjust{\goodbreak}
\]
Furthermore, and this part is very similar as in the chordal case, one
can obtain a ``weak'' control on the dependencies within $\Spec_{f_R}$;
namely\vadjust{\goodbreak} that for any
$r$-square $Q$ and any subset $W$ such that $Q\cap W =\varnothing$, one has
\[
{\hat\PP}_{f_R} [ |\Spec\cap\bar Q| \geq r^2 \alpha_4(r)
\mid \Spec\cap Q \neq\varnothing, \Spec\cap W = \varnothing
] \geq c ,
\]
where $\bar Q$ is the concentric square in $Q$ of side-length $r/2$ and
$c>0$ is some absolute constant.

As in Section \ref{s.geo}, if $S\subset[-R,R]^2$, we denote by
$S_{(r)}$ the set of those $r\times r$ squares which intersect $S$.
Thanks to the above weak independence statement, it remains (see
Section \ref{s.geo}) to study the lower tail of the number of
mesoscopic squares
touched by the spectral sample, that is, the lower tail of $|\Spec
_{(r)}|$. (Only the very bottom-part of this distribution needs to be
understood in a sharp way.) Hence,
one has to understand how $\Spec_{(r)}$ typically looks when it is
conditioned to be of very small size [say of size less than $\log R/r$,
i.e., much smaller than the average size of $\Spec_{(r)}$ which is of
order $(R/r)^{3/4}$]. The intuition differs here from the chordal case
(analyzed in Section \ref{s.geo}); indeed recall that in the chordal
case, if $|\Spec_{(r)}|$ was conditioned to be of small size, then it
was typically ``concentrated in space'' somewhere ``randomly'' in the
square. In our radial setting, if $|\Spec_{(r)}|$ is conditioned to be
of small size, it also tends to be localized, but instead of being
localized somewhere randomly through the disk of radius $R$, it is
localized around the origin. The reason of this localization around the
origin comes from the following estimate:
if $Q_0$ denotes the $r$-square surrounding the origin
\begin{eqnarray*}
\hat\PP\bigl[{Q_0 \notin\Spec_{(r)}}\bigr] & \le& \hat\PP
\biggl[{\Spec_{f_R} \subset B\biggl(0,\frac r 2\biggr)^c}\biggr]\\
& = & \frac{1} {\alpha_1(R)} \hat\Q\biggl[{\Spec_{f_R} \subset
B\biggl(0,\frac r 2\biggr)^c}\biggr] \\
& = & \frac1 {\alpha_1(R)} {\mathbb E}\bigl[{{\mathbb E}[{f_R
\mid B(0,R)\setminus B(0,r)}]^2}\bigr]\qquad \mbox{(by Lemma
\ref{l.A})} \\
& \le& \frac{O(1)} {\alpha_1(R)} \alpha_1(r,R) \alpha_1(r)^2
\le
O(1) \alpha_1(r).
\end{eqnarray*}
(In the last line, we used the fact that ${\mathbb E}[{f_R
\mid B(0,R)\setminus B(0,r)}] \le\alpha_1(r) 1_{r
\leftrightarrow R}$.)
Since $\alpha_1(r)$ goes to 0 as $r\to\infty$, this means that for a
mesoscopic scale $r$ such that $1\ll r \ll R$,
$\Spec_{(r)}$ ``does not like'' to avoid the origin. In some sense,
this shows that the origin is an attractive point for $\Spec_{f_R}$.

Still, this does not yet imply that the origin will be attractive also
for~$\Spec_{(r)}$ conditioned to be of small size. If one assumes that,
as in the chordal case,~$\Spec_{(r)}$ tends to be localized when it is
conditioned to be of small size (i.e., that for $k$ small, $\hat\PP
[{|\Spec_{(r)}| =k}] \approx\hat\PP[{|\Spec
_{(r)}| = 1}]$), then it not hard to see
that it has to be localized around the origin: indeed, one can estimate
by hand that once conditioned on $|\Spec_{(r)}|=1$, $\Spec_{(r)}$ will
be (with high conditional probability) close to the origin. To
summarize, assuming that $\Spec_{(r)}$ tends to ``clusterize'' when it
is conditioned to be of small size, one has for $k$ ``small''
\begin{eqnarray*}
\hat\Q\bigl[{\bigl|\Spec_{(r)}\bigr| = k}\bigr] & \approx& \hat\Q
\bigl[{\bigl|\Spec_{(r)}\bigr| =1}\bigr]
\asymp \hat\Q\bigl[{\Spec_{(r)} = Q_0}\bigr] \\
& \asymp& \hat\Q[{\varnothing\neq\Spec_{f_R} \subset
B(0,r)}]
\\
&\asymp& {\mathbb E}[{{\mathbb E}[{f_R \mid
B(0,r)}]^2}] - \hat\Q[{\Spec
_{f_R} = \varnothing}]
\\
&\asymp& O(1) \alpha_1(r) \alpha_1(r,R)^2 - \alpha_1(R)^2
\\
&\asymp& \frac1 {\alpha_1(r)} \alpha_1(R)^2 \qquad\mbox{(using
quasi-multiplicativity)} .
\end{eqnarray*}

Using an inductive proof
(similar to the one used in the chordal case, but where the origin
plays a special role), one can prove that it is indeed the case that
$\Spec_{(r)}$
tends to be localized when it is conditioned to be of small size. More
quantitatively, the inductive proof leads to the following estimate on the
lower tail of $|\Spec_{(r)}|$:
\begin{proposition}[(\cite{GPS}, Section 4)]\label{p.roughRadial}
There is a sub-exponentially fast growing function $g(k), k\geq1$,
such that for any
$1\leq r \leq n $
%
%
\begin{equation}\label{e.roughS}
{\hat\Q}_{f_R} \bigl[ \bigl| \Spec_{(r)} \bigr| = k \bigr] \leq g(k)
\frac1 {\alpha
_1(r)} \alpha_1(R)^2 .
\end{equation}
\end{proposition}

Using the same technology as in Section \ref{s.geo} (i.e., the weak
independence control, the above proposition on the lower tail of the
mesoscopic sizes
and a~large deviation lemma which helps implementing the ``scanning
procedu\-re''), one obtains that $\hat\Q[{0<|\Spec_{f_R}|< r^2
\alpha_4(r)}]
\asymp{\hat\Q}_{f_R} [ |\Spec_{(r)}|=1 ]
\asymp\frac{1} {\alpha_1(r)} \alpha_1(R)^2$, which is exactly a sharp
lower-tail estimate on $\hat f_R$. More exactly one has the following theorem
(analog of Theorem \ref{th.GPSprecise} in the chordal case):
\begin{theorem}[\cite{GPS}]
Both on the triangular lattice $\T$ and on $\Z^2$, if $f_R$ denotes the
radial event up to radius $R$, one has
\[
\hat\Q[{| \Spec_{f_R} | < r^2 \alpha_4(r)}] \asymp
\frac1 {\alpha_1(r)}
\alpha_1(R)^2 ,
\]
where the constants involved in $\asymp$ are absolute constants.
\end{theorem}

On the triangular lattice $\T$, using the knowledge of the critical
exponents, this can be read as follows:
\[
\hat\Q[{| \Spec_{f_R}| < r^{3/4}}] \approx r^{5/48}
\alpha_1(R)^2
\]
or equivalently for any $ u \lesssim R^{3/4}$,
\[
\hat\Q[{| \Spec_{f_R}| < u}] \approx u^{5/36} \alpha
_1(R)^2 .
\]
Recall our goal was to obtain a sharp bound on the correlation between~$f_R(\omega_0)$ and $f_R(\omega_t)$.
From our control on the lower tail of $\hat f_R$ and equation~(\ref
{e.CorrelSQb}) one obtains
\begin{eqnarray*}
{\mathbb P}\bigl[{0\stackrel{\omega_0}{\longleftrightarrow}R,
0\stackrel{\omega_t}{\longleftrightarrow}R}\bigr] &\lesssim& \hat
\Q[{0< | \Spec
_{f_R}| < 1/t}]
\\
& \lesssim& t^{-5/{36}} \alpha_1(R)^2 ,
\end{eqnarray*}
which, as desired, implies that the exceptional set $\Exc$ is indeed
a.s. of dimension $\frac{31}{36}$ (as claimed above, the upper bound
is much easier
using the knowledge on the density function $\theta_{\T}$).

\makeatletter\write@toc@ignorecontentsline\makeatother
\section*{\texorpdfstring{Acknowledgments.}{Acknowledgments}}
\makeatletter\write@toc@restorecontentsline\makeatother

I wish to warmly thank Jeff Steif who helped tremendously in improving
the readability of this survey.
The week he invited me to spend in G\"oteborg was a great source of
inspiration for writing this paper.
I would also like to thank G\'abor Pete: since we started our
collaboration three years ago with Oded, it has always been a pleasure
to work with.
Finally, I would like to thank Itai Benjamini who kindly invited
me to write the present survey.

%

%
\printaddresses

\end{document}